\documentclass[11pt]{article}

\usepackage{latexsym}
\usepackage{amsfonts}
\usepackage{amsmath}

\usepackage[dvips]{color}

\newtheorem{thrm}{Theorem}[section]
\newtheorem{lemma}[thrm]{Lemma}
\newtheorem{prop}[thrm]{Proposition}

\newtheorem{cor}[thrm]{Corollary}
\newtheorem{remark}[thrm]{Remark}

\numberwithin{equation}{section}

\usepackage{enumerate}

\usepackage{amssymb}
\usepackage{mathtools}
\mathtoolsset{showonlyrefs}
\usepackage{mathrsfs}
\usepackage{comment}
\usepackage{hyperref}
\usepackage{xcolor}

\def\P{\mathbb{P} }

\def\R{\mathbb{R} }
\def\N{\mathbb{N} }

\makeatletter
\begin{document}
\allowdisplaybreaks

\title{\Large \bf{Invariance principle for the maximal
position process of
branching Brownian motion in random environment}\footnote{The research of this project is supported by the National Key R\&D Program of China (No. 2020YFA0712900).}}
\author{ \bf  Haojie Hou \hspace{1mm}\hspace{1mm}
Yan-Xia Ren\footnote{The research of this author is supported by NSFC (Grant Nos. 11731009 and  12071011) and LMEQF.\hspace{1mm} } \hspace{1mm}\hspace{1mm} and \hspace{1mm}\hspace{1mm}
Renming Song\thanks{Research supported in part by a grant from the Simons
Foundation (\#429343, Renming Song).}
\hspace{1mm} }

\date{}
\maketitle

\begin{abstract}
In this paper we study the maximal position process of  branching Brownian motion in random spatial environment.  The random
environment is given by a  process $\xi = \left(\xi(x)\right)_{x\in\mathbb{R}}$ satisfying certain conditions.  We show that
the maximum position $M_t$ of particles alive at time $t$ satisfies a quenched strong law of large numbers and an annealed invariance principle.
\end{abstract}

\medskip

\noindent\textbf{AMS 2010 Mathematics Subject Classification:} 60J80; 60G70; 82B44.

\medskip

\noindent\textbf{Keywords and Phrases}: Branching Brownian motion, random environment, law of large numbers,  invariance principle.

\section{Introduction}

\subsection{Background and assumptions}
A binary branching Brownian motion (BBM) in $\R$ can be described as follows:
There is an initial particle starting from the origin, and the particle moves as a standard Brownian motion. After an independent exponential amount of time with parameter one,  the particle dies and splits into two new particles.  The offspring evolve independently as their parent.

Let $M_t$ denote the maximal position among all the particles alive at time $t$.
Kolmogorov,  Petrovskii and Piskounov\cite{KPP} showed that $\frac{M_t}{t}$ converges to $\sqrt{2}$ in probability as $t\to\infty$.
In \cite{Bramson1, Bramson2}, Bramson proved that, as $t\to\infty$,
$M_t-(\sqrt{2}t - \frac{3}{2\sqrt{2}}\ln t)$ converges weakly to a limit related to a travelling wave solution.  In \cite{LS}, Lalley and Sellke gave a probabilistic representation of the limit. Thus in this classical case, $M_t = \sqrt{2} t + O(\ln t)$ as $t\to\infty$.

\v{C}ern\'{y} and Drewitz \cite{JD} studied continuous-time branching random walks in random environment and
their results show that the asymptotic behavior of branching random walk in random environment is different from the classical case.
A continuous-time branching random walk in random environment can be described as follows:
Suppose that $\xi=\left\{\xi(x): x\in\mathbb{Z}\right\}$ is a family of iid random variables on a probability space $(\Omega,\mathcal{F},\mathbb{P})$ with
$0 < \textup{essinf}\  \xi(0) < \textup{esssup}\  \xi(0) < \infty$.
Let $\Omega_0:=\{\omega: \left\{\xi(x,\omega), x\in\mathbb{Z}\right\} \text{ is a sequence of positive numbers}\}$. Then $\mathbb{P}(\Omega_0)=1$.
For any $\omega\in \Omega_0$ and $x\in \mathbb{Z}$,
let $\mathrm{P}_x^\xi$ denote the law of a
continuous-time binary branching random walk starting from $x$
with branching rate $\xi(\cdot,\omega)$.
There is a particle at $x$ at time 0. As time evolves, the particle moves as a continuous time random walk with jump rate 1.  In addition, and independently of everything else,
while at site $y$, a particle splits into two at rate $\xi(y, \omega)$,
and when it does so, the two new particles evolve independently according to the same diffusion mechanism as their parent. $\mathrm{P}_x^\xi$  is referred to as the quenched law and $\mathbb{P}\times \mathrm{P}_x^\xi$ is referred to as the annealed law.

We still use $M_t$ to denote the maximal displacement at time $t$ of branching random walk in random environment. Under some conditions, \v{C}ern\'{y} and Drewitz \cite{JD} proved an annealed invariance principle for $M_t$.
That is, there exist $v_0>0$ and $\sigma_0 > 0$
such that, under $\mathbb{P}\times \mathrm{P}_{0}^\xi$,  the sequence of processes
$$[0, \infty) \ni t \mapsto \frac{M_{nt}-v_{0} n t}{\sigma_{0} \sqrt{n}}, \quad n \in \mathbb{N},
$$
converges, as $n\to\infty$, weakly
to a standard Brownian motion.
Thus $\mathbb{P}$-almost surely, $M_t \neq v_0t + O(\ln t)$ as $t\to\infty$ under $\mathrm{P}_0^\xi$.

In this paper, we study branching Brownian motion in random environment (BBMRE). As in \cite{JD},  we suppose that $\left\{\xi(x), x\in\mathbb{R}\right\}$ is a family of non-negative random variables on
$(\Omega,\mathcal{F},\mathbb{P})$.
We will add conditions on $\left\{\xi(x), x\in\mathbb{R}\right\}$ below.
For any fixed $\omega\in\Omega$, $\left\{\xi(x,\omega),x\in\mathbb{R}\right\}$ can be regarded as a branching rate.
 Branching Brownian motion in random environment $\xi$  started at $x\in\R$ can be defined as follows: Conditionally
on a realization of $\xi$, we place one particle at $x$ at time 0. As time evolves, all particles
move independently according to a standard Brownian motion. In addition, and independently of
everything else, while at $y$, a particle splits at rate $\xi(y)$. Once a particle splits, this particle is
removed and, randomly and independently from everything
else, replaced by $k$ new particles,  with probability $p_k$,
at the position of the removed particle.
These $k$ new particles evolve
independently according to the same diffusion-branching mechanism
as their parent.
  For a given $\xi(\cdot,\omega)$, let $\mathrm{P}_x^\xi$ denote the law of branching Brownian motion starting from $x$ with spatial-dependent branching rate $\xi(\cdot,\omega)$ and offspring distribution $\left\{p_k: k = 0,1,2,...\right\}$. More precisely, let $N(t)$ be the set of particles alive at time $t$ and $X_t^\nu$ be the position of the particle $\nu\in N(t).$  Define
   $$X_t:= \sum_{\nu\in N(t)} \delta_{X_t^\nu},\quad t\geq 0,$$ which is  the point process generated by the position of the particles alive at time $t$.
   We call $\{X_t, t\geq 0\}$ a branching Brownian motion in random in random environment $\xi$.

For any $\nu\in N(t)$ and $s\in [0, t]$,  we use $X_s^\nu$ to denote the position of $\nu$ or its ancestor at time $s$.
$X_{\cdot}^\nu$ is referred to as the genealogy of $\nu\in N(t)$.

   For any $t\geq 0$, put
   $$M_t:= \max_{\nu\in N(t)} X_t^\nu.$$
The purpose of this paper is to study the limit behaviour of $M_t$ as $t\to\infty$ under the  quenched law $\mathrm{P}_x^\xi$, $x\in\R$, and prove an invariance principle for $M_t$ under the annealed law $\P\times \mathrm{P}_0^\xi$.
Throughout this paper, we assume that $p_0 = 0$ and $\sum_{k=1}^\infty kp_k = m > 1$.
When $p_0>0$,
we can consider $M_t$ under $\mathrm{P}_x^\xi(\cdot | \mbox{ survival})$ and $\P \times \mathrm{P}_x^\xi(\cdot | \mbox{ survival})$.

The first  assumption {\bf(H1)} on  the environment $\xi$ already appeared in \cite{Nolen, Nadin, BN, CDS, DS}.
\begin{itemize}
	\item [{\bf(H1)}]  {\bf(a)}
$\mathbb{P}$-almost surely, $\xi$ is uniformly H\"{o}lder continuous, i.e., there exist constants $C(\xi), \alpha(\xi) > 0$ such that
	$$|\xi(x,\omega)-\xi(y,\omega)| \leq C(\xi)|x-y|^{\alpha(\xi)},\quad \forall \ x,y\in\mathbb{R}.$$

	{\bf(b)} There exist constants $0 < \textup{ei} \leq \textup{es} < \infty$ such that $\mathbb{P}$-a.s.,
	$$\textup{ei} \leq \xi(x) \leq \textup{es}, \quad \forall\ x\in\mathbb{R}.$$

	{\bf(c)} There exists a group of measure-preserving transformations $\{\theta_y \}_{y\in\mathbb{R}}$, acting ergodically on $(\Omega, \mathcal{F},\mathbb{P})$, such that, for almost every $\omega \in\Omega$,  $\xi(x+y,\omega) = \xi(x,\theta_y\omega)$ holds for all $x,y\in\mathbb{R}.$

	{\bf(d)}
$\xi$ satisfies the $\zeta$-mixing condition:
Let $\mathcal{F}_x:= \sigma\left(\xi(t): t\leq x \right)$,  $\mathcal{F}^y : = \sigma\left(\xi(r): r\geq y \right)$ and $\zeta:[0,\infty) \to [0, \infty)$ be a continuous non-increasing
function with $\sum_{k=0}^{\infty} \zeta(k)<\infty$.
It holds
that for all $x,y\in \mathbb{R}$ with $x\leq y $ and all $X\in L^1(\Omega, \mathcal{F}_x,\mathbb{P})$ and $Y\in L^1(\Omega, \mathcal{F}^y,\mathbb{P}),$
	\begin{equation}\label{Mixing}
		\begin{array}{l}
					\displaystyle\left|\mathbb{E}\left(\left(X-\mathbb{E}(X)\right) \mid \mathcal{F}^{y}\right)\right| \leq \mathbb{E}\left(|X|\right) \cdot \zeta(y-x), \\
			\displaystyle\left|\mathbb{E}\left(\left(Y-\mathbb{E}(Y)\right) \mid \mathcal{F}_{x}\right)\right| \leq \mathbb{E}
\left(|Y|\right) \cdot \zeta(y-x).
		\end{array}
	\end{equation}
\end{itemize}
The mixing condition above is given in \cite{DS}, and is  weaker than the one given in \cite{Nolen}.
For simplicity, we assume that {\bf(c)} holds for all $\omega \in \Omega$ rather than $\P$-almost surely.

Our second assumption is \cite[(VEL)]{DS}, but we will describe it using principal eigenvalues.
Define
$$\mathcal{A}_\infty : = \left\{\phi\in \mathcal{C}^2(\mathbb{R}), \frac{\phi_x}{\phi} \in L^\infty(\mathbb{R}), \phi > 0\ \textup{in}\ \mathbb{R}, \lim_{|x|\to +\infty} \frac{1}{x} \ln \phi(x) = 0 \right\}.$$
For $\lambda \in \mathbb{R}$ and $\phi(x) \in \mathcal{C}^2(\mathbb{R}),$ define
$$
L_\lambda^\omega \phi := \frac{1}{2} \phi_{xx} - \lambda\phi_x +\left(\frac{\lambda^2}{2} + (m-1)\xi(x,\omega) \right)\phi.
$$
For any $\lambda\in\mathbb{R}$ and non-empty open interval $I \subset \mathbb{R}$,
define
\begin{align*}
	\underline{\gamma}(L_\lambda^\omega)& : =\sup\left\{\gamma \big|\  \exists\  \phi\in \mathcal{A}_\infty \ \textup{such that}\ L_\lambda^\omega \phi \geq \gamma\phi \ \textup{in}\ \mathbb{R} \right\},\\
	\overline{\gamma}(L_\lambda^\omega)& : = \inf\left\{\gamma \big|\  \exists\  \phi\in \mathcal{A}_\infty \ \textup{such that}\ L_\lambda^\omega \phi \leq \gamma\phi \ \textup{in}\ \mathbb{R} \right\},\\
	\Lambda_1(L_0^\omega, I)
	 &:= \inf\left\{\gamma \big| \ \exists \phi\in \mathcal{C}^2(I)\cap \mathcal{C}^0(\overline{I}),\ \phi > 0\ \textup{in}\ I,\ \phi=0 \ \textup{in}\ \partial I,\ L_0^\omega \phi \leq \gamma \phi \ \textup{in}\ I \right\},\\
	\Lambda_1(L_0^\omega)& := \lim_{R \to +\infty} \Lambda_1\left(L_0^\omega, (-R, +R) \right).
\end{align*}
The definitions of $\underline{\gamma}(L_\lambda^\omega)$ and $ \overline{\gamma}(L_\lambda^\omega)$
are special cases of Berestycki and Nadin \cite[(10) and (11)]{BN} with $R = -\infty$,
and the definitions of $\Lambda_1(L_0^\omega, I)$ and $\Lambda_1(L_0^\omega)$
 are given in \cite[(45) and (47)]{BN}. By \cite[Theorem 5.1 (1)]{BN}, there exists
 $\Omega_1 \subset \Omega$ with $\P(\Omega_1 )=1$ such that, for any $\omega \in \Omega_1$
 and  $\lambda \in \R$, $\underline{\gamma}(L_\lambda^\omega) = \overline{\gamma}(L_\lambda^\omega)$
 and this common value does not depend on $\omega \in \Omega_1$.
 We will write $\underline{\gamma}(L_\lambda^\omega)$ as $\gamma(\lambda)$ for simplicity.  By \cite[(51) and Lemma 5.1]{BN}, we know that there exist two constants
 $\rho_L\leq 0$ and $\rho_R\geq 0$
 such that when $\lambda \in [\rho_L, \rho_R]$, $\gamma(\lambda) = \Lambda_1(L_0^\omega)$.
  Nadin \cite[Lemma 3.2]{Nadin} proved that
  $$\rho_R = -\rho_L =:\rho \geq 0.$$
From the argument of \cite[Theorem 5.1]{BN}, we know that
when $\lambda \notin [\rho_L, \rho_R]$,  $\gamma(\lambda) > \Lambda_1(L_0^\omega)$,
and in  this case, by \cite[Theorem 5.1 (2)]{BN},  there exists $\phi = \phi(\cdot, \lambda, \omega) \in \mathcal{A}_\infty$ satisfying
\begin{equation}\label{Principal-Engenfunction}
	L_\lambda^\omega \phi = \frac{1}{2} \phi_{xx} - \lambda\phi_x +\left(\frac{\lambda^2}{2} + (m-1)\xi(x,\omega) \right)\phi= \gamma(\lambda) \phi.
\end{equation}
Moreover, the proof of  \cite[Theorem 5.1 (1)]{BN} shows that $\gamma(\cdot)$ is positive and strictly increasing in
$(\rho, +\infty)$.
Also, Nadin \cite[Lemma 3.2]{Nadin} shows that $\gamma(\cdot)$ is even.
By  \cite[Lemma 3.3]{Nadin}, $\phi$ is unique up to
multiplicative constant.
Thus, we can assume that $\phi(0,\lambda) = 1$ for all $|\lambda|> \rho$ and
$\omega \in \Omega_1$
without loss of generality. Also \cite[Proposition 1 and Proposition 3]{Nadin} imply that for all $\lambda \in \mathbb{R},$
$$(m-1)\textup{es} + \lambda^2/2 \geq \gamma(\lambda) \geq (m-1)\textup{ei} + \lambda^2/2,$$
so $v^*$ and $\lambda^*$ below are well-defined:
\begin{equation}\label{Critical-value}
	v^*:= \min_{\lambda > 0}\frac{\gamma(\lambda)}{\lambda},\quad \lambda^*:= \mathop{\textup{argmin}}\limits_{\lambda > 0} \frac{\gamma(\lambda)}{\lambda}.
\end{equation}
Our second assumption is as follows:
\begin{itemize}
	\item [{\bf(H2)}] $\gamma(\lambda^*) >  (m-1)\textup{es}$.
\end{itemize}
Since $\gamma(\rho)=\gamma(0) \leq (m-1)\textup{es}$, {\bf(H2)} implies $\lambda^* > \rho$.
We mention in passing
that, by \cite[Proposition 4.10]{DS}, there exists a $\xi$ satisfying {\bf(H1)} but not {\bf(H2)}.

Our final assumption {\bf(H3)} will be used
when we deal with the invariance principle for $M_t$.
\begin{itemize}
	\item [{\bf(H3)}]
		$m_2:= \sum_{k=1}^\infty k^2 p_k < +\infty.$
\end{itemize}

BBMRE is related to  the random F-KPP:
\begin{equation}\label{KPP}
	w_t= \frac{1}{2}w_{xx} + \xi(x,\omega)\left(1-w-\sum_{k=1}^\infty p_k (1-w)^k\right).
\end{equation}
 Solutions to \eqref{KPP} can be written as
\begin{equation}\label{Prob-representation}
	w(t,x)= 1 - \mathrm{E}_x^\xi \left(\prod_{\nu\in N(t)}\left(1- w(0,X_t^\nu)\right)\right),
\end{equation}
see  \cite[Proposition 2.1]{DS}. \eqref{Prob-representation} is
referred to as  McKean's representation, see \cite{McKean} for  the case of homogeneous branching Brownian motion. In particular, $w(t,x)=\mathrm{P}^\xi_x(M_t\geq 0)$ solves \eqref{KPP} with $w(0,x)=1_{[0,\infty)}(x)$, and $w(t,x) = \mathrm{P}_x^\xi\left(\min_{\nu\in N(t)} X_t^\nu \leq 0\right)$ solves \eqref{KPP} with $w(0,x) = 1_{(-\infty, 0] }(x)$.

Freidlin \cite[Chapter 7]{FM} studied branching Brownian motion in random environment by considering
\eqref{KPP}.
The asymptotic wave front propagation velocity for
\eqref{KPP} is of special interest. By \cite[Theorem 7.6.1]{FM}, under
suitable assumptions, the solution $w(t, x)$ to
\eqref{KPP} converges
to 0 (resp. 1), uniformly for all $x\geq vt$ with $v >v^*$ (resp. for all $x\leq vt$ with $v\in(0, v^*)$), as $t\to\infty$.
 In particular, using \cite[Theorem 7.6.1]{FM} with  $w(t,x) = \mathrm{P}_x^\xi\left(\min_{\nu\in N(t)} X_t^\nu \leq 0\right)$ (which implies $w(0,x) = 1_{(-\infty, 0] }(x)$), we have that $\P$-a.s. $\lim_{t\to\infty}\frac{m(t)}{t}=v^*$, where
$$
m(t):= \sup\left\{x: \mathrm{P}_x^\xi\left(\min_{\nu\in N(t)} X_t^\nu \leq 0\right) = \frac{1}{2}\right\}
$$
is known as the front of \eqref{KPP}.
In particular when $\xi(x, \omega)\equiv c > 0$ with $c$ being a positive constant,  $v^* = \sqrt{2c(m-1)}$.
This result is due to Kolmogorov-Petrovskii-Piskunov \cite{KPP}.

Nolen \cite {Nolen} proved a
central limit theorem for $\mathcal{X}_t:= \sup \{x: w(t,x) = \frac{1}{2}\}$, where the assumption
on the initial value $w(0,x)$  of \eqref{KPP}
falls into the supercritical regime --- the limit of $\frac{\mathcal{X}_t}{t}$ is larger than the minimal speed $v^*$.
Recently, Drewitz and Schmitz \cite{DS} studied the case when $w(0,x)$ satisfies
\begin{equation}\label{KPP-INI}
	0\leq w(0,x)\leq 1_{(-\infty,0]}(x),\quad \mbox{ and } \int_{[x-N,x]}w(0,y) dy\geq \delta,\quad  \forall x\leq -N',
\end{equation}
for some fixed $\delta, N, N'>0$.
This case corresponds to the critical regime.
\cite{DS} proved an  invariance principle
 for  $m^\epsilon(t):= \sup\left\{x\in\R: \mathrm{P}_x^\xi\left(\min_{\nu\in N(t)} X_t^\nu \leq 0\right) \geq \varepsilon\right\}$ with $\varepsilon \in (0,1)$.
 Note that
  $$\mathrm{P}_x^\xi\left(\min_{\nu\in N(t)} X_t^\nu \leq 0\right) =\mathrm{P}_x^\xi\left(\max_{\nu\in N(t)} (-X_t^\nu) \geq 0\right),$$
  and that
 $\{\sum_{\nu\in N(t)} \delta_{-X_t^\nu}, t\geq 0\}$ is a BBMRE starting at $-x$ with branching rate $\tilde\xi(x)=\xi(-x)$.
   Therefore,
  $$m^\epsilon(t)=\sup\{x\in\R: \mathrm{P}^{\tilde\xi}_{-x}(M_t\geq 0)\geq\epsilon\}.$$
  By the translation and reflection invariance of the environment $\xi$, we have for $x\in\R$,
   $$\mathrm{P}^{\xi}_{-x}(M_t\geq 0)=\mathrm{P}^{\xi}_{0}(M_t\geq x) \quad\mbox{in distribution under }\P.$$
This says  that  the invariance principle for $m^\epsilon(t)$ in \cite{DS} is related to our invariance principle for $M_t$. In  the homogeneous medium, i.e., $\xi(x,\omega)\equiv c$,
$\mathrm{P}_x^c\left(\min_{\nu\in N(t)} X_t^\nu \leq 0\right) =\mathrm{P}_0^c\left(M_t\geq x\right),$  and then
 $m^\epsilon(t)$, the front of the solution to \eqref{KPP} coincides with the  $\epsilon$ median
 of the distribution of the maximal particle of the BBM. In the random medium case, $m^\epsilon(t)$ has the same distribution as the $\epsilon$ median
of the  distribution of the maximal particle of BBMRE.
In this paper, we are  mainly interested in the behavior of $M_t, t\geq 0$.

 According to \cite[Proposition 2.3]{DS}, we have the following many-to-one and many-to-two formulae:
\begin{prop}\label{Many-to-one}
	 Let $\varphi_{1}, \varphi_{2}:[0, \infty) \rightarrow[-\infty, \infty]$ be c\`{a}dl\`{a}g functions with $\varphi_{1} \leq \varphi_{2}$. Then the first and second moments of the number of particles in $N(t)$ with genealogy staying between $\varphi_{1}$ and $\varphi_{2}$ in the time interval $[0, t]$ are given by
	$$
	\begin{array}{l}
		\displaystyle\mathrm{E}_{x}^{\xi}\left(\#\left\{\nu \in N(t): \varphi_{1}(s) \leq X_{s}^{\nu} \leq \varphi_{2}(s),\  \forall s \in[0, t]\right\}\right) \\
		\displaystyle\quad=\Pi_{x}\left(\exp \left\{\int_{0}^{t} (m-1)\xi\left(B_{r}\right) \mathrm{d} r\right\} ; \varphi_{1}(s) \leq B_{s} \leq \varphi_{2}(s),\  \forall s \in[0, t]\right)
	\end{array}
	$$	
and
$$
	\begin{array}{l}
	\displaystyle	\mathrm{E}_{x}^{\xi}\left(\#\left\{\nu \in N(t): \varphi_{1}(s) \leq X_{s}^{\nu} \leq \varphi_{2}(s)\  \forall s \in[0, t]\right\} ^{2}\right)
\\ \displaystyle =\Pi_{x}\left(\exp \left\{\int_{0}^{t}(m-1) \xi\left(B_{r}\right) \mathrm{d} r\right\} ; \varphi_{1}(s) \leq B_{s} \leq \varphi_{2}(s)\  \forall s \in[0, t]\right)
 \\ \displaystyle +\left(m_{2}-m\right) \int_{0}^{t} \Pi_{x}\Big(\exp \left\{\int_{0}^{s} (m-1)\xi\left(B_{r}\right) \mathrm{d} r\right\} \xi\left(B_{s}\right) {1}_{ \{\varphi_{1}(r) \leq B_{r} \leq \varphi_{2}(r),\  \forall 0 \leq r \leq s\}}
\\ \displaystyle	\quad \times\left(\Pi_{y}\left(\exp \left\{\int_{0}^{t-s}(m-1) \xi\left(B_{r}\right) \mathrm{d} r\right\} {1}_{\{\varphi_{1}(r+s) \leq B_{r} \leq \varphi_{2}(r+s)\  \forall 0 \leq r \leq t-s\}}\right)\right)_{\mid y=B_{s}}^{2}\Big) \mathrm{d} s
	\end{array}
$$
	respectively.
\end{prop}

For $a \in \mathbb{R}$, letting $\varphi_{1}\equiv-\infty$, $\varphi_{2}(s)=\left\{\begin{array}{ll}+\infty, &\mbox{ if }s < t,\\
a,  &\mbox{ if }s=t\end{array}\right.$ in Proposition \ref{Many-to-one}, we get that
\begin{align}\label{Many-to-one2}
	\mathrm{E}_x^\xi \left(\sum_{\nu\in N(t)} f\left(X_t^\nu\right)\right)= \Pi_x \left(\exp\left\{\int_0^t (m-1)\xi(B_r)\mathrm{d}r \right\} f(B_t)\right),
\end{align}
first for $f(x)=1_{(-\infty, a]}(x)$, and then  for any non-negative Borel function $f$.

Recall that for $|\lambda| > \rho$ and $\omega \in \Omega_1$,
$\phi(x, \lambda,\omega)$ solves  \eqref{Principal-Engenfunction} with $\phi(0,\lambda,\omega)=1$.
Define
\begin{equation}\label{def-psi}
\psi(x,\lambda) := e^{-\lambda x}\phi(x,\lambda),\quad x\in\R
\end{equation}
 and
 $$ u_{\psi(\lambda)}(t,x) := e^{\gamma(\lambda)t} \psi(x,\lambda) = e^{\gamma(\lambda)t}e^{-\lambda x}\phi(x,\lambda),\quad x\in \R.$$
Then $\psi(\cdot,\lambda)$ solves the problem
\begin{equation}\label{Principal-Engenfunction-2}
	\left\{\begin{array}{rl}&\displaystyle\frac{1}{2}\psi_{xx} + (m-1)\xi(x,\omega) \psi = \gamma(\lambda) \psi,\quad x\in\R,\\
&\displaystyle\psi(0,\lambda,\omega)=1,
\end{array}\right.
\end{equation}
and  $u_{\psi(\lambda)}$ solves the problem
\begin{equation}\label{Linear-Equa}
\left\{	\begin{array}{rl}&\displaystyle u_t = \frac{1}{2}u_{x x} + (m-1)\xi(x,\omega) u, \quad t>0, \quad x\in\mathbb{R},\\
&u(0,x,\lambda) = \psi(x,\lambda).
\end{array}\right.
\end{equation}

We will always use $\left\{(B_t)_{t\geq 0}; \Pi_x\right\}$ to denote a standard Brownian motion starting from $x$ at $t=0$, and  also use  $\Pi_x$ for  expectation  with respect to $\Pi_x$ for simplicity. By the Feynman-Kac formula, $u_{\psi(\lambda)}(t,x)$ has the following representation
\begin{equation}\label{F-K-u}
	u_{\psi(\lambda)}(t,x) = \Pi_x\left(\exp\left\{ \int_0^t(m-1) \xi(B_s,\omega )\mathrm{d}s \right\}e^{-\lambda B_t}\phi(B_t,\lambda)\right).
\end{equation}
Thus
\begin{align*}
		& e^{-\lambda x}\phi(x,\lambda)  =e^{-\gamma(\lambda)t}u_{\psi(\lambda)}(t,x)  \nonumber
		\\  =  &e^{-\gamma(\lambda)t}\Pi_x\left(\exp\left\{ \int_0^t(m-1) \xi(B_s,\omega)\mathrm{d}s \right\}e^{-\lambda B_t}\phi(B_t,\lambda)\right).
\end{align*}
By
\eqref{Many-to-one2},  one has
\begin{equation}\label{step_1}
		e^{-\lambda x}\phi(x,\lambda) =e^{-\gamma(\lambda)t}\mathrm{E}_{x}^{\xi} \left(\sum_{\nu \in N(t)} e^{-\lambda X_t^\nu}\phi(X_t^\nu,\lambda) \right).
\end{equation}
Let
$$W_t(\lambda):= e^{-\gamma(\lambda)t}\sum_{\nu \in N(t)} e^{-\lambda X_t^\nu}\phi(X_t^\nu,\lambda), \quad t\geq 0,$$
and $\widetilde{\mathcal{F}}_t$  be the $\sigma$-field generated by all information of the  branching Brownian motion up to time $t$.
Using  \eqref{step_1} and the Markov property of $\{X_t, t\geq 0\}$, we have the following lemma:
\begin{lemma}
	For any $|\lambda| > \rho$ and $\omega \in \Omega_1,$ $\left(W_t(\lambda), \widetilde{\mathcal{F}}_t, \mathrm{P}_x^\xi\right)$ is a positive martingale with mean $e^{-\lambda x}\phi(x,\lambda)$.
\end{lemma}

We omit the proof here, see \cite[Lemma 1.2]{RS} for a proof in a more general case.

The first purpose of the present paper is to study
the quenched limit of $W_t(\lambda)$ as $t\to\infty$  under Assumptions {\bf (H1)} and {\bf(H2)}, see Theorem \ref{Additive martingale}.
As a consequence, we will get that,
when  $\sum_{k=1}^\infty (k\ln k) p_k < \infty$,
$\mathbb{P}$-almost surely,
$$M_t/t \to v^*,\ \mathrm{P}_x^\xi-\mbox{a.s.}\quad \mbox{as } t \to \infty,$$
see Corollary \ref{Cor3.2}.

The second purpose of this paper is to prove an annealed invariance principle for $M_t$ under Assumptions {\bf (H1)}, {\bf(H2)} and {\bf (H3)}.

Conditioned on $\xi$, for any non-negative Borel function $f$ on $\R$,
$u(t,x):=u_f(t,x)=\mathrm{E}_x^\xi\langle X_t,f\rangle$
solves the following parabolic Anderson problem:
 \begin{equation}\label{PAM}
	\left\{\begin{array}{rl}&\frac{\partial}{\partial t}u(t,x)=\displaystyle\frac{1}{2}u_{xx} + (m-1)\xi(x,\omega) u(t,x),\quad x\in\R,\\
&\displaystyle u(0,x)=f(x).
\end{array}\right.
\end{equation}
According to \cite[Section 1]{Nolen},
\cite[Section 6]{JD} and  \cite[Section 4]{DS}, the maximal position
$M_t$ of $X_t$
is related to the front of the solution $u_f$. We write $u^{(-\infty,0]}$  for the solution to \eqref{PAM}  with initial condition $u(0,x)=1_{(-\infty,0]}(x)$. The front of $u^{(-\infty,0]}$ is defined as
\begin{equation}\label{front-linear}
	\overline m^\epsilon(t)=\sup\{x\in\R:  u^{(-\infty,0]}(t,x)\geq \epsilon\}.
\end{equation}
\cite[Theorem 1.4]{DS} proved that there exist a constant $c\in(0,\infty)$  and a $\P$-a.s. finite random time $T(\omega)$ such that for all $t \geq T(\omega)$,
$$
\overline m^\epsilon(t)-m^\epsilon(t) \leq c \ln t.
$$
The processes $\{\overline m^\epsilon(t), t\geq 0\}$ and $\{m^\epsilon(t), t\geq 0\}$ satisfy invariance principles,  see  \cite[Theorem 1.3, Corollary 1.5]{DS}.

Nolen \cite{Nolen} studied the front of $u_f$ with $f=\psi(x,-\lambda^*)$.
Note that
\begin{equation}\label{u-psi}
u_{\psi(-\lambda^*)}(t,x)= e^{\gamma(-\lambda^*)t} \psi(x,-\lambda^*,\omega) = e^{\gamma(\lambda^*)t}e^{\lambda^* x}\phi(x,-\lambda^*,\omega)
\end{equation}
 solves \eqref{PAM} with $f=\psi(x,-\lambda^*)$.
The position of the wave $u_{\psi(-\lambda^*)}$ at time $t$ is defined  in \cite{Nolen} by
$$
\overline m^{\epsilon}_\psi(t,\omega):=
\sup \left\{
x\in \mathbb{R}: u_{\psi(-\lambda^*)}(t,x,\omega)=\epsilon\right\}.$$
Using the fact that $\log(\psi(x, -\lambda^*, \omega)) \sim \lambda^* x$ as $x\to\infty$,
we have
\begin{equation}\label{limit-overline-m}
\lim_{t\to\infty}\frac{\overline m^{\epsilon}_\psi(t,\omega)}{t}=-\frac{\gamma(\lambda^*)}{\lambda^*}=-v^*\quad \P\mbox{-a.s.},
\end{equation}
see a similar argument in \cite[(5)]{Nolen}.

In this paper, we will introduce
a process $V_t$ to play the role of $-\overline m^{\epsilon}_\psi(t,\omega)$.
 {\bf(H2)} implies that
$|-\lambda^*| > \rho$.
Taking $\lambda=-\lambda^*$ in  \eqref{step_1}
and using the fact that $\gamma(-\lambda^*)= \gamma(\lambda^*)$, we have
\begin{equation}\label{psi-lambda*}
e^{\gamma(\lambda^*)t}\psi(x,-\lambda^*) =\mathrm{E}_{x}^{\xi} \left(\sum_{\nu \in N(t)} \psi( X_t^\nu,-\lambda^*) \right).
\end{equation}
We claim that the main contribution to the sum  above is given by
\begin{equation}\label{equality-V_t}
V_t:= \sup\left\{x\in\R: \psi(x,-\lambda^*,\omega) =  e^{\gamma(\lambda^*)t} \right\}=\sup\left\{x\in\R: {\lambda^* x}+ \ln \phi(x,-\lambda^*,\omega) =\gamma(\lambda^*)t \right\}.
\end{equation}
Since $\psi(\cdot, -\lambda^*,\omega)$ is continuous, the supremum above is attained. Thus,
$
\psi(V_t, -\lambda^*)=e^{\gamma(\lambda^*)t}.
$
If $\nu\in N(t)$ is
such that $X_t^\nu=V_t$, then, since
$\psi(\cdot, -\lambda^*,\omega)$ is strictly increasing (see Remark \ref{remark3} below),  $\psi(X_t^\nu, -\lambda^*)$ is the main contribution to the sum of the right hand side of \eqref{psi-lambda*}.
By \eqref{u-psi}, we may rewrite our $V_t$ as
	$$V_t(\omega)= \sup\left\{x: u_{\psi(-\lambda^*)}(t,x,\omega) =  e^{2\gamma(\lambda^*)t} \right\}.$$
By \eqref{limit-overline-m},  we have
$$\lim_{t\to\infty}\frac{V_t}{t}=v^*\quad \P\mbox{-a.s.},$$
Although both $V_t$ and $\overline m^\epsilon(t)$ are  fronts of the linear problem \eqref{PAM}
at different levels,
 we find that $V_t$ is easier to handle than $\overline m^\epsilon (t)$ when we investigate the main contribution to $\sum_{\nu \in N(t)} \psi( X_t^\nu,-\lambda^*)$.

Our strategy to prove an invariance principle of $M_t$ is to first prove that  there exists a non-random constant
$\Gamma>0$
 such that for $\P$-a.s. $\omega$,
\begin{equation}\label{differ-V-M}
	\limsup_{t \to \infty} \frac{|M_t - V_t|}{\ln t} \leq \Gamma ,\qquad\mathrm{P}_0^\xi\mbox{-a.s.}
\end{equation}
see Theorem \ref{Log-distance} below.
Then we use \cite[Lemma 3.1]{Nolen} to get  an invariance principle for $V_t$. The invariance principle for $M_t$ will be obtained from the invariance principle of $V_t$ by \eqref{differ-V-M}.
It seems using $V_t$ as a tool to establish invariance principle for $M_t$ is new.

\subsection{Main results}\label{s:mainresults}
\begin{thrm}\label{Additive martingale}
If {\bf(H1)} and {\bf(H2)} hold, then
there exists $\Omega_2\subset \Omega_1$ with $\mathbb{P}(\Omega_2)=1$ such that for any $\omega \in \Omega_2$, the limit $W_\infty(\lambda) := \lim_{t\to\infty} W_t(\lambda)$ exists $\mathrm{P}_x^\xi$-a.s.  Moreover,

(i) if
$ |\lambda |\ge \lambda^*$
then $W_\infty(\lambda) = 0, \ \mathrm{P}_x^\xi$-a.s.;

	(ii) if $|\lambda| \in (\rho, \lambda^* ),$ then $W_\infty(\lambda)$ is an $L^1(\mathrm{P}_x^\xi)$ limit or $W_\infty(\lambda)=0$ according to $\sum_{k=1}^\infty (k\ln k) p_k < \infty$ or $\sum_{k=1}^\infty (k\ln k) p_k =\infty;$

\end{thrm}
\begin{cor}\label{Cor3.2}
Assume that {\bf(H1)} and {\bf(H2)} hold.
If  $\sum_{k=1}^\infty (k\ln k) p_k < \infty$, then $\mathbb{P}$-almost surely, $M_t/t \to v^*,\ \mathrm{P}_x^\xi$-a.s. as $t \to \infty$.
\end{cor}

\begin{thrm}\label{Log-distance}
Assume that  {\bf(H1)}, {\bf(H2)} and {\bf(H3)} hold. There exists a non-random constant $\Gamma>0$
such that for $\P$-a.s. $\omega$,
	$$ \limsup_{t \to \infty} \frac{|M_t - V_t|}{\ln t} \leq \Gamma ,\qquad\mathrm{P}_0^\xi\mbox{-a.s.}$$
\end{thrm}
\begin{thrm}\label{Invariance_of_M_t}
	Assume that {\bf(H1)}, {\bf(H2)} and {\bf(H3)} hold.

(i) Under $\mathbb{P} \times \mathrm{P}_0^\xi$, we have
	$$\frac{M_t - v^* t}{\sqrt{t}} \stackrel{\mathrm{d}}{\Rightarrow} \mathcal{N}(0, \widetilde{\sigma}_{-\lambda^*}^2),\qquad \textup{as}\ t\to+\infty,$$
where $\widetilde{\sigma}_{-\lambda^*}^2:=  (\sigma_{-\lambda^*}')^2 v^*/(\lambda^*)^2$ and $(\sigma_{-\lambda^*}')^2$ is defined in \eqref{Another_var}.

(ii) If $(\sigma_{-\lambda^*}')^2 > 0$, then the sequence of processes
		$$
	[0, \infty) \ni t \mapsto \frac{M_{nt}-v^* n t}{\tilde{\sigma}_{-\lambda^*} \sqrt{n}}, \quad n \in \mathbb{N}
	$$
	converges weakly as $n\to\infty$ to a standard Brownian motion on $[0,\infty)$, in the topology of $C[0,\infty)$.
\end{thrm}

\section{Preliminaries}

\subsection{Martingale change of measure and spine decomposition}

It follows from the relationship between $\phi$ and $\psi$ that
	\begin{align}\label{Psi_and_Phi}
		\ln \psi(x,\lambda) & = -\lambda x +\ln \phi(x,\lambda),\quad
		\frac{\psi_x(x,\lambda)}{\psi(x,\lambda)} = -\lambda + \frac{\phi_x(x,\lambda)}{\phi(x,\lambda)},\nonumber\\ \frac{\psi_\lambda(x,\lambda)}{\psi(x,\lambda)}& = -x +\frac{\phi_\lambda(x,\lambda)}{\phi(x,\lambda)},\quad   \frac{1}{x} \ln \psi(x,\lambda) \to -\lambda \textup{ as } |x| \to \infty.
	\end{align}

For any $\omega \in \Omega_1,$  using the non-negative martingale $\left(W_t(\lambda), \widetilde{\mathcal{F}}_t, t\geq 0; \mathrm{P}_x^\xi\right)$,  we define a new probability measure
$\mathrm{Q}_{x}^{\xi, \lambda}$ by
$$ \frac{\mathrm{d}\mathrm{Q}_{x}^{\xi, \lambda}}{\mathrm{d} \mathrm{P}_x^\xi}\bigg|_{\widetilde{\mathcal{F}}_t}(X) = \frac{W_t (\lambda)}{e^{-\lambda x}\phi(x,\lambda)}.$$

To help us understand $\{X_t, t\geq 0\}$ under $\mathrm{Q}_{x}^{\xi, \lambda}$, we
introduce a martingale  change of measure for the law $\Pi_x$ of Brownian motion $\{B_t, t\geq 0\}$.
For fixed $\omega \in \Omega_1$,
it follows from It\^o's formula that
\begin{align*}
	\mathrm{d} \ln \psi(B_t,\lambda) & = \frac{\psi_x(B_t,\lambda)}{\psi(B_t,\lambda)} \mathrm{d} B_t + \frac{1}{2} \frac{\psi_{xx}(B_t,\lambda) \psi(B_t,\lambda) -
		\left(\psi_x(B_t,\lambda)\right)^2}{\psi^2(B_t,\lambda)} \mathrm{d} t \\
	& = \frac{\psi_x(B_t,\lambda)}{\psi(B_t,\lambda)} \mathrm{d} B_t + \frac{1}{2} \frac{2\left( \gamma(\lambda) - (m-1)\xi(B_t) \right) \psi^2(B_t,\lambda) -
	\left(\psi_x(B_t,\lambda)\right)^2}{\psi^2(B_t,\lambda)} \mathrm{d} t,\quad\Pi_x\mbox{-a.s.}
\end{align*}
Since  $\Pi_x(B_0 =x)=1$,
$$\ln \frac{\psi(B_t,\lambda)}{\psi(x,\lambda)} = \int_0^t  \frac{\psi_x(B_s,\lambda)}{\psi(B_s,\lambda)} \mathrm{d} B_s- \frac{1}{2}\int_0^t \left(\frac{\psi_x(B_s,\lambda)}{\psi(B_s,\lambda)}\right)^2 \mathrm{d}s+ \int_0^t \left(\gamma(\lambda)-(m-1)\xi(B_s)\right) \mathrm{d}s.$$
Now we define
\begin{align*}
	\Upsilon_t
	:= &\exp\left\{\int_0^t  \frac{\psi_x(B_s,\lambda)}{\psi(B_s,\lambda)} \mathrm{d} B_s -\frac{1}{2}\int_0^t \left(\frac{\psi_x(B_s,\lambda)}{\psi(B_s,\lambda)}\right)^2 \mathrm{d}s \right\} \\
	=&\exp\left\{\ln \frac{\psi(B_t,\lambda)}{\psi(x,\lambda)}- \int_0^t\left(\gamma(\lambda)-(m-1)\xi(B_s)\right)\mathrm{d}s \right\}\\
=& \exp\left\{-\gamma(\lambda) t + \int_0^t (m-1)\xi(B_s)\mathrm{d}s \right\} \frac{\psi(B_t,\lambda)}{\psi(x,\lambda)},\quad t\geq 0.
\end{align*}
Then $\{\Upsilon_t, t\geq 0; \Pi_x\}$ is a martingale with respect to $\sigma(B_r, r\leq t)$. In fact,  by \eqref{F-K-u} and the  Markov property of $\{B_t, t\geq 0\}$, we have that for $t,s \geq 0$,
\begin{align*}
	\Pi_x [\Upsilon_{t+s}| \sigma(B_r: r \leq t)]&= \Upsilon_t \cdot \Pi_x \left(\exp\left\{\ln \frac{\psi(B_{t+s},\lambda)}{\psi(B_t,\lambda)}- \int_t^{t+s}\left(\gamma(\lambda)-(m-1)\xi(B_r)\right)\mathrm{d}r \right\} \bigg| B_t\right)\\
	& = \Upsilon_t\cdot \frac{e^{-\gamma(\lambda)s}u(s,B_t,\lambda)}{ \psi(B_t,\lambda)} = \Upsilon_t.
\end{align*}
 Define a probability measure $\widetilde{\Pi}_x^{\xi, \lambda}$ by
\begin{align}\label{Change-of-measure}
 	\frac{\mathrm{d}\widetilde{\Pi}_x^{\xi, \lambda}}{\mathrm{d} \Pi_x}\bigg|_{\sigma(B_s: s\leq t)}
 & = \Upsilon_t.\quad t\geq 0.
\end{align}
Suppose that $\left\{\Xi_t, t\geq 0;\widetilde{\Pi}_x^{\xi,\lambda}\right\}$ is a diffusion
such that
$$ \mathrm{d} \Xi_t = \mathrm{d} W_t +\frac{\psi_x(\Xi_t,\lambda)}{\psi(\Xi_t,\lambda)} \mathrm{d} t=\mathrm{d} W_t +\left(\frac{\phi_x(\Xi_t,\lambda)}{\phi(\Xi_t,\lambda)} -\lambda\right)\mathrm{d} t,$$
where $\{W_t, t\geq 0; \widetilde{\Pi}_x^{\xi, \lambda}\}$ is a Brownian motion starting from $x$ at time $t=0$.
Applying Girsanov's theorem, we get
\begin{equation}\label{=distribution}\left\{B_t, t\ge0;  \widetilde{\Pi}_x^{\xi, \lambda}\right\}\stackrel{\mathrm{d}}{=} \left\{\Xi_t,t\geq0; \widetilde{\Pi}_x^{\xi, \lambda}\right\}.\end{equation}

Now we  define a new process $\{\widetilde{X}_t, t\geq 0\}$
with probability $ \widetilde{\mathrm{P}}_{x}^{\xi, \lambda}$ such that $\{\widetilde{X}_t, t\geq 0;  \widetilde{\mathrm{P}}_{x}^{\xi, \lambda}\}$ has the same  law as $\{{X}_t, t\geq 0;\mathrm{Q}_{x}^{\xi, \lambda}\}$.
More precisely,  we consider a branching particle system in which

 (i) there is an initial marked particle at $x \in \mathbb{R}$ which moves according to $(\Xi_t, \widetilde{\Pi}_x^{\xi, \lambda})$;

(ii) the branching rate of this marked particle is $m\xi(y)$ at site $y$;

 (iii) when a marked particle dies at site $y$,
it gives birth to $k$ children with $\tilde{p}_k = k p_k/m, k\geq 1$;

(iv)  one of these children is uniformly selected and marked, and
the marked child evolves independently like its parent
and the other children evolve independently with law $\mathrm{P}_y^\xi.$
Define
   $$\widetilde X_t:= \sum_{\nu\in \widetilde N(t)} \delta_{\widetilde X_t^\nu},\quad t\geq 0,$$
where
$\widetilde N(t)$ be the set of particles alive at time $t$ and $\widetilde X_t^\nu$ is the position of particle $\nu\in \widetilde N(t).$  Then by \cite[Theorem 2.9]{RS},
$$\{{X}_t, t\geq 0;\mathrm{Q}_{x}^{\xi, \lambda}\}\stackrel{\mathrm{d}}{=} \{\widetilde{X}_t, t\geq 0;  \widetilde{\mathrm{P}}_{x}^{\xi, \lambda}\}.$$
The set of the marked particles along with their trajectories is called a \emph{spine}.

\subsection{Properties of $\gamma(\lambda)$ and $\phi(x,\lambda)$}

We only assume that {\bf(H1)} holds in this subsection. We will give some basic properties of $\gamma(\lambda)$ and $\phi(\cdot,\lambda)$ for $|\lambda| > \rho$.
The case that $\lambda > \rho$ and $\lambda < -\rho$ are similar.
We will state our results for $|\lambda|>\rho$, but only prove the case $\lambda > \rho$.
By \cite[Theorem 5.1, \S 7.4]{FM}, $\gamma(\lambda)$ is differentiable, strictly convex
and $\gamma'(\lambda) > 0$ for $\lambda > \rho$,
and thus $\gamma'(\lambda)$ is strictly  increasing.

\begin{lemma}\label{Engenvalue-Property}
(1) $\gamma(\lambda^*)=\lambda^* \gamma'(\lambda^*)$, $v^*=\gamma'(\lambda^*)$.
(2)	If $|\lambda| > \lambda^*$, then $\lambda\gamma'(\lambda) > \gamma(\lambda)$
(3) if $\rho< |\lambda| < \lambda^*$, then $\lambda\gamma'(\lambda) < \gamma(\lambda)$.
\end{lemma}
\textbf{Proof:} By the definition of $\lambda^*$, we know that
\begin{align}\label{step_52}
	0=\frac{\mathrm{d}}{\mathrm{d}\lambda}\left(\frac{\gamma(\lambda)}{\lambda}\right)\bigg|_{\lambda=\lambda^*}=\frac{\gamma'(\lambda^*)\lambda^* - \gamma(\lambda^*)}{\left(\lambda^*\right)^2}\ \Leftrightarrow\ \gamma(\lambda^*)=\lambda^* \gamma'(\lambda^*).
\end{align}
For $\lambda > \rho$,
\begin{align}\label{equal_1}
	\lambda \gamma'(\lambda) - \gamma(\lambda) &= \lambda \gamma'(\lambda) - \gamma(\lambda) -\left(\lambda^* \gamma'(\lambda^*) - \gamma(\lambda^*) \right)\nonumber\\
	&= \lambda^* \left(\gamma'(\lambda) - \gamma'(\lambda^*)\right)+ \int_{\lambda*}^\lambda \left(\gamma'(\lambda) -\gamma'(y)\right)\mathrm{d}y.
\end{align}
Since $\gamma'(\lambda)$ is strictly increasing, the conclusions of the lemma follow immediately.
\hfill$\Box$
\bigskip

 Recall that $\psi(x,\lambda) = e^{-\lambda x}\phi(x,\lambda)$ satisfies \eqref{Principal-Engenfunction-2}. Thus, for fixed $\lambda > 0$, $\psi(x,\lambda) \to 0 $ as $x \to +\infty$.
 For any $y\in\R$, define
 \begin{equation}\label{define-H}H_y:= \inf\left\{t>0:\ B_t=y \right\}.\end{equation}
  According to the Feynman-Kac formula,
 we have
\begin{equation}\label{Psi_2}
	\psi(x, \lambda) = \Pi_x\left(\exp\left\{ \int_0^{H_y}\left((m-1)\xi(B_s) -\gamma(\lambda) \right)\mathrm{d}s \right\} \right)\psi(y,\lambda),\quad
x>y.
\end{equation}
In particular,
for $x>0$,
\begin{equation}\label{Psi}
	\psi(x, \lambda) = \Pi_x\left(\exp\left\{ \int_0^{H_0}\left((m-1)\xi(B_s) -\gamma(\lambda) \right)\mathrm{d}s \right\} \right).
\end{equation}
Recall that
${\cal F}_x=\sigma\{\xi(y): y\le x\}$ and ${\cal F}^x= \sigma\{\xi(y): y\ge x\}$.
By \eqref{Psi_2}, we have for any $x>y$,
\begin{equation}\label{measurable-psi}\
\frac{\psi(x, \lambda)}{\psi(y, \lambda)}\in {\cal F}^y.
\end{equation}
Similarly, for fixed $\lambda<0$, $\psi(x, \lambda)\to 0$ as $x\to-\infty$.
It follows from the Feynman-Kac formula that
\begin{equation}\label{Psi_2-2}
	\psi(x, \lambda) = \Pi_x\left(\exp\left\{ \int_0^{H_y}\left((m-1)\xi(B_s) -\gamma(\lambda) \right)\mathrm{d}s \right\} \right)\psi(y,\lambda),\quad
x<y.
\end{equation}
Thus
\begin{equation}\label{measurable-psi-2}\
\frac{\psi(x, \lambda)}{\psi(y, \lambda)}\in {\cal F}_y,
\quad \mbox{ for } x<y, \lambda<0.
\end{equation}

\begin{lemma}
For all $|\lambda|>\rho$ and $\omega\in \Omega_1$, $\phi(\cdot,\lambda,\omega)$ is differentiable in $x$ and $\lambda$.
\end{lemma}

\textbf{Proof:}
Without loss of generality, we assume $\lambda\in (\rho,\infty)$.
We claim that, for fixed $x > 0$,
$\psi$ is differentiable as a function of $\lambda$  and
\begin{equation}\label{Psi_Lambda}
	\psi_\lambda(x,\lambda):= \frac{\partial\psi(x, \lambda)}{\partial \lambda} = -\gamma'(\lambda) \Pi_x\left(H_0\exp\left\{ \int_0^{H_0}\left((m-1)\xi(B_s) -\gamma(\lambda) \right)\mathrm{d}s \right\} \right).
\end{equation}
In fact,  we can choose $\lambda_1\in (\rho,  \lambda)$ and constant $K>0$ such that $h \exp\left\{-\left(\gamma(\lambda) -\gamma (\lambda_1) \right)h\right\}\leq K$ for all $h >0$,
and thus
$$
\Pi_x\left(H_0\exp\left\{ \int_0^{H_0}\left((m-1)\xi(B_s) -\gamma(\lambda) \right)\mathrm{d}s \right\} \right) \leq K \Pi_x\left(\exp\left\{ \int_0^{H_0}\left((m-1)\xi(B_s) -\gamma(\lambda_1) \right)\mathrm{d}s \right\} \right).
$$
We can now use the dominated convergence theorem in \eqref{Psi} to get \eqref{Psi_Lambda}.

The same argument together with \eqref{Psi_2} shows that, for $x>y$,
$\psi(x,\lambda) /\psi(y,\lambda)$ is differentiable
in $\lambda\in (\rho,\infty)$.
By choosing $x>0\geq  y$, we get that $\psi$ is differentiable in $\lambda$ for all $x\in \mathbb{R}$.
The desired conclusion follows immediately from \eqref{def-psi}.
\hfill$\Box$
\bigskip

To compare with the constant case in which $\phi\equiv1$,  we will use $\phi$  in the statements of the results, but use $\psi$ in the proof for convenience.
By the uniqueness of $\phi$ ,
we have $\P$-almost surely $\phi(z+y,\lambda,\omega) = \phi(y,\lambda ,\omega) \phi(z, \lambda,\theta_y \omega).$
Thus by \eqref{Psi_and_Phi}, for any $z,y\in\R$,
\begin{align}\label{step_3-1}
\ln \psi(z+y,\lambda,\omega ) - \ln \psi(y,\lambda,\omega ) & =  -\lambda z + \ln \phi(z+y,\lambda,\omega )  - \ln \phi(y,\lambda,\omega ) \nonumber \\
	& = -\lambda z+ \ln \phi(z,\lambda, \theta_{y}\omega).
	\end{align}
Taking derivative with respect to $\lambda$, we get
\begin{align}\label{step_3-2}
		\frac{\psi_\lambda(z+y,\lambda,\omega)}{\psi(z+y,\lambda,\omega)}-\frac{\psi_\lambda(y,\lambda,\omega)}{\psi(y,\lambda,\omega)}  & = -z + \frac{\phi_\lambda(z,\lambda,\theta_y \omega)}{\phi(z,\lambda,\theta_y \omega)}.
\end{align}
Note also that  \eqref{step_3-1} is equivalent to
\begin{align}\label{step_3-1'}
\frac{\psi(z+y,\lambda,\omega)}{\psi(y,\lambda,\omega)}=e^{-\lambda z}\phi(z,\lambda,\theta_y\omega).
\end{align}

\begin{lemma}\label{Engenfunction-Property1}
Suppose $|\lambda| > \rho$ and  $\gamma(\lambda) > (m-1)\textup{es}$.

(1) For all $|\lambda|>\rho$, there exist
positive constants $C_1(\lambda),C_2(\lambda)$
depending on $\lambda$ only
such that
\begin{equation}\label{lambda>rho}
-C_1(\lambda) \leq
\left(\frac{\psi_\lambda(x,\lambda)}{\psi(x,\lambda)}\right)_x= \left(-x +\frac{\phi_\lambda(x,\lambda)}{\phi(x,\lambda)}\right)_x
\leq -C_2(\lambda), \quad x\in \mathbb{R};
\end{equation}

 (2)
If $\lambda > \rho$, there exist
positive constants  $K_1(\lambda),K_2(\lambda)$ depending on $\lambda$ only
such that
\begin{equation}\label{lambda>rho2}
-K_1(\lambda) \leq
\frac{\psi_x(x,\lambda)}{\psi(x,\lambda)} =
-\lambda +\frac{\phi_x(x,\lambda)}{\phi(x,\lambda)} \leq
-K_2(\lambda), \quad x\in \mathbb{R};
\end{equation}
and  if $\lambda < -\rho$,
	\begin{equation}\label{lambda<rho2}
K_2(\lambda) \leq
\frac{\psi_x(x,\lambda)}{\psi(x,\lambda)} =
  -\lambda +\frac{\phi_x(x,\lambda)}{\phi(x,\lambda)} \leq
  K_1(\lambda), \quad x\in \mathbb{R}.
\end{equation}

(3)	$$\lim_{|x|\to\infty} \frac{\phi_\lambda(x,\lambda)}{x\phi(x,\lambda)} = 0,\qquad \mathbb{P}\textup{-a.s.}$$

(4) $\gamma'(\lambda)$ is continuous in the open set $\{\lambda\in\R: |\lambda|> \rho, \gamma(\lambda)> (m-1)\textup{es}\}$.
	
\end{lemma}
\textbf{Proof:} Recall that, by \eqref{Psi_and_Phi},  $\frac{\psi_\lambda(x,\lambda)}{\psi(x,\lambda)} = -x +\frac{\phi_\lambda(x,\lambda)}{\phi(x,\lambda)}$. Put $g(x,\lambda):= \frac{\psi_\lambda(x,\lambda)}{\psi(x,\lambda)}$ for simplicity.

(1) Fix $\lambda > \rho$.
For $x > y$, taking logarithm
and then differentiating in $\lambda$ in \eqref{Psi_2},
 we get
\begin{equation}\label{Psi_3}
	\frac{\psi_\lambda(x,\lambda)}{\psi(x,\lambda)} =-\gamma'(\lambda) \frac{\Pi_x\left(\exp\left\{ \int_0^{H_y}\left((m-1)\xi(B_s) -\gamma(\lambda) \right)\mathrm{d}s \right\}H_y \right)}{\Pi_x\left(\exp\left\{ \int_0^{H_y}\left((m-1)\xi(B_s) -\gamma(\lambda) \right)\mathrm{d}s \right\} \right)} + \frac{\psi_\lambda(y,\lambda)}{\psi(y,\lambda)}.
\end{equation}
Hence, when $x>y$,
\begin{align}\label{step_14}
 g(x,\lambda) - g(y,\lambda)
 = &-\gamma'(\lambda) \frac{\Pi_x\left(\exp\left\{ \int_0^{H_y}\left((m-1)\xi(B_s) -\gamma(\lambda) \right)\mathrm{d}s \right\}H_y \right)}{\Pi_x\left(\exp\left\{ \int_0^{H_y}\left((m-1)\xi(B_s) -\gamma(\lambda) \right)\mathrm{d}s \right\} \right)} \nonumber \\
		 \leq &-\gamma'(\lambda) \frac{\Pi_x\left(\exp\left\{ \left((m-1)\textup{ei} -\gamma(\lambda) \right)H_y \right\}H_y \right)}{\Pi_x\left(\exp\left\{ \left((m-1)\textup{es} -\gamma(\lambda) \right)H_y \right\} \right)} \nonumber \\
		 = &\frac{-\gamma'(\lambda)(x-y)}{\sqrt{2\left(\gamma(\lambda) - (m-1)\textup{ei}\right)}}\cdot \frac{e^{-(x-y)\sqrt{2\left(\gamma(\lambda) - (m-1)\textup{ei}\right)}} }{e^{-(x-y)\sqrt{2\left(\gamma(\lambda) - (m-1)\textup{es}\right)}}}.
\end{align}
In the first inequality above we used the fact that $-\gamma'(\lambda) <0$ and in  last equality the fact that for any $u > 0$ and $x,y\in\mathbb{R}$,
\begin{align}\label{Exp-of-H}
	\Pi_x\left(e^{-u H_y}\right) = e^{-|y-x|\sqrt{2u}} \qquad \textup{ and }\qquad \Pi_x\left( H_y e^{-uH_y}\right) = \frac{|y-x|}{\sqrt{2u}}e^{-|y-x|\sqrt{2u}}.
\end{align}
Dividing both sides of \eqref{step_14} by $x-y$ and letting $y \uparrow x$, we get
$$g_x(x,\lambda) \leq \frac{-\gamma'(\lambda)}{\sqrt{2\left(\gamma(\lambda) - (m-1)\textup{ei}\right)}}=:
-C_2(\lambda), \quad x\in\R.$$
Similarly, for $x>y$,
\begin{align}
	 g(x,\lambda) - g(y,\lambda)
		& \geq -\gamma'(\lambda) \frac{\Pi_x\left(\exp\left\{ \left((m-1)\textup{es} -\gamma(\lambda) \right)H_y \right\}H_y \right)}{\Pi_x\left(\exp\left\{ \left((m-1)\textup{ei} -\gamma(\lambda) \right)H_y \right\} \right)} \nonumber \\
		& = \frac{-\gamma'(\lambda)(x-y)}{\sqrt{2\left(\gamma(\lambda) - (m-1)\textup{es}\right)}}\cdot \frac{e^{-(x-y)\sqrt{2\left(\gamma(\lambda) - (m-1)\textup{es}\right)}} }{e^{-(x-y)\sqrt{2\left(\gamma(\lambda) - (m-1)\textup{ei}\right)}}},
\end{align}
and thus
$$g_x(x,\lambda) \geq \frac{-\gamma'(\lambda)}{\sqrt{2\left(\gamma(\lambda) - (m-1)\textup{es}\right)}} =:
-C_1(\lambda), \quad x\in\R.$$
Hence \eqref{lambda>rho} is valid.

When $\lambda < -\rho$,  \eqref{Psi_3} holds for all $x < y$. In this case we have that for $x <y$,
$$
g(y,\lambda) - g(x,\lambda)
 = \gamma'(\lambda) \frac{\Pi_x\left(\exp\left\{ \int_0^{H_y}\left((m-1)\xi(B_s) -\gamma(\lambda) \right)\mathrm{d}s \right\}H_y \right)}{\Pi_x\left(\exp\left\{ \int_0^{H_y}\left((m-1)\xi(B_s) -\gamma(\lambda) \right)\mathrm{d}s \right\} \right)}<0.$$
An argument similar as above shows that \eqref{lambda>rho} holds.

(2) The argument is very similar to that of (1).  We only prove the case of $\lambda > \rho$. By \eqref{Psi_2}, for $x > y$ we have
\begin{align*}
	&\int_y^x \frac{\psi_x(z,\lambda)}{\psi(z,\lambda)} \mathrm{d}z  = \ln \psi(x,\lambda) - \ln \psi(y,\lambda)
	 = \ln \Pi_x\left(\exp\left\{ \int_0^{H_y}\left((m-1)\xi(B_s) -\gamma(\lambda) \right)\mathrm{d}s \right\} \right) \\
	& \leq \ln \Pi_x\left(\exp\left\{ \int_0^{H_y}\left((m-1)\textup{es} -\gamma(\lambda) \right)\mathrm{d}s \right\} \right)  = -(x-y)\sqrt{2\left(\gamma(\lambda)-(m-1)\textup{es} \right)}.
\end{align*}
Similarly, for any $x >y$,
$$ \int_y^x \frac{\psi_x(z,\lambda)}{\psi(z,\lambda)} \mathrm{d}z \geq -(x-y)
\sqrt{2\left(\gamma(\lambda)-(m-1)\textup{ei} \right)}.
$$
Letting
$K_2(\lambda):= \sqrt{2\left(\gamma(\lambda)-(m-1)\textup{es} \right)}$
and
$K_1(\lambda):= \sqrt{2\left(\gamma(\lambda)-(m-1)\textup{ei} \right)}$,
we get the results of (2).

(3) For $\lambda > \rho,$ in view of (1), for any $x,y\in\mathbb{R}$ with $ |x-y|\leq 1$,
\begin{equation}\label{step_45}
	|g(x,\lambda) -g(y,\lambda)|=\left|\frac{\psi_\lambda(x,\lambda)}{\psi(x,\lambda)}-\frac{\psi_\lambda(y,\lambda)}{\psi(y,\lambda)}\right|= \left|\frac{\phi_\lambda(x,\lambda)}{\phi(x,\lambda)}- \frac{\phi_\lambda(y,\lambda)}{\phi(y,\lambda)}\right|
\leq C_1(\lambda).
\end{equation}
Thus, it suffices to prove that for $k\in\mathbb{Z}$,
$$\lim_{|k|\to\infty} \frac{g(k,\lambda)}{k} = -1,\qquad \mathbb{P}\textup{-a.s.}$$
By
\eqref{step_3-2}
 and \eqref{step_45}, $\left\{g(k+1,\lambda,\omega) - g(k,\lambda ,\omega): k\in \mathbb{Z}\right\}$ is a stationary and ergodic sequence with uniform bound
$C_1(\lambda).$
By \eqref{step_3-1} and \eqref{lambda>rho2}, we also have that
 $\{\ln \psi(k+1, \lambda,\omega) -\ln \psi(k,\lambda,\omega):k\in\mathbb{Z} \}$ is a stationary and ergodic sequence with uniform bound
$K_1(\lambda)$.
Thus,
$$\mathbb{E}\ln \psi(1,\lambda,\omega)= \lim_{k\to\infty} \frac{\ln \psi(k,\lambda,\omega)}{k} = -\lambda.$$
Taking derivative with respect to  $\lambda$  in the display above, using  the boundness of $g$, we get that
\begin{equation}\label{Exp-g}
	\mathbb{E}g(1,\lambda,\omega)= -1.
\end{equation}
By Birkhoff's ergodic theorem,
$$ \lim_{|k|\to\infty} \frac{g(k,\lambda)}{k} = \mathbb{E}g(1,\lambda ,\omega)= -1,\qquad \mathbb{P}\textup{-a.s.}$$
This completes the proof of  (3).

 (4) Suppose $\lambda> \rho$ and $\gamma(\lambda)> (m-1)\textup{es}$. Using \eqref{Psi_3}, we have
$$
g(1,\lambda)=\frac{\psi_\lambda (1,\lambda)}{\psi(1,\lambda)} = -\gamma'(\lambda)
\left(\frac{\Pi_1\left(\exp\left\{ \int_0^{H_0}\left((m-1)\xi(B_s) -\gamma(\lambda) \right)\mathrm{d}s \right\}H_0 \right)}{\Pi_1\left(\exp\left\{ \int_0^{H_0}\left((m-1)\xi(B_s) -\gamma(\lambda) \right)\mathrm{d}s \right\} \right)} \right).
$$
Taking expectation with respect to $\mathbb{P}$, and using \eqref{Exp-g}, we get
$$
1= \gamma'(\lambda)\mathbb{E}
\left(\frac{\Pi_1\left(\exp\left\{ \int_0^{H_0}\left((m-1)\xi(B_s) -\gamma(\lambda) \right)\mathrm{d}s \right\}H_0 \right)}{\Pi_1\left(\exp\left\{ \int_0^{H_0}\left((m-1)\xi(B_s) -\gamma(\lambda) \right)\mathrm{d}s \right\} \right)} \right).
$$
Using the bounded convergence theorem and the continuity of $\gamma(\lambda)$, we get $\gamma'(\lambda)$ is continuous in $\lambda$.
\hfill$\Box$

\begin{remark}\label{remark3}
By Lemma \ref{Engenfunction-Property1}(2),
if {\bf(H2)} holds, then
 $\psi(\cdot ,-\lambda^*)$ is increasing and for any $y<z$,
	\begin{equation}\label{difference-psi}
		\ln \psi(z, -\lambda^*) - \ln \psi(y,-\lambda^*) =
	\int_y^z \left(\ln \psi(x,-\lambda^*)\right)_x \mathrm{d} x
		 \in \left[K_2(-\lambda^*)(z-y), K_1(-\lambda^*)(z-y)\right].
	\end{equation}
Thus, for any $y<z$,
	\begin{equation}\label{difference-psi'}
K_2(-\lambda^*)(z-y)\leq \ln\psi(z,-\lambda^*)- \ln \psi(y,-\lambda^*)\leq K_1(-\lambda^*)(z-y).
	\end{equation}
Recalling the definition of $V_t$ in
\eqref{equality-V_t}, and using \eqref{difference-psi'},
 we get that for any $s<t$,
	 \begin{equation}\label{difference-V} \frac{\gamma(\lambda^*)(t-s)}{K_1(-\lambda^*)}\leq V_t -V_s \leq  \frac{\gamma(\lambda^*)(t-s)}{K_2(-\lambda^*)}.
\end{equation}
\end{remark}

Throughout this paper,
for any real $x$,
we use $[x]$ to denote the integer part of $x$, and $\lceil x \rceil$ to denote the smallest integer larger than $x$.
The following lemma is from \cite[Theorem 5.5]{HH}, see also  \cite[Theorem 2.8]{Nolen}.

\begin{lemma}\label{CLT}
Suppose that  $\{\eta_k \}_{k=0}^\infty \subset L^2(\Omega, \mathcal{F}, \mathbb{P})$ is a stationary sequence with $\mathbb{E}[\eta_k] = 0$ and that
$$
\sum_{k=1}^\infty \left(\eta_0 - \mathbb{E}(\eta_0 \vert \mathcal{F}_k) \right),\quad and \quad \sum_{k=1}^\infty \mathbb{E}(\eta_k \vert \mathcal{F}_0)
$$
converge in $L^2(\Omega, \mathcal{F}, \mathbb{P})$. Then the limit
$$
\sigma^2:= \lim_{N\to\infty} \frac{1}{N}\mathbb{E}\left(\sum_{i=0}^{N-1} \eta_k \right)^2
$$
exists and is finite. If $\sigma^2>0$,
	then $$\frac{\sum_{i=0}^{N-1} \eta_k }{ \sqrt{N}}\stackrel{\mathrm{d}}{\Longrightarrow}\mathcal{N}(0,\sigma^2)$$
and the process
$$
\frac{1}{\sigma\sqrt{n}}\left(\sum_{i=0}^{[nt]-1} \eta_k +(nt-[nt])\eta_{[nt]} \right),\qquad t\in [0,1],
$$
converges weakly to a standard Brownian motion on $[0,1].$
\end{lemma}

\begin{lemma}\label{Engenfunction-Property2}
	Suppose $|\lambda| > \rho$ and $\gamma(\lambda) > (m-1)\textup{es}$, then there exists
$\sigma_\lambda^2\geq 0$ such that under $\mathbb{P}$,
	\begin{equation}\label{convergence-in-d}\frac{\ln \phi(x,\lambda) - \lambda \frac{\phi_\lambda(x,\lambda)}{\phi(x,\lambda)}}{\sqrt{|x|}} \stackrel{\mathrm{d}}{\Longrightarrow} \mathcal{N}(0,\sigma_\lambda^2),\qquad \textup{as } |x|\to\infty.
	\end{equation}
	Moreover, if $\sigma_\lambda^2 > 0$, then
	the process
	$$\frac{\ln \phi(nt,\lambda) - \lambda \frac{\phi_\lambda(nt,\lambda)}{\phi(nt,\lambda)}}{\sigma_\lambda\sqrt{n}},\quad t\in [0,\infty),$$
	converge weakly to a standard Brownian motion on $[0,\infty)$ in the Skorohod topology and
	$$\limsup_{|x|\to+\infty} \frac{\ln \phi(x,\lambda) - \lambda \frac{\phi_\lambda(x,\lambda)}{\phi(x,\lambda)}}{\sqrt{|x|}} = +\infty ,\qquad \mathbb{P}\textup{-a.s.}$$
	If $\sigma_\lambda^2 =0$, then
	$$\sup_{x\in\R} \left\Vert \ln \phi(x,\lambda) - \lambda \frac{\phi_\lambda(x,\lambda)}{\phi(x,\lambda)} \right\Vert_\infty < \infty .$$
\end{lemma}
\textbf{Proof:} By \eqref{Psi_and_Phi}, we have that
\begin{equation}\label{relation-lnpsi-lnphi}\ln \phi(x,\lambda) - \lambda \frac{\phi_\lambda(x,\lambda)}{\phi(x,\lambda)} = \ln \psi(x,\lambda) - \lambda \frac{\psi_\lambda(x,\lambda)}{\psi(x,\lambda)}
\end{equation}
holds for all $|\lambda|>\rho$ and $x\in\mathbb{R}.$ In view of Lemma \ref{Engenfunction-Property1}(1)-(2), to prove \eqref{convergence-in-d}, it suffices to prove it for the case of $x\in \mathbb{Z}.$

Fix $\lambda > \rho$. Let $x = n\in \mathbb{Z}$
 and define
\begin{equation}\label{U_i}
	U_i:= \ln \psi(i+1,\lambda,\omega ) - \ln \psi(i,\lambda,\omega ) - \lambda
	\left(\frac{\psi_\lambda(i+1,\lambda,\omega)}{\psi(i+1,\lambda,\omega)}-\frac{\psi_\lambda(i,\lambda,\omega)}{\psi(i,\lambda,\omega)} \right),
\quad i\in\mathbb{Z}.
\end{equation}
Then by \eqref{lambda>rho} and \eqref{lambda>rho2},
$\{U_i\}_{i\in \mathbb{Z}}$ is a stationary and ergodic sequence with uniform bound
$C_1(\lambda) +\lambda K_1(\lambda)$,
and
$$\ln \psi(n,\lambda) - \lambda \frac{\psi_\lambda(n,\lambda)}{\psi(n,\lambda)} = \sum_{i=0}^{n-1}U_i,\qquad \ln \psi(-n,\lambda) - \lambda \frac{\psi_\lambda(-n,\lambda)}{\psi(-n,\lambda)} = \sum_{i=1}^{n}-U_{-i}.$$
Since
$$\ln \psi(n,\lambda) - \lambda \frac{\psi_\lambda(n,\lambda)}{\psi(n,\lambda)} \stackrel{\mathrm{d}}{=}\lambda \frac{\psi_\lambda(-n,\lambda)}{\psi(-n,\lambda)}- \ln \psi(-n,\lambda)$$
by the stationarity of $U_i$, we only need to
show \eqref{convergence-in-d}  for $x=n\in \mathbb{Z}^+$.
By Lemma \ref{Engenfunction-Property1}(3) and Birkhoff's ergodic theorem, $\mathbb{P}$-almost surely,
$$\mathbb{E}U_0 = \lim_{n\to+\infty} \frac{\sum_{i=0}^{n-1}U_i }{n} = \lim_{n\to+\infty} \left(\frac{\ln \psi(n,\lambda)}{n} - \lambda \frac{\psi_\lambda(n,\lambda)}{n\psi(n,\lambda)}\right) = -\lambda +\lambda = 0.$$

We will prove \eqref{convergence-in-d} using Lemma \ref{CLT}.
First we check the conditions of Lemma \ref{CLT} are satisfied using the argument of
\cite[Lemma A.2]{DS}.
For $j \geq 2$, set
\begin{align*}
A_1(\lambda)
& = A_1(\lambda, \omega;j) := \Pi_{1}
\left(\exp\left\{ \int_0^{H_{0}}\left((m-1)\xi(B_s) -\gamma(\lambda) \right)\mathrm{d}s \right\}; \sup_{0\leq s\leq H_0} (B_s -1) <  [j/2] \right); \\
A_2(\lambda)&= A_2(\lambda, \omega;j) := \Pi_{1}
\left(\exp\left\{ \int_0^{H_{0}}\left((m-1)\xi(B_s) -\gamma(\lambda) \right)\mathrm{d}s \right\} ;\sup_{0\leq s\leq H_0} (B_s -1) \geq  [j/2] \right)
\end{align*}
and $A_i'(\lambda) := \frac{\partial A_i(\lambda)}{\partial \lambda}$ for $i=1,2.$  According to
\eqref{Psi_2} and \eqref{Psi_3}, we have
$$U_0 = \ln
\left(A_1(\lambda)+ A_2(\lambda)\right)
-\lambda \frac{A_1'(\lambda)}{A_1(\lambda)+A_2(\lambda)} -\lambda \frac{A_2'(\lambda)}{A_1(\lambda)+A_2(\lambda)}.$$
Note that $A_1(\lambda)$
is $\mathcal{F}_{[j/2]+1}\cap \mathcal{F}^0$- measurable.
By the assumption that $\gamma(\lambda) > (m-1)\textup{es}$, we have
$$
A_1 (\lambda)\leq \Pi_{1}
\left(\exp\left\{ \left((m-1)\textup{es} -\gamma(\lambda) \right)H_{0} \right\} \right)
=:c_1(\lambda)<1,
$$
and for all $j\geq 2$,
\begin{align*}
A_1(\lambda) & \geq \Pi_{1}
\left(\exp\left\{ \left((m-1)\textup{ei} -\gamma(\lambda) \right)H_0 \right\}; \sup_{0\leq s\leq H_0} (B_s -1) <  [j/2] \right) \\ & \geq  \Pi_{1}
\left(\exp\left\{ \left((m-1)\textup{ei} -\gamma(\lambda) \right)H_0 \right\}; \sup_{0\leq s\leq H_0} (B_s -1) <  1\right)
=:c_2(\lambda) > 0.
\end{align*}
Therefore, for $j\geq 2$, we have
\begin{equation}\label{domi-A1}0<c_2(\lambda)\leq A_1(\lambda)\leq c_1(\lambda)<1.
\end{equation}
Next, if $H_0 >j $, then
$$
\exp\left\{ \int_0^{H_{0}}\left((m-1)\xi(B_s) -\gamma(\lambda) \right)\mathrm{d}s \right\} \leq \exp\left\{\left((m-1)\textup{es} -\gamma(\lambda) \right)j \right\};
$$
if $H_0 \leq j$, then
$
\exp\left\{ \int_0^{H_{0}}\left((m-1)\xi(B_s) -\gamma(\lambda) \right)\mathrm{d}s \right\} \leq 1
$, and $$\left\{\sup_{0\leq s\leq H_0} (B_s -1) \geq  [j/2]  \right\} \subset \left\{\sup_{0\leq s\leq j} (B_s -1) \geq  [j/2]  \right\}.$$
Thus,
$$ 0 \leq A_2(\lambda) \leq \exp\left\{\left((m-1)\textup{es} -\gamma(\lambda) \right)j \right\} + \Pi_{1} \left(\sup_{0\leq s\leq j} (B_s -1) \geq  [j/2]  \right).
$$
By the reflection principle and
Markov's inequality, there exists $\delta\in\left(0, (\gamma(\lambda) - (m-1)\textup{es})\wedge \frac{1}{8}\right)$ such that
for all $j\geq 2$,
\begin{align}\label{step_25}
	0\leq A_2(\lambda) &\leq  e^{-j\delta} + 2 \Pi_0(B_j \geq [j/2]) \leq e^{-j \delta}+ 2 \frac{\Pi_0 \left( e^{B_j/2}\right)
}{e^{[j/2]/2}}\nonumber\\& \leq e^{-j \delta}+ 2 e^{j/8-(j/2-1)/2}\leq 5e^{-j\delta}.
\end{align}
Since  $\mathcal{F}_{[j/2]+1}\subset \mathcal{F}_j$ for $j\geq 2$,  we have $A_1(\lambda)\in \mathcal{F}_j$. Thus
\begin{align}\label{step_24}
&	\left|\ln\phi(1,\lambda) -
\mathbb{E}\left(\ln\phi(1,\lambda) \big|\mathcal{F}_{j}  \right)\right|
 = \left|\ln\psi(1,\lambda) - \mathbb{E}
\left(\ln\psi(1,\lambda) \big|\mathcal{F}_{j}  \right)\right|
\nonumber \\
 &=	\left|\ln(A_1(\lambda) + A_2(\lambda)) -\mathbb{E}
\left(\ln(A_1(\lambda) + A_2(\lambda)) \big| \mathcal{F}_{j} \right)
\right|\nonumber\\
  &\leq \left|\ln(A_1(\lambda)) -\mathbb{E}
\left(\ln(A_1(\lambda)) \big| \mathcal{F}_{j} \right)
  \right|
+ \frac{10}{c_2(\lambda)} e^{-j\delta}=\frac{10}{c_2(\lambda)} e^{-j\delta},
\end{align}
where in the last inequality we used the following estimate:
\begin{equation}\label{step_26}
	0\leq
	\ln\left(1+ A_2(\lambda)/ A_1(\lambda)\right)
	\leq A_2(\lambda) /A_1 (\lambda)\leq
 \frac{5}{c_2(\lambda)} e^{-j\delta}.
\end{equation}
Take $\alpha>0$
small so that $\gamma(\lambda) > (m-1)\textup{es} + \alpha$.
Noticing that
$\sup_{h>0} \left(h \exp\left\{-\alpha h\right\} \right)<\infty,$
we have that there is a constant
$c_3$
depending only on $\alpha$ such that
$h \exp\left\{-\alpha h\right\}\leq c_3, \forall h>0$ and thus
\begin{align*}
	0\leq -A_1'(\lambda) & \leq c_3 \gamma'(\lambda)\Pi_{1}
\left(\exp\left\{ \int_0^{H_{0}}\left((m-1)\xi(B_s) -\gamma(\lambda) +\alpha \right)\mathrm{d}s \right\}; \sup_{0\leq s\leq H_0} (B_s -1) < [j/2] \right),\\
	0\leq -A_2'(\lambda)& \leq c_3 \gamma'(\lambda)\Pi_{1}
\left(\exp\left\{ \int_0^{H_{0}}\left((m-1)\xi(B_s) -\gamma(\lambda) +\alpha \right)\mathrm{d}s \right\}; \sup_{0\leq s\leq H_0} (B_s -1) \geq [j/2] \right).
\end{align*}
Define $\gamma(\tilde{\lambda}):= \gamma(\lambda) - \alpha$.
Then we get
$$0\leq -A_i'(\lambda)\leq c_3 \gamma'(\lambda){A}_i(\tilde\lambda),\quad i=1,2.$$
By \eqref{step_25}, there exists $\delta \in\left(0,(\gamma(\tilde\lambda) - (m-1)\textup{es} - \alpha)\wedge\frac{1}{8}\right)$ such that
$\label{domi-A2}0\leq A_2(\tilde\lambda)\leq 5 e^{-j\delta}$, and thus
\begin{equation}\label{domo-A'}
	0\leq -A_1'(\lambda)\leq c_3 \gamma'(\lambda)c_1(\tilde\lambda),\quad 0\leq -A_2'(\lambda)\leq 5 c_3 \gamma'(\lambda) e^{-j\delta}.
\end{equation}
In the display below,
$A_j(\lambda)$ and $A_j'(\lambda)$ will be simply denoted as $A_j$ and $A_j'$ for $j=1,2$.
Hence,
\begin{align}\label{step_27}
& \left|\frac{\phi_\lambda(1,\lambda)}{\phi(1,\lambda)} - \mathbb{E}
\left(\frac{\phi_\lambda(1,\lambda)}{\phi(1,\lambda)} \Big| \mathcal{F}_j\right)
\right|= \left|\frac{\psi_\lambda(1,\lambda)}{\psi(1,\lambda)} - \mathbb{E}
\left(\frac{\psi_\lambda(1,\lambda)}{\psi(1,\lambda)} \Big| \mathcal{F}_j\right)
\right|\nonumber\\
&=\left|-\lambda \frac{A_1'}{A_1+A_2} -\lambda \frac{A_2'}{A_1+A_2} - \mathbb{E}
\left(-\lambda \frac{A_1'}{A_1+A_2} -\lambda \frac{A_2'}{A_1+A_2} \Big|\mathcal{F}_j \right)
 \right|\nonumber\\
& \leq \lambda \left| \frac{A_1'}{A_1+A_2} - \mathbb{E}
\left( \frac{A_1'}{A_1+A_2} \Big|\mathcal{F}_j \right)
 \right| + \lambda \left| \frac{A_2'}{A_1+A_2} \right| + \lambda\mathbb{E}
\left(\left| \frac{A_2'}{A_1+A_2} \right| \bigg|\mathcal{F}_j \right)
\nonumber\\
 & \leq \lambda \frac{|A_1'|}{A_1} \left\{ 1-\frac{A_1}{A_1+A_2} + \mathbb{E}
 \left( 1-\frac{A_1}{A_1+A_2} \Big|\mathcal{F}_j \right)
 \right\}	+ 10\lambda  c_3 \gamma'(\lambda) \frac{e^{-j\delta}}
{c_2(\lambda)}\nonumber
	\\ &\leq  10\lambda c_3 \gamma'(\lambda)\frac{{c}_1(\tilde\lambda)}{c_2(\lambda)}  \frac{e^{-j\delta}}{c_2(\lambda)}  +10 \lambda c_3 \gamma'(\lambda)\frac{e^{-j\delta}}
{c_2(\lambda)}
= :c_4
e^{-j\delta},
\end{align}
here in the second inequality above we used the fact that
$A_1$ and $A_1'$ are $\mathcal{F}_j$-measurable, and in the last inequality we used \eqref{domi-A1}, \eqref{domo-A'}, and the fact that $A_2\geq 0$.

Now \eqref{step_24} and \eqref{step_27}
imply that there exists a non-random constant
$c_5$
  independent of $j$ such that for all $j\geq 2$,
\begin{equation}\label{step_40}
\left|\mathbb{E} \left(U_0 \big| \mathcal{F}_{j} \right)
- U_0 \right| \leq c_5 e^{-j\delta},
\end{equation}
which implies that $\sum_{j=1}^\infty \left( \mathbb{E}
\left(U_0 \big| \mathcal{F}_{j} \right)
- U_0 \right) $ converges in $L^2(\Omega, \mathcal{F}, \mathbb{P}).$

Next, we estimate $\left|\mathbb{E} \left(U_j \big| \mathcal{F}_0 \right) \right|$ for $j\geq 1$.
Note that when $\lambda >\rho$, $U_j \in \mathcal{F}^j$ by \eqref{Psi_2}. Using the $\zeta$-mixing condition and the fact that $\mathbb{E} U_j = 0$, we have
\begin{equation}\label{step_41}
	\left|\mathbb{E}(U_j | \mathcal{F}_0) \right|\leq \mathbb{E}\left(|U_j|\right)
\zeta(j)\leq \Vert U_0 \Vert_\infty \zeta(j).
\end{equation}
Here and later in this proof, we use $\|\cdot\|_\infty$ to denote the supremum over $\omega$.
By the triangle inequality,
$$
\sqrt{\mathbb{E}\left(\sum_{j=m}^n \mathbb{E}(U_j | \mathcal{F}_0)
\right)^2} \leq \sum_{k=m}^n \sqrt{\mathbb{E}
\left(\left( \mathbb{E}(U_j | \mathcal{F}_0)\right)^2\right)}
\leq \sum_{j=m}^n \Vert U_0\Vert_\infty \zeta(j),
$$
which implies that $\sum_{j=1}^\infty  \mathbb{E}
(U_j | \mathcal{F}_0)
$ converges in $L^2(\Omega,\mathcal{F}, \mathbb{P}).$
Applying Lemma \ref{CLT} and note that
$\left|\mathbb{E}(U_0U_j)\right| \leq
\mathbb{E}\left(|U_0| \left|\mathbb{E}\left(U_j|\mathcal{F}_0\right)\right|\right)
\leq
\Vert U_0\Vert_\infty^2
\zeta(j)$,
\begin{equation}\label{Variation}
\sigma_\lambda^2 := \lim_{N\to\infty} \frac{1}{N} \mathbb{E}
\left( \sum_{i=0}^{N-1} U_i \right)^2
=\mathbb{E}(U_0^2)+ 2\sum_{j=1}^\infty \mathbb{E}(U_0 U_j)
\leq \Vert U_0\Vert_\infty^2 + 2\Vert U_0 \Vert_\infty^2\sum_{j=1}^\infty \zeta(j) < \infty.
\end{equation}
Thus
$$\frac{\ln \phi(n,\lambda) - \lambda \frac{\phi_\lambda(n,\lambda)}{\phi(n,\lambda)}}{\sqrt{n}}=\frac{\ln \psi(n,\lambda) - \lambda \frac{\psi_\lambda(n,\lambda)}{\psi(n,\lambda)}}{\sqrt{n}} \stackrel{\mathrm{d}}{\Longrightarrow} \mathcal{N}(0,\sigma_\lambda^2),\qquad \textup{as } n\to+\infty,$$
which concludes the proof of the first part.

If $\sigma_\lambda^2> 0$, set $R_x = \ln \psi(x,\lambda) - \lambda \frac{\psi_\lambda(x,\lambda)}{\psi(x,\lambda)}$. For any $K>0$, by the central limit theorem for $R_n$ we have
$$
\mathbb{P}\left(\limsup_{x\to+\infty} R_x /\sqrt{x}>K \right) = \lim_{x\to+\infty} \mathbb{P}\left(\sup_{y\geq x}R_y/\sqrt{y} >K \right)\geq \lim_{x\to\infty} \mathbb{P}(R_x/\sqrt{x} > K) > 0.
$$
Since $\{\limsup_{x\to+\infty}R_x/{\sqrt{x}}> K \}$ is an invariant set,  we have
$$
\mathbb{P}\left(\limsup_{x\to\infty}{R_x}/{\sqrt{x}} =+\infty \right) =1.
$$
When $x \to -\infty$,
a similar argument also shows that
$\mathbb{P}(\limsup_{x\to-\infty} R_x/\sqrt{|x|} = +\infty ) = 1$.

For any $M>0$, we modify the definition of $U_i$ in \eqref{U_i}  by
$$	
U_i= \ln \psi\left(M(i+1),\lambda,\omega \right) - \ln \psi\left( Mi,\lambda,\omega \right) - \lambda
\left(\frac{\psi_\lambda\left(M(i+1),\lambda,\omega\right)}
{\psi\left(M(i+1),\lambda,\omega\right)}-\frac{\psi_\lambda\left(Mi,\lambda,\omega\right)}{\psi\left(Mi,\lambda,\omega\right)} \right).
$$
By Lemma \ref{CLT}, the process
$$\eta^n_t:=\frac{\ln \phi(nt,\lambda) - \lambda \frac{\phi_\lambda(nt,\lambda)}{\phi(nt,\lambda)}}{\sigma_\lambda\sqrt{n}},\quad t\in [0,M],$$
converges weakly to a standard Brownian motion on $[0, M]$.
Since $M$ is arbitrary, by \cite[Lemma 3, p.173]{Billingsley},  we get the weak convergence of $\{\eta^n_t,t\geq 0\}$ to a
Brownian motion on $[0,\infty)$.

If $\sigma_\lambda^2 = 0$, then by \cite[\S 5.4]{HH} or \cite[(19.8), p.198]{Billingsley}, we may write the increments $\{U_i, i\in\mathbb{Z}\}$ of the stationary ergodic process $\{\ln \psi(n,\lambda) - \lambda \frac{\psi_\lambda(n,\lambda)}{\psi(n,\lambda)}, n\in \mathbb{Z}\}$ as:
 $$U_i= Y_i + (T^*)^i Z_0 - (T^*)^{i+1}Z_0,\quad i\in \mathbb{Z},$$
 here $\{Y_i, i\in\mathbb{Z}\}$ is a stationary ergodic sequence of  martingale differences with $Y_i\in \mathcal{F}_i$ and  $
 \mathbb{E}\left(Y_{i+1}|\mathcal{F}_i\right)
 =0$,
  $Z_0$ is a random variable and $T^*$ is the unitary operator associate with $\theta_1$, that is, $T^*Z_0(\omega) = Z_0(\theta_1 \omega).$
  In fact, we may take
\begin{align}\label{step_47}
	Z_0 = \sum_{k=0}^\infty \mathbb{E}\left(U_{k}\big| \mathcal{F}_{-1} \right) - \sum_{k=1}^\infty \left(U_{-k} -\mathbb{E}\left(U_{-k} \big| \mathcal{F}_{-1} \right)\right),
\end{align}
$Y_i=(T^*)^i Y_0$ with
$$
Y_0=\sum^\infty_{l=-\infty}
\left(\mathbb{E}(U_l \big| \mathcal{F}_0)-\mathbb{E}(U_l \big| \mathcal{F}_{-1})\right),
$$
see \cite[(5.17) and (5.18), p.137]{HH}.
From the proof of \cite[Theorem 5.5]{HH}, we know that $\mathbb{E}[Y_0^2] = \sigma_\lambda^2=0$ (see \cite[the fourth paragraph, on p.142]{HH}), and thus $Y_i=0$, $\mathbb{P}$-a.s. for any $i\in\mathbb{Z}$.
Note that
\begin{equation}
\|U_{-k}-\mathbb{E}(U_{-k}\big|\mathcal{F}_{-1})\|_\infty=\|(T^*)^{-k}\left(U_0-\mathbb{E}(U_0\big|\mathcal{F}_{k-1}) \right)\|_\infty=
\|U_0-\mathbb{E}(U_0\big|\mathcal{F}_{k-1})\|_\infty.
\end{equation}
Thus by \eqref{step_40},
\begin{equation}\label{tran-U}
\|U_{-k}-\mathbb{E}(U_{-k}\big|\mathcal{F}_{-1})\|_\infty\leq
c_5 e^{-(k-1)\delta},\quad k\geq 3.
\end{equation}
Now using \eqref{tran-U}, \eqref{step_41} and \eqref{step_47}, we have that
$$\Vert Z_0 \Vert_\infty \leq \Vert U_0 \Vert_\infty \sum_{k=0}^\infty\zeta(k+1) + c_5 \sum_{k=2}^\infty  e^{-k\delta} + 4 \Vert U_0 \Vert_\infty < \infty.$$
Therefore, for all $n\in\mathbb{Z}^+$,
$$\left|\ln \phi(n,\lambda) - \lambda \frac{\phi_\lambda(n,\lambda)}{\phi(n,\lambda)}\right| = \left| \sum_{i=0}^{n-1}U_i\right|= \left|Z_0 - (T^*)^n Z_0 \right|\leq 2\Vert Z_0 \Vert_\infty,$$
which implies that
$$\sup_{x\in\R} \left\Vert \ln \phi(x,\lambda) - \lambda \frac{\phi_\lambda(x,\lambda)}{\phi(x,\lambda)} \right\Vert_\infty \leq 2\Vert Z_0 \Vert_\infty +
C_1(\lambda)+\lambda K_1(\lambda)<\infty.
$$
\hfill$\Box$

\begin{remark}\label{CLT2} Suppose that $|\lambda| > \rho$ and $\gamma(\lambda)> (m-1)\textup{es}$.
	By the proof of Lemma \ref{Engenfunction-Property2}, and using \eqref{step_24} and \eqref{step_27},
the limits
	\begin{equation}\label{Another_var}
(\sigma_\lambda')^2:= \lim_{x\to +\infty}\frac{1}{|x|} \mathbb{E} \left(\ln \psi(x,\lambda) + \lambda x\right)^2,
\qquad  (\sigma_\lambda'')^2:= \lim_{x\to +\infty} \frac{1}{|x|}\mathbb{E}\left( \frac{\psi_\lambda(x,\lambda)}{\psi(x,\lambda)} +  x\right)^2,
	\end{equation}
exist under $\mathbb{P}$, and as $|x|\to \infty$,
$$
\frac{\ln \phi(x,\lambda)}{\sqrt{|x|}}=\frac{\ln \psi(x,\lambda)+ \lambda x}{\sqrt{|x|}} \stackrel{\mathrm{d}}{\Longrightarrow} \mathcal{N}(0,(\sigma_\lambda')^2),\qquad \frac{ \phi_\lambda(x,\lambda)}{\phi(x,\lambda)\sqrt{|x|}}=\frac{ \frac{\psi_\lambda(x,\lambda)}{\psi(x,\lambda)}+  x}{\sqrt{|x|}} \stackrel{\mathrm{d}}{\Longrightarrow} \mathcal{N}(0,(\sigma_\lambda'')^2).
$$
If $(\sigma_\lambda')^2>0$, the process
	$$\frac{\ln \phi(nt,\lambda) }{\sigma_\lambda'\sqrt{n}},\quad t\in [0,\infty)$$
	converge weakly to a standard Brownian motion on $[0,\infty)$ in the Skorohod topology.
If $(\sigma_\lambda')^2=0$,
then $\sup_{x\in\R}\Vert \phi(x,\lambda)\Vert_\infty<\infty.$ The same result also holds for $\frac{ \phi_\lambda(x,\lambda)}{\phi(x,\lambda)}$.
\end{remark}

\subsection{Limit theorems for the spine}
Define
	\begin{equation}\label{define-T}T_k: = \inf\{t> 0 :\Xi_t - \Xi_0 = k \},\qquad k\in \mathbb{Z}.
	\end{equation}
\begin{prop}\label{lemma 4.1}
If {\bf(H1)} holds, then for any $|\lambda| > \rho$ and  $\omega \in \Omega_1$,
\begin{equation}\label{spine-limit}\lim_{t\to\infty}\frac{\Xi_t}{t} \to -\gamma'(\lambda),\ \widetilde{\mathrm{P}}_x^{\xi, \lambda}\textup{-a.s.}
\end{equation}
If we further assume that $\gamma(\lambda) > (m-1)\textup{es}$, then there exist
$\Omega_\lambda $ with
$\P(\Omega_\lambda)=1$ and a constant $\Sigma_\lambda^2> 0$ such that for $\omega \in \Omega_\lambda$,  when $\lambda > \rho$,
	$$
 \frac{T_{-k} - \widetilde{\mathrm{E}}_x^{\xi, \lambda} T_{-k}}{\sqrt{k}} \stackrel{\mathrm{d}}{\Rightarrow}
\mathcal{N}\left(0,\Sigma_\lambda^2\right) \quad \mbox{ under }
\widetilde{\mathrm{P}}_x^{\xi,\lambda};
$$
and when $\lambda < -\rho$,
	$$ \frac{T_{k} - \widetilde{\mathrm{E}}_x^{\xi, \lambda} T_{k}}{\sqrt{k}} \stackrel{\mathrm{d}}{\Rightarrow}
\mathcal{N}\left(0,\Sigma_\lambda^2\right) \quad \mbox{ under }
\widetilde{\mathrm{P}}_x^{\xi,\lambda}.
$$
\end{prop}
\textbf{Proof} The idea of the proof is from \cite[Lemma 2.6 and Corollary 2.7]{LTZ}.
(i)  Fix $\lambda, x$ and $\omega$. For any $\varepsilon > 0$, by the definition of $\gamma'(\lambda)$, there exists
$\delta: = \delta(\varepsilon, \lambda,\rho) \in (0, |\lambda| -\rho)$
such that  for any $0< |\eta| < \delta$,
$$ \left|\frac{\gamma(\lambda + \eta) - \gamma(\lambda)}{\eta} - \gamma'(\lambda) \right| < \varepsilon.$$
Our first goal is to show that
\begin{equation}\label{step_8}
	\lim_{t\to\infty} \frac{1}{t}\ln \widetilde{\mathrm{E}}_x^{\xi, \lambda}
\left(e^{\eta \Xi_t} \right)
= \gamma(\lambda -\eta)-\gamma(\lambda).
\end{equation}
By \eqref{Change-of-measure} and \eqref{=distribution},
\begin{align}\label{step_4}
& \frac{1}{t}\ln \widetilde{\mathrm{E}}_x^{\xi, \lambda}\left(e^{\eta \Xi_t} \right)=\frac{1}{t} \ln \Pi_x
\left(\exp\left\{\eta B_t -\gamma(\lambda)t + \int_0^t (m-1)\xi(B_s) \mathrm{d} s - \lambda (B_t -x)\right\} \frac{\phi(B_t,\lambda)}{\phi(x,\lambda)}  \right)
\nonumber \\
& = -\gamma(\lambda) - \frac{- x\lambda +\ln \phi(x,\lambda)}{t} + \frac{1}{t} \ln \Pi_x
\left(\exp\left\{(\eta - \lambda) B_t+ \int_0^t (m-1)\xi(B_s) \mathrm{d}s\right\}\phi(B_t, \lambda)\right)
\nonumber\\
&= -\gamma(\lambda) - \frac{-x\lambda +\ln \phi(x,\lambda)}{t}\nonumber \\
&\quad + \frac{1}{t} \ln \Pi_x
\left(\exp\left\{(\eta - \lambda) B_t+ \int_0^t (m-1)\xi(B_s) \mathrm{d}s\right\}\phi(B_t, \lambda-\eta)\frac{\phi(B_t,\lambda)}{\phi(B_t,\lambda - \eta)} \right).
\end{align}
Now fix $\eta$. Since $\ln \phi(y, \lambda)/y$ and
$\ln\phi(y, \lambda -\eta) /y$ converge to $0$ as $|y|\to\infty$, for any $\epsilon > 0$, there exists an $M_1 = M_1(\epsilon,\delta,\lambda)$ such that for all $y\in \mathbb{R},$
$$M_1^{-1}e^{-\epsilon |y|} \leq  \phi(y, \lambda) \leq M_1e^{\epsilon |y|},\quad M_1^{-1}e^{-\epsilon |
	y|} \leq \phi(y, \lambda-\eta) \leq M_1e^{\epsilon |y|}.$$
Thus, for all $y
\in \mathbb{R}$,
$$M_1^{-2}e^{-2\epsilon|y| }\leq \frac{\phi (y, \lambda)}{\phi(y, \lambda - \eta)} \leq M_1^2 e^{2\epsilon |y|}.$$
Applying the above to \eqref{step_4}, we get
\begin{align}\label{step_5}
&	\frac{1}{t}\ln \widetilde{\mathrm{E}}_x^{\xi, \lambda}
\left(e^{\eta \Xi_t} \right)
\nonumber \\
& \geq  -\gamma (\lambda) + \frac{1}{t} \ln \Pi_x
  \left(\exp\left\{(\eta - \lambda) B_t- 2\epsilon |B_t|+ \int_0^t (m-1)\xi(B_s) \mathrm{d}s\right\}\phi(B_t, \lambda-\eta) \right)
+ O\left(\frac{1}{t}\right) \nonumber\\
& =  -\gamma(\lambda) + \gamma(\lambda - \eta) + \frac{1}{t}\ln \widetilde{\mathrm{E}}_x^{\xi,\lambda -\eta}
\left(e^{-2\epsilon |\Xi_t|} \right)
+  O\left(\frac{1}{t}\right),
\end{align}
here  $O(1/t)=
\left(\ln (M_1^{-2})- \left(-x\lambda +\ln \phi(x,\lambda){-\ln \phi(x,\lambda-\eta)}\right)\right)
/t.$
Similarly
\begin{equation}\label{step_5'}
\frac{1}{t}\ln \widetilde{\mathrm{E}}_x^{\xi, \lambda}
 \left(e^{\eta \Xi_t} \right)
\leq -\gamma(\lambda) + \gamma(\lambda - \eta) + \frac{1}{t}\ln \widetilde{\mathrm{E}}_x^{\xi, \lambda -\eta}
 \left(e^{2\epsilon |\Xi_t|} \right)
+  O\left(\frac{1}{t}\right)
\end{equation}
with $O(1/t)= \left(\ln (M_1^{2})- \left(-x\lambda +\ln \phi(x,\lambda){-\ln \phi(x,\lambda-\eta)}\right)\right)
/t.$
Note  under $\widetilde{\mathrm{P}}_x^{\xi,\lambda-\eta}$, $\Xi_t$  satisfies
$$
\Xi_t = W_t + \int_0^t \frac{\psi_x(\Xi_s,\lambda-\eta)}{\psi(\Xi_s,\lambda-\eta)}\mathrm{d} s.$$
Since $\phi(\cdot,\lambda-\eta)\in \mathcal{A}_\infty,$ $\phi_x(\cdot,\lambda-\eta)/\phi(\cdot, \lambda-\eta)\in L^\infty(\R)$, by \eqref{Psi_and_Phi}, there exists a constant $M_{\lambda-\eta}$, depending only on $\lambda-\eta$, such that for all $y\in\R$,
\begin{equation}\label{domi-frac-psi}\left|\frac{\psi_x(y,\lambda-\eta)}{\psi(y,\lambda-\eta)} \right|
= \left|-(\lambda-\eta)
 +\frac{\phi_x(y,\lambda-\eta)}{\phi(y,\lambda-\eta)} \right|
\leq M_{\lambda-\eta}.\end{equation}
Therefore, under $\widetilde{\mathrm{P}}_x^{\xi,\lambda - \eta}$, $|\Xi_t| \leq |W_t-x| + |x| + M_{\lambda-\eta}t.$ Since $\left(W_t-x, \widetilde{\mathrm{P}}_x^{\xi,\lambda - \eta}\right) \stackrel{\mathrm{d}}{=} \left(B_t, \Pi_0\right)$,
by \eqref{step_5},  we have
\begin{align}\label{liminf-Xi}
\liminf_{t\to\infty} \frac{1}{t}\ln \widetilde{\mathrm{E}}_x^{\xi, \lambda}
\left(e^{\eta \Xi_t} \right)
& \geq -\gamma(\lambda) + \gamma(\lambda - \eta) + \liminf_{t\to\infty} \frac{1}{t}\ln \widetilde{\mathrm{E}}_x^{\xi, \lambda -\eta}
\left(e^{-2\epsilon |\Xi_t|} \right)
\nonumber\\
& \geq -\gamma(\lambda) + \gamma(\lambda - \eta)- 2\epsilon M_{\lambda-\eta}
+ \liminf_{t\to\infty} \frac{1}{t}\ln \Pi_0
\left(e^{-2\epsilon |B_t|} \right).
\end{align}
For fixed $\epsilon$, when $t$ is large enough,
\begin{equation}\label{step_7}
\Pi_0\left(e^{-2\epsilon|B_t|} \right)
= \sqrt{\frac{2}{\pi}}\int_0^\infty e^{-2\epsilon\sqrt{t}y - y^2/2}\mathrm{d}y  = e^{-2\epsilon^2 t}\sqrt{\frac{2}{\pi}} \int_{2\epsilon \sqrt{t}}^\infty e^{-y^2 /2}\mathrm{d}y \geq  e^{-8\epsilon^2t},
\end{equation}
where in the last inequality we used the fact that
$$\left(\frac{1}{z}-\frac{1}{z^3}\right)e^{-z^2/2} \leq \int_z^\infty e^{-y^2/2}\mathrm{d}y,\quad z > 0$$
and the inequality
$\sqrt{2/\pi}(1/z - 1/z^3) \geq e^{-z^2}$ for $z>0$ large.
Combining  \eqref{liminf-Xi} and \eqref{step_7}, we get
\begin{align}\label{liminf-Xi2}
\liminf_{t\to\infty} \frac{1}{t}\ln \widetilde{\mathrm{E}}_x^{\xi, \lambda} \left(e^{\eta \Xi_t} \right)\geq  -\gamma(\lambda) + \gamma(\lambda - \eta)- 2\epsilon M_{\lambda-\eta} - 8\epsilon^2.
\end{align}
Similarly,
\begin{equation}\label{step_6}
\Pi_0\left(e^{2\epsilon|B_t|} \right)
= \sqrt{\frac{2}{\pi}}\int_0^\infty e^{2\epsilon\sqrt{t}y - y^2/2}\mathrm{d}y  \leq \sqrt{\frac{2}{\pi}}\int_{-\infty}^\infty e^{2\epsilon\sqrt{t}y - y^2/2}\mathrm{d}y = \sqrt{\frac{2}{\pi}} e^{2\epsilon^2 t}.
\end{equation}
Using \eqref{step_5'} and \eqref{step_6}, we get
\begin{equation}\label{limsup-Xi2}
\limsup_{t\to\infty} \frac{1}{t}\ln \widetilde{\mathrm{E}}_x^{\xi, \lambda} \left(e^{\eta \Xi_t} \right)\leq -\gamma(\lambda) + \gamma(\lambda - \eta) + 2\epsilon M_{\lambda-\eta} + 2\epsilon^2.
\end{equation}
Since $M_{\lambda-\eta}$
 is independent of $\epsilon$, let $\epsilon \to 0+$ in \eqref{liminf-Xi2} and \eqref{limsup-Xi2}, we get \eqref{step_8}.

 Our second goal is to show that, for any $\epsilon_1 > 0$,
there exist $M_3 = M_3(\lambda ,\epsilon_1)>0$ and $ \delta_1 = \delta_1(\epsilon_1,\lambda)>0$
 such that for all $t$,
\begin{equation}\label{step_9}
\widetilde{\mathrm{P}}_x^{\xi, \lambda}\left(|\Xi_t + \gamma'(\lambda)t | > \epsilon_1 t\right) \leq M_3 e^{-\delta_1 t}.
\end{equation}
For any $\epsilon_1>0$, let $\varepsilon =  \epsilon_1/4$ and
let $\delta$ satisfy the condition at the beginning of the proof.
By \eqref{step_8}, for fixed $\eta \in (0, \delta)$, there exists $M_4 = M_4 (\eta, \lambda,\epsilon_1)$ such that for all $t$,
\begin{align*}
\ln \widetilde{\mathrm{E}}_x^{\xi, \lambda}  \left(e^{-\eta \Xi_t} \right)
& \leq  M_4+ \left(\gamma(\lambda + \eta) -\gamma(\lambda) +\epsilon_1 \eta/4\right)t \leq M_4 + \gamma'(\lambda)\eta t + \epsilon_1 \eta t/2,\\
\ln \widetilde{\mathrm{E}}_x^{\xi, \lambda}
\left(e^{\eta \Xi_t} \right)
& \leq  M_4+ \left(\gamma(\lambda - \eta) -\gamma(\lambda) +\epsilon_1 \eta/4\right)t \leq M_4 -\gamma'(\lambda)\eta t +\epsilon_1\eta t/2.
	\end{align*}	
Thus, by Markov's inequality,
\begin{align*}
	\widetilde{\mathrm{P}}_x^{\xi, \lambda}\left(\Xi_t < \left(- \gamma'(\lambda) - \epsilon_1\right) t\right)&  \leq \exp\left\{ \ln \widetilde{\mathrm{E}}_x^{\xi, \lambda}
	 \left(e^{-\eta \Xi_t} \right)
	 -\left(\gamma'(\lambda)+\epsilon_1  \right)\eta t \right\}
	\leq \exp\left\{M_4 -\eta\epsilon_1 t/2  \right\}, \nonumber\\
	\widetilde{\mathrm{P}}_x^{\xi, \lambda}\left(\Xi_t >\left( -\gamma'(\lambda) + \epsilon_1\right) t\right)&  \leq \exp\left\{ \ln \widetilde{\mathrm{E}}_x^{\xi, \lambda}
	\left(e^{\eta \Xi_t} \right)
	 +\left(\gamma'(\lambda)-\epsilon_1  \right)\eta t \right\}	\leq \exp\left\{M_4 -\eta\epsilon_1 t/2  \right\}.
\end{align*}
Taking $M_3 = 2e^{M_4}$ and $\delta_1 = \eta \epsilon_1/2$,  we obtain \eqref{step_9}.

Finally, for any $\epsilon_1 > 0$ and $n > 2\left(M_\lambda
 + |\gamma'(\lambda)|\right)/\epsilon_1$
where $M_\lambda$ is defined in \eqref{domi-frac-psi},
   by Markov's inequality, we have
\begin{align}\label{step_51}
	&\widetilde{\mathrm{P}}_x^{\xi, \lambda}\left( \sup_{n \leq t \leq n+1} |\Xi_t - \Xi_n + \gamma'(\lambda)(t-n)|> \epsilon_1 n \right)\nonumber\\
& = \widetilde{\mathrm{P}}_x^{\xi, \lambda}\left( \sup_{n \leq t \leq n+1} \left|W_t - W_n + \int_n^t \psi_x(\Xi_s,\lambda)/\psi(\Xi_s,\lambda) \mathrm{d}s +\gamma'(\lambda)(t-n)\right| > \epsilon_1 n \right)\nonumber\\
 & \leq \Pi_0 \left(\sup_{t\in[0,1]} |B_t| > \epsilon_1 n -
M_\lambda-|\gamma'(\lambda)|  \right) \leq \Pi_0 \left(\sup_{t\in[0,1]} |B_t| > \epsilon_1 n /2 \right) \nonumber\\
& = 2 \Pi_0(B_1 > \epsilon_1 n/2) \leq
2\frac{\Pi_0 \left(e^{\epsilon_1 n B_1/4}\right)
}{e^{\epsilon_1^2 n^2/8}} = 2e^{-\epsilon_1^2 n^2/16},
\end{align}
where in the first inequality we used \eqref{domi-frac-psi}.
Together with \eqref{step_9} and \eqref{step_51}, we conclude that
\begin{align*}
	&\sum_{n=1}^\infty \widetilde{\mathrm{P}}_x^{\xi, \lambda}\left( \sup_{n \leq t \leq n+1} |\Xi_t
	+ \gamma'(\lambda)t|>2 \epsilon_1 n \right)  \\
 & \leq  \sum_{n=1}^\infty \widetilde{\mathrm{P}}_x^{\xi, \lambda}\left(  |\Xi_n
	 +  \gamma'(\lambda)n|> \epsilon_1 n \right) +
		\sum_{n=1}^\infty \widetilde{\mathrm{P}}_x^{\xi, \lambda}\left( \sup_{n \leq t \leq n+1} |\Xi_t - \Xi_n
	+ \gamma'(\lambda)(t-n)|> \epsilon_1 n \right)< \infty.
\end{align*}
Now \eqref{spine-limit} follows from the  Borel-Cantelli lemma.

(ii)
Assume $\lambda > \rho$.
For $i \in \mathbb{N}$, define
\begin{align*}
	F_{2}(x,- i, \omega) & := \widetilde{\mathrm{E}}_x^{\xi, \lambda}
	\left(T_{-i} - T_{-i+1} - \widetilde{\mathrm{E}}_x^{\xi, \lambda} )T_{-i} - T_{-i+1})\right)^2,\\
	F_{4}(x, -i, \omega) & := \widetilde{\mathrm{E}}_x^{\xi, \lambda} \left(T_{-i} - T_{-i+1} - \widetilde{\mathrm{E}}_x^{\xi, \lambda} (T_{-i} - T_{-i+1})\right)^4.
\end{align*}
We claim that
$\mathbb{P}$-almost surely for all $i\geq 1,\ F_2(x,i,\omega) = F_2(x,1, \theta_{-i+1}\omega)$ and $F_4(x,i,\omega) = F_4(x,1, \theta_{-i+1}\omega)$.
By the binomial theorem, to prove the claim above, it suffices to show that,
for any $m=1, \dots, 4$ and $i\in \mathbb{N}$,
$\mathrm{E}_x (T_{-i} - T_{-i+1})^m
$ is stationary.
By \eqref{Change-of-measure}
 and the strong Markov property of $\Xi_t$,
\begin{align*}
&\widetilde{\mathrm{E}}_x^{\xi, \lambda}
\left(T_{-i} - T_{-i+1}\right)^m =
\widetilde{\mathrm{E}}_{-i+1+x}^{\xi, \lambda}
\left({T^m_{-1}}\right) \\
& =\Pi_{-i+1+x} \left(\exp\left\{-\gamma(\lambda) H_{-i+x} + \int_0^{H_{-i+x}} (m-1)\xi(B_s,\omega)\mathrm{d}s \right\} \frac{\psi(-i+x,\lambda,\omega)}{\psi(-i+1+x,\lambda,\omega)}H_{-i+x}^m\right)\\
 &= \Pi_{x}
 \left(\exp\left\{-\gamma(\lambda) H_{-1+x} + \int_0^{H_{-1+x}} (m-1)\xi(B_s,\theta_{-i+1}\omega)\mathrm{d}s \right\} \frac{\psi(-1+x,\lambda,\theta_{-i+1}\omega)}{\psi(x,\lambda,\theta_{-i+1}\omega)}H_{-1+x}^m\right),
\end{align*}
here in the last equality, we used the fact that
$$
\frac{\psi(-i+x,\lambda,\omega)}{\psi(-i+1+x,\lambda,\omega)}=e^\lambda\phi(-1,\lambda,\theta_{-i+1+x}\omega)=
 \frac{\psi(-1+x,\lambda,\theta_{-i+1}\omega)}{\psi(x,\lambda,\theta_{-i+1}\omega)},
 $$
 which is true by \eqref{step_3-1'}.
Thus $\mathrm{E}_x (T_{-i} - T_{-i+1})^m
$ is stationary.
Note that we have assumed $\lambda >\rho $ and $\gamma(\lambda) > (m-1)\textup{es}$, so by Jensen's inequality,
the trivial inequality
$\mathbb{E}(X-\mathbb{E}X)^4 \leq 16 \mathbb{E}(X^4)$
and \eqref{Change-of-measure},
\begin{align*}
	0 &< \left(\mathbb{E} F_2(x,-1,\omega)\right)^2  \leq {\mathbb{E} F_2^2(x,-1,\omega)} \leq {\mathbb{E} F_4(x,-1,\omega)} \leq
	16\mathbb{E} \widetilde{\mathrm{E}}_x^\lambda T_{-1}^4 \\ &=
	16\mathbb{E}\left\{ \frac{\Pi_x\left(\exp\left\{ \int_0^{H_{x-1}}\left((m-1)\xi(B_s, \omega ) -\gamma(\lambda) \right)\mathrm{d}s \right\} H_{x-1}^4\right)}{\Pi_x\left(\exp\left\{ \int_0^{H_{x-1}}\left((m-1)\xi(B_s, \omega ) -\gamma(\lambda) \right)\mathrm{d}s \right\} \right)}\right\}\\
	&\leq 16\frac{\Pi_x\left(\exp\left\{- \left(\gamma(\lambda) - (m-1)\textup{es} \right)H_{x-1}\right\} H_{x-1}^4\right)}{\Pi_x\left(\exp\left\{- \left(\gamma(\lambda) - (m-1)\textup{ei} \right)H_{x-1}\right\} \right)} < \infty.
\end{align*}
Applying Birkhoff's ergodic theorem, we get that $\mathbb{P}$-almost surely, both $\frac{1}{k}\sum_{i=1}^k F_2(x,-i,\omega)$ and  $\frac{1}{k}\sum_{i=1}^k F_4(x,-i,\omega)$ converge as $k\to\infty$. Let
\begin{align}\label{Omega_Lambda}
		\Omega_\lambda:= \Omega_1 & \bigcap \left\{\omega: \lim_{k\to+\infty} \frac{\sum_{i=1}^k F_2(x,-i,\omega) }{k}= \mathbb{E}F_2(x,-1,\omega)  \right\}\nonumber\\
		& \bigcap \left\{\omega: \lim_{k\to+\infty} \frac{\sum_{i=1}^k F_4(x,-i,\omega) }{k}= \mathbb{E}F_4(x,-1,\omega)  \right\}.
\end{align}
Then $\mathbb{P}(\Omega_\lambda) = 1$ and for $\omega \in \Omega_\lambda, $
$$ \lim_{k\to\infty} \frac{1}{\left(\sum_{i=1}^k F_2(x,-i,\omega) \right)^2} \sum_{i=1}^k F_4(x,-i,\omega) =0,$$
which implies by Lyapunov's theorem that
$(T_{-k} - \widetilde{\mathrm{E}}_x^{\xi, \lambda} T_{-k})/\sqrt{k}$ converges in
$\widetilde{\mathrm{P}}_x^{\xi, \lambda}$-
distribution to $\mathcal{N}\left(0, \mathbb{E}F_2(x,-1,\omega)\right)$.
Setting $\Sigma_\lambda^2 : = \mathbb{E}F_2(x,-1,\omega)$, we get the desired conclusion.
\hfill$\Box$

\begin{remark}\label{remark}
By the strong Markov property of $\Xi$, $\{T_{-i}-T_{-i+1}, i\ge 1\}$
are  independent
(see \eqref{strong-Markov-3} below).
For any $M>0$ and any subsequence $\{T_{-n_k}\}$ of $\{T_{-k} \}$,
$\{\liminf_{k\to\infty} (T_{-n_k} - \widetilde{\mathrm{E}}_x^{\xi, \lambda} T_{-n_k})/ \sqrt{n_k} <-M \}$ belongs to the tail $\sigma$-field of
$\{T_{-i} - T_{-i+1},i\geq 1 \}$.
By the central limit theorem, we have
\begin{align*}
\widetilde{\mathrm{P}}_x^{\xi, \lambda}
\left(\liminf_{t\to\infty}\frac{T_{-n_m} - \widetilde{\mathrm{E}}_x^{\xi, \lambda} T_{-n_m}}{\sqrt{n_m}} < -M \right)
=&\lim_{k\to\infty} \widetilde{\mathrm{P}}_x^{\xi, \lambda}
\left( \inf_{m\geq k} \frac{T_{-n_m} - \widetilde{\mathrm{E}}_x^{\xi, \lambda} T_{-n_m}}{\sqrt{n_m}} < -M \right)\\
 \geq &\lim_{k\to\infty}\widetilde{\mathrm{P}}_x^{\xi, \lambda}
 \left( \frac{T_{-n_k} - \widetilde{\mathrm{E}}_x^{\xi, \lambda} T_{-n_k}}{\sqrt{n_k}} < -M \right) >0.
\end{align*}
Thus by Kolmogorov's 0-1 law,
$$
\widetilde{\mathrm{P}}_x^{\xi, \lambda}
\left( \liminf_{k\to\infty} \frac{T_{-n_k} - \widetilde{\mathrm{E}}_x^{\xi, \lambda} T_{-n_k}}{\sqrt{n_k}} <-M \right) = 1.
$$
Letting $M\to +\infty$, we get
$$
\widetilde{\mathrm{P}}_x^{\xi, \lambda}
\left( \liminf_{k\to\infty} \frac{T_{-n_k} - \widetilde{\mathrm{E}}_x^{\xi, \lambda} T_{-n_k}}{\sqrt{n_k}} = -\infty \right) = 1.
$$
\end{remark}

\section{Proof of the main results}

\subsection{Proof of Theorem \ref{Additive martingale}}

\textbf{Proof of Theorem \ref{Additive martingale}:} We only consider the case $\lambda > \rho.$  Define
$$\widetilde{W}_t(\lambda) :=e^{-\gamma(\lambda)t}\sum_{\nu \in N(t)} e^{-\lambda \widetilde{X}_t^\nu}\phi(\widetilde{X}_t^\nu,\lambda) = e^{-\gamma(\lambda)t}\sum_{\nu \in N(t)} \psi(\widetilde{X}_t^\nu,\lambda)$$
 and $\widetilde{W}_\infty (\lambda) := \limsup_{t\to\infty} \widetilde{W}_t(\lambda)$.  By \cite[Theorem 4.3.5. p.227]{Durrett},
\begin{align*}
		\widetilde{W}_\infty(\lambda) = +\infty ,\ \widetilde{\mathrm{P}}_x^{\xi,\lambda}\textup{-a.s.} \ &\Longleftrightarrow \ W_\infty(\lambda) =0,\ \mathrm{P}_x^\xi\textup{-a.s.};\nonumber\\
		\widetilde{W}_\infty(\lambda) < +\infty ,\ \widetilde{\mathrm{P}}_x^{\xi,\lambda}\textup{-a.s.} \ &\Longleftrightarrow \ \mathrm{E}_x^\xi W_\infty(\lambda) =e^{-\lambda x}\phi(x,\lambda).
\end{align*}
Define
\begin{align*}
	\Omega_2:=
	\Omega_{\lambda^*}\cap \Omega_{-\lambda^*}
	& \bigcap\left\{ \limsup_{n\to\infty}\frac{\ln \psi(-n+x,
		\lambda^*)
	-\lambda^*\frac{\psi_{\lambda}(-n+x,\lambda^*)}{\psi(-n+x,\lambda^*)} }{\sqrt{n}} \geq 0\right\}\\ & \bigcap \left\{ \limsup_{n\to\infty}\frac{\ln \psi(n+x,-\lambda^*)
-\lambda^*\frac{\psi_{\lambda}(n+x,-\lambda^*)}{\psi(n+x,-\lambda^*)}}{\sqrt{n}} \geq 0\right\},
\end{align*}
here $\Omega_{\lambda^*}$  and $\Omega_{-\lambda^*}$  are defined in \eqref{Omega_Lambda}.
Then
$\mathbb{P}(\Omega_2)=1$ by {\bf(H1)}, \eqref{relation-lnpsi-lnphi}
 and Lemma \ref{Engenfunction-Property2}. We will fix $\omega\in \Omega_2$ in the rest of the proof.

(i)
We first consider the case $\lambda > \lambda^*$. In this case,
$$
\widetilde{W}_t(\lambda) \geq e^{-\gamma(\lambda)t}\psi(\Xi_t,\lambda) = \exp\left\{
\left(\frac{\ln \psi(\Xi_t,\lambda)}{\Xi_t}\times \frac{\Xi_t}{t} -\gamma(\lambda)\right)
t\right\}.
$$
By Lemma \ref{lemma 4.1} and the asymptotic behavior of $\ln \psi(\cdot, \lambda)$, we have
$$
\lim_{t\to\infty} \frac{\ln \psi(\Xi_t,\lambda)}{\Xi_t} \frac{\Xi_t}{t} = (-\lambda) \cdot \left(-\gamma'(\lambda)\right) = \lambda \gamma'(\lambda),\ \widetilde{\mathrm{P}}_x^{\xi,\lambda}\textup{-a.s.}
$$
By Lemma \ref{Engenvalue-Property}, when $|\lambda| > \lambda^*$, $\lambda \gamma'(\lambda) > \gamma(\lambda).$
Thus, $\widetilde{W}_\infty(\lambda) = +\infty,\  \widetilde{\mathrm{P}}_x^{\xi,\lambda}$-a.s.  or equivalently, $W_\infty(\lambda)=0,\  \mathrm{P}_x^\xi$-a.s.

Now we consider the case $\lambda =\lambda^*$. In this case,
we first prove that for any $M > 0$,
$$
\widetilde{\mathrm{P}}_x^{\xi,\lambda^*} \left(\ln \psi(\Xi_t,\lambda^*) - \gamma(\lambda^*)t > M\  \textup{i.o.}\right)=1.
$$
Note that by  \eqref{Change-of-measure},
$$
\gamma(\lambda^*)\widetilde{\mathrm{E}}_x^{\xi,\lambda} T_{-n} =  \gamma(\lambda^*)\Pi_x
\left(\exp\left\{\int_0^{H_{-n+x}}\left((m-1)\xi(B_s)-\gamma(\lambda^*)\right)\mathrm{d}s \right\} H_{-n+x}\right)
\frac{\psi(-n+x,\lambda^*)}{\psi(x,\lambda^*)}.
$$
By Lemma \ref{Engenvalue-Property},  $\gamma(\lambda^*)=\lambda^*\gamma'(\lambda^*)$.
Thus, using \eqref{Psi_2} and \eqref{Psi_3}, we have
\begin{align*}
	\gamma(\lambda^*)\widetilde{\mathrm{E}}_x^{\xi,\lambda} T_{-n} &  = \lambda^*\gamma'(\lambda^*)\frac{\Pi_x
\left(\exp\left\{\int_0^{H_{-n+x}}\left((m-1)\xi(B_s)-\gamma(\lambda^*)\right)\mathrm{d}s \right\} H_{-n+x}\right)
}{\Pi_x
\left(\exp\left\{\int_0^{H_{-n+x}}\left((m-1)\xi(B_s)-\gamma(\lambda^*)\right)\mathrm{d}s \right\} \right)}\\
& = \lambda^*
\left(\frac{\psi_{\lambda}(-n+x,\lambda^*)}{\psi(-n+x,\lambda^*)}-\frac{\psi_{\lambda}(x,\lambda^*)}{\psi(x,\lambda^*)}\right).
\end{align*}
Therefore,
\begin{align*}
& \widetilde{\mathrm{P}}_x^{\xi,\lambda^*}
\left(\ln \psi(\Xi_t,\lambda^*) - \gamma(\lambda^*)t > M\  \textup{i.o.}\right)  \geq \widetilde{\mathrm{P}}_x^{\xi,\lambda^*} \left(\ln \psi(\Xi_{T_{-n}},\lambda^*) - \gamma(\lambda^*)T_{-n} > M\  \textup{i.o.}\right) \\
&= \widetilde{\mathrm{P}}_x^{\xi,\lambda^*}
\left( \gamma(\lambda^*)\frac{ T_{-n} - \widetilde{\mathrm{E}}_x^{\xi,\lambda^*} T_{-n}}{\sqrt{n}}  < \frac{\ln\psi(-n+x,\lambda) - \gamma(\lambda^*)\widetilde{\mathrm{E}}_x^{\xi,\lambda^*} T_{-n}-M}{\sqrt{n}}\  \textup{i.o.}\right)\\
& = \widetilde{\mathrm{P}}_x^{\xi,\lambda^*} \left( \gamma(\lambda^*)\frac{ T_{-n} - \widetilde{\mathrm{E}}_x^{\xi,\lambda^*} T_{-n}}{\sqrt{n}}  < \frac{\ln\psi(-n+x,\lambda^*) -  \lambda^* \left(\frac{\psi_{\lambda^*}(-n+x,\lambda^*)}{\psi(-n+x,\lambda^*)}-\frac{\psi_{\lambda^*}(x,\lambda^*)}{\psi(x,\lambda^*)}\right)-M}{\sqrt{n}}\  \textup{i.o.}\right).
\end{align*}
Thanks to our choice of $\Omega_2$, there exists a subsequence
 $\{n_k=n_k(\omega),k\geq 1\}$
such that for all $\omega \in \Omega_2,$
$$
\lim_{k\to\infty} \frac{\ln\psi(-n_k+x,\lambda^*) -  \lambda^*
\left(\frac{\psi_{\lambda^*}(-n_k+x,\lambda^*)}{\psi(-n_k+x,\lambda^*)}-\frac{\psi_{\lambda^*}(x,\lambda^*)}{\psi(x,\lambda^*)}\right)-M}{\sqrt{n_k}} \geq 0,
$$
so
$$
\inf_{k\geq1} \frac{\ln\psi(-n_k+x,\lambda^*) -  \lambda^*
\left(\frac{\psi_{\lambda^*}(-n_k+x,\lambda^*)}{\psi(-n_k+x,\lambda^*)}-\frac{\psi_{\lambda^*}(x,\lambda^*)}{\psi(x,\lambda^*)}\right)-M}{\sqrt{n_k}} > -\infty.$$
Together with Remark \ref{remark}, we have
\begin{align*}
&\widetilde{\mathrm{P}}_x^{\xi,\lambda^*}
\left(\ln \psi(\Xi_t,\lambda^*) - \gamma(\lambda^*)t > M\  \textup{i.o.}\right)  \\
	&\geq \widetilde{\mathrm{P}}_x^{\xi,\lambda^*}
\left( \gamma(\lambda^*)\frac{ T_{-n_k} - \widetilde{\mathrm{E}}_x^{\xi,\lambda^*} T_{-n_k}}{\sqrt{n_k}}  < \frac{\ln\psi(-n_k+x,\lambda^*) -  \lambda^* \left(\frac{\psi_{\lambda^*}(-n_k+x,\lambda^*)}{\psi(-n_k+x,\lambda^*)}-\frac{\psi_{\lambda^*}(x,\lambda^*)}{\psi(x,\lambda^*)}\right)-M}{\sqrt{n_k}}\  \textup{i.o.}\right)\\
&=1.
\end{align*}
Letting $M \to \infty$,
$$
\widetilde{\mathrm{P}}_x^{\xi,\lambda^*}
\left(\limsup_{t\to\infty} \left(\ln \psi(\Xi_t,\lambda^*) - \gamma(\lambda^*)t\right) = +\infty\right)=1,
$$
this implies  $\widetilde{W}_\infty(\lambda^*) = +\infty,\  \widetilde{\mathrm{P}}_x^{\xi,\lambda^*}$-a.s. or equivalently, $W_\infty(\lambda^*)=0,\  \mathrm{P}_x^\xi$-a.s..

(ii) When $|\lambda| \in (\rho, \lambda^*)$, by Lemma \ref{Engenvalue-Property}, $\lambda \gamma'(\lambda) - \gamma(\lambda) < 0$.
Let $O$ be the law of the offspring of the spine particle,  that is, $\widetilde{\mathrm{P}}_x^{\xi,\lambda}(O=k) = kp_k/m$ for $k=1,2,...$.
If $\sum_{k=1}^\infty (k\ln k) p_k = +\infty$,
then
$$
\infty = \widetilde{\mathrm{E}}_x^{\xi,\lambda}
\left(\frac{\ln^+ O}{M}\right)
= \int_0^\infty \widetilde{\mathrm{P}}_x^{\xi,\lambda}(\ln^+ O > yM)\mathrm{d} y \leq \sum_{k=0}^\infty
\widetilde{\mathrm{P}}_x^{\xi,\lambda}(\ln^+ O > kM), \quad M>0.
$$
Let $\mu_k$ be the time of the $k$th fission  of the spine and $O_k$ be the number of offspring of the spine at time $\mu_k$. Then by the Borel-Cantelli lemma and the arbitrariness of $M$,  we have
$$
\limsup_{k\to\infty}\frac{\ln O_k}{k} = +\infty,\ \widetilde{\mathrm{P}}_x^{\xi,\lambda}\textup{-a.s.}
$$
Note that
$$
\widetilde{W}_{\mu_k} \geq (O_k -1)e^{-\gamma(\lambda)\mu_k}\psi(\Xi_{\mu_k},\lambda) = \exp\left\{
\left(\frac{\ln^+(O_k-1)}{\mu_k} +\frac{\ln \psi(\Xi_{\mu_k},\lambda)}{\mu_k} -\gamma(\lambda)\right)
\mu_k\right\}.
$$
Since $(m-1)\xi(x)\geq (m-1)\textup{ei}$, we have
$$\widetilde{\mathrm{P}}_x^{\xi,\lambda}(\mu_k \leq t) = \widetilde{\mathrm{P}}_x^{\xi,\lambda} (N_t \geq k) \geq  \widetilde{\mathrm{P}}_x^{\xi,\lambda} (\tilde{N}_t \geq k) = \widetilde{\mathrm{P}}_x^{\xi,\lambda}(\tilde{\mu}_k \leq t), $$
where, given the trajectory of $\Xi$, $N_t$ is a Poisson process with rate $(m-1)\xi(\Xi_t)$ and $\widetilde{N}_t$ is a Poisson process with rate $(m-1)\textup{ei}$,
$\tilde{\mu}_k:=\inf\{t: \widetilde{N}_t = k\}$.
We use a coupling of $(N_t, \widetilde{N}_t)$ such that $N_t -\tilde{N}_t$ is a Poisson process with rate $(m-1)\left( \xi(\Xi_t)-\textup{ei}\right)$ and is independent of $\tilde{N}_t$. So $\mu_k \leq \tilde{\mu}_k$.
Since $k/\tilde{\mu}_k  \to \left((m-1)\textup{ei}\right),\ \widetilde{\mathrm{P}}_x^{\xi,\lambda}$-a.s., we have $\liminf_{k\to\infty} k/\mu_k \geq \left((m-1)\textup{ei}\right),\ \widetilde{\mathrm{P}}_x^{\xi,\lambda}$-a.s.. Similarly we also have $\limsup_{k\to\infty} k/\mu_k \leq \left((m-1)\textup{es}\right),\ \widetilde{\mathrm{P}}_x^{\xi,\lambda}$-a.s.. Thus, $\mu_k \to \infty\ \widetilde{\mathrm{P}}_x^{\xi,\lambda}$-a.s. and
$$
\limsup_{k\to\infty}\widetilde{W}_{\mu_k} \geq  \limsup_{k\to\infty} \exp\left\{
\left(\frac{\ln^+(O_k-1)}{k} \frac{k}{\mu_k} +\frac{\ln \psi(\Xi_{\mu_k} ,\lambda)}{\mu_k} -\gamma(\lambda)\right)
\mu_k\right\} = +\infty.
$$

If $\sum_{k=1}^\infty (k\ln k) p_k < +\infty$, then for any $c>0$,
$$
\infty > \widetilde{\mathrm{E}}_x^{\xi,\lambda}
\left(\frac{\ln^+ O}{c}\right)
= \int_0^\infty \widetilde{\mathrm{P}}_x^{\xi,\lambda}(\ln^+ O > yc)\mathrm{d} y \geq \sum_{k=1}^\infty \widetilde{\mathrm{P}}_x^{\xi,\lambda}(\ln^+ O > kc) . $$
Since $c$ is arbitrary, we have
 $$\limsup_{k\to\infty}\frac{\ln O_k}{k} = 0,\quad \widetilde{\mathrm{P}}_x^{\xi,\lambda}\textup{-a.s.}$$
Define $\widetilde{\mathcal{G}}$ to be the $\sigma$-field of the genealogy along the spine,
then we have
$$
\widetilde{\mathrm{E}}_{x}^{\xi, \lambda}
\left(\widetilde{W}_t(\lambda) \Big|\widetilde{\mathcal{G}}  \right)
 = \psi (\Xi_t,\lambda) e^{-\gamma(\lambda)t} + \sum_{k=1}^{N_t-1} (O_k - 1)\psi(\Xi_{\mu_k},\lambda) e^{-\gamma(\lambda)\mu_k}.
$$
Here $O_k, \mu_k$ and $N_t$
have the same meaning as in  the case  $\sum_{k=1}^\infty( k\ln k) p_k = +\infty.$
Then, $\widetilde{\mathrm{P}}_x^{\xi,\lambda}$-a.s.,
$$\limsup_{k\to\infty} \frac{\ln(O_k -1)}{\mu_k} = \limsup_{k\to\infty} \frac{\ln(O_k -1)}{k}\frac{k}{\mu_k}=0.$$
Note that $\widetilde{\mathrm{P}}_x^{\xi,\lambda}$-a.s.,
$$\lim_{t\to\infty}\frac{\ln\psi (\Xi_t,\lambda)}{t}=\lim_{t\to\infty}\frac{\ln\psi (\Xi_t,\lambda)}{\Xi_t}\frac{\Xi_t}{t}=\lambda\gamma'(\lambda)<\gamma(\lambda),$$
and thus
$$\lim_{t\to\infty}\frac{\ln\psi (\Xi_t,\lambda)}{t}<\frac{\gamma(\lambda) + \lambda \gamma'(\lambda)}{2},$$
which implies
$$\limsup_{t\to\infty} \psi (\Xi_t,\lambda) \exp\left\{-\left(\gamma(\lambda) + \lambda \gamma'(\lambda)\right)t/2\right\}  = 0.$$
So there exists a $\widetilde{\mathrm{P}}_x^{\xi,\lambda}$-a.s. finite random variable
$\eta$
such that for all $k$ and $t$,
$$\ln(O_k -1)\leq \eta
+ \frac{\gamma(\lambda) - \lambda \gamma'(\lambda)}{4}\mu_k
,\ \textup{and}\  \psi (\Xi_t,\lambda) \exp\left\{-\left(\gamma(\lambda) + \lambda \gamma'(\lambda)\right)t/2\right\} \leq \eta.
$$
Therefore,
\begin{align*}
&	\widetilde{\mathrm{E}}_{x}^{\xi,\lambda}
\left(\widetilde{W}_t(\lambda) \Big|\widetilde{\mathcal{G}}  \right)
 = \psi (\Xi_t,\lambda) e^{-\gamma(\lambda)t} + \sum_{k=1}^{N_t-1} (O_k - 1)\psi(\Xi_{\mu_k},\lambda) e^{-\gamma(\lambda)\mu_k} \\  & \leq
\eta
\exp\left\{-\left(\gamma(\lambda) - \lambda \gamma'(\lambda)\right)t/2\right\} + \sum_{k=1}^\infty
\eta e^\eta
\exp\left\{ \exp\left\{-\left(\gamma(\lambda) - \lambda \gamma'(\lambda)\right)\mu_k/4\right\}\right\} < \infty,
\end{align*}
which implies $\limsup_{t \to \infty}\widetilde{\mathrm{E}}_{x}^{\xi,\lambda}
\left(\widetilde{W}_t(\lambda) \big|\widetilde{\mathcal{G}}  \right) < \infty, \ \widetilde{\mathrm{P}}_x^{\xi,\lambda}$-a.s.
Fatou's lemma for conditional probability implies that $\liminf_{t \to \infty}\widetilde{W}_t(\lambda)<\infty$, $\widetilde{\mathrm{P}}_{x}^{\xi,\lambda}$-a.s.
Since $\widetilde{W}_t(\lambda)^{-1}$ is a non-negative supermartingale under $\widetilde{\mathrm{P}}_x^{\xi,\lambda}$, we have $\lim_{t \to \infty}\widetilde{W}_t(\lambda)=\widetilde{W}_\infty(\lambda)<\infty$,
$\widetilde{\mathrm{P}}_x^{\xi,\lambda}$-a.s.
Thus we have shown that
when $\lambda \in (\rho, \lambda^*)$ and $\sum_{k=1}^\infty (k\ln k) p_k < \infty$, $W_t(\lambda)$ converges $\mathrm{P}_x^\xi$-a.s. and in $L^1$ to $W_\infty(\lambda)$.
\hfill$\Box$

\textbf{Proof of Corollary \ref{Cor3.2}:}
For any $\epsilon > 0$, there exists constant $c_\epsilon>0$
such that for all $x\in\mathbb{R}$,
$$
\psi(x, -\lambda^*) \geq c_\epsilon e^{\lambda^* x -\epsilon|x| }.
$$
Taking  $\lambda = -\lambda^*$ and using Theorem \ref{Additive martingale},
 we have, $\mathrm{P}_x^\xi$-a.s.,
$$
0 = W_\infty(-\lambda^*) \geq \limsup_{t\to\infty} \psi(M_t, -\lambda^*) e^{-\gamma(\lambda^*)t}\geq c_\epsilon
\limsup_{t \to \infty} \exp\left\{\lambda^* M_t - \epsilon|M_t| -\gamma(\lambda^*)t \right\}\geq 0.
$$
Thus, $\mathrm{P}_x^\xi$-a.s.,
\begin{equation*}
	\lim_{t\to\infty} \lambda^* M_t - \epsilon|M_t| -\gamma(\lambda^*)t = -\infty ,
\end{equation*}
 which implies that
 \begin{equation}\label{step_10}
 \limsup_{t \to \infty} \frac{\lambda^* M_t -\epsilon|M_t|}{t} \leq \gamma(\lambda^*),\quad \mathrm{P}_x^\xi\textup{-a.s.}
\end{equation}

When $\lambda \in (\rho, \lambda^*)$,  we will first show that for any $\omega\in \Omega_2$,
$\widetilde{\mathrm{P}}_x^{\xi,\lambda}$ and $\mathrm{P}_x^\xi$ are mutually absolutely continuous.
Define
$$q(x,\lambda,\omega) : = \mathrm{P}_x^\xi \left(W_\infty(\lambda)=0\right).$$
We only need to prove that
\begin{equation}\label{q=0} q(x,\lambda,\omega)=0,\quad \forall x\in\mathbb{R}.
\end{equation}
By the branching property, for $s <t$, we have
$$W_t(\lambda) = e^{\gamma(\lambda)s}\sum_{\nu\in N(s)} W_{t-s}^\nu(\lambda).$$
Here conditioned on
$\widetilde{\mathcal{F}}_s$,  $\{W_{t-s}^\nu(\lambda), \nu\in N(s)\}$
are independent. Letting $t\to\infty$ , we get
\begin{equation}\label{step_11}
	q(x,\lambda,\omega) = \mathrm{E}_x^\xi
	\left(\prod_{\nu\in N(s)} q(X_s^\nu, \lambda,\omega)\right).
\end{equation}
By the branching property and the Markov property, $\left\{\prod_{\nu\in N(s)} q(X_s^\nu, \lambda,\omega)\right\}_{s\geq 0}$ is a bounded martingale under $\mathrm{P}_x^\xi$. Since $W_\infty(\lambda)$ is the $L^1$ limit of $W_t(\lambda)$,  $q(x,\lambda,\omega) < 1$ for all $x\in \mathbb{R}.$ By \eqref{Prob-representation}, $q(\cdot,\lambda,\omega)$ satisfies the equation
\begin{align}\label{step_48}
	q_t(x,\lambda,\omega) = \frac{1}{2}q_{xx}(x,\lambda,\omega) + \xi(x)\left(\sum_{k=1}^\infty p_k q^k(x,\lambda,\omega) - q(x,\lambda,\omega) \right).
\end{align}
By the Feynman-Kac formula, we have
$$
q(x,\lambda,\omega) = \Pi_x
\left(q(B_t,\lambda,\omega)
\exp\left\{\int_0^t \xi(B_s)\left(\sum_{k=1}^\infty p_k q^{k-1}(B_s,\lambda,\omega) -1\right)\mathrm{d}s \right\} \right).
$$
Let $A(\varepsilon):= \{x : q(x,\lambda,\omega) \leq 1-\varepsilon \}$, then
\begin{align}\label{domi-q}
&q(x,\lambda,\omega)  \leq \Pi_x
\left(\exp\left\{\int_0^t \xi(B_s)\left(\sum_{k=1}^\infty p_k q^{k-1}(B_s,\lambda,\omega) -1\right)\mathrm{d}s \right\}\right)\\
& \leq \Pi_x
\left(\exp\left\{\int_0^t \textup{ei}(1-p_1)
\left( q(B_s,\lambda,\omega) -1\right)\mathrm{d}s \right\} \right)
\leq \Pi_x \left(exp\left\{-\varepsilon\textup{ei}(1-p_1)\int_0^t 1_{A(\varepsilon)}(B_s) \mathrm{d}s \right\} \right),
\end{align}
where in the second inequality, we used the domination $\left(\sum_{k=1}^\infty p_k q^{k-1}(B_s,\lambda,\omega) -1\right)\le p_1+(\sum^\infty_{k=2}p_i)q(B_s,\lambda,\omega)-1=(p_1-1)+(1-p_1)q(B_s,\lambda,\omega)$.
Let $m$ be the Lebesgue measure in $\mathbb{R}$. Note that if $m\left(A(\varepsilon)\right) > 0$, then as $t\to\infty$,
$$\int_0^t 1_{A(\varepsilon)}(B_s) \mathrm{d}s \to +\infty, \ \Pi_x\textup{-a.s.}$$
Since $q(x,\lambda,\omega)<1$, we can find
$\varepsilon=\varepsilon(\omega)>0$
such that $m\left(A(\varepsilon)\right) > 0.$
Letting $t\to\infty$ in \eqref{domi-q},
we obtain \eqref{q=0}, and
thus $\widetilde{\mathrm{P}}_x^{\xi,\lambda}$ and $\mathrm{P}_x^\xi$ are mutually absolute continuous.
Noticing that
$$\widetilde{\mathrm{P}}_x^{\xi,\lambda}\left(\liminf_{t\to\infty} \frac{\widetilde{M}_t}{t}\geq \gamma'(\lambda) \right) \geq \widetilde{\mathrm{P}}_x^{\xi,\lambda}\left(\liminf_{t\to\infty} \frac{\Xi_t}{t}\geq \gamma'(\lambda) \right)
=1,$$
we have$$ \mathrm{P}_x^\xi\left(\liminf_{t\to\infty} \frac{M_t}{t}\geq \gamma'(\lambda) \right) = 1. $$
Let $\lambda \uparrow \lambda^*$ we get that
\begin{equation}\label{step_12}
	\mathrm{P}_x^\xi\left(\liminf_{t\to\infty} \frac{M_t}{t}\geq v^* \right) =1.
\end{equation}
{\bf(H2)} implies that $\lambda^*>\rho$ and $\gamma(\lambda^*)>(m-1)
\textup{es}
.$ By Lemma \ref{Engenfunction-Property1} (4), $\gamma'(\lambda)\uparrow \gamma'(\lambda^*)=v^*$  as $\lambda\uparrow \lambda^*$.
In view of \eqref{step_12}, we know that $M_t \to +\infty,\ \mathrm{P}_x^\xi$-a.s. and thus \eqref{step_10} can be rewritten as
$$\limsup_{t \to \infty} \frac{ M_t}{t} \leq \frac{\gamma(\lambda^*)}{\lambda^* -\epsilon},\ \mathrm{P}_x^\xi\textup{-a.s.}$$
Letting $\epsilon \to 0+$ and using \eqref{step_12}, we get the conclusion of the corollary.
\hfill$\Box$

\subsection{Proof of Theorem \ref{Log-distance} }

In this subsection,
unless explicitly mentioned otherwise, we always assume
$\lambda = -\lambda^*$ and that {\bf(H1)}, {\bf(H2)} and {\bf(H3)} hold.
The proof of Theorem \ref{Log-distance} is inspired by
\cite[Section 6]{JD} and \cite[Section 4]{DS}.
Since the proof is long and complicated, we
divide it into several parts.
In this subsection,
$$K_1 := K_1(-\lambda^*),\quad K_2 := K_2(-\lambda^*)$$ are the two constants  in Lemma \ref{Engenfunction-Property1}(2).

For $s\in[0,t]$ and
$\Gamma_0>0$,  define
\begin{align*}
	&\varphi_t(s):=  \sup\left\{k\in \mathbb{Z}:\ s \geq  \frac{1}{\lambda^*v^*}\ln \psi(k, -\lambda^*)- \Gamma_0\ln t \right\}.
 \end{align*}
Later we will take $\Gamma_0>2C(\ell')+1$ with $C(\ell')$ being defined in Lemma \ref{lemma 4.8} (i) below.
We claim that for any $\Gamma_0>0$,
$\varphi_t(s)$ is c\`{a}dl\`{a}g and non-decreasing as a function of $s\in [0,t].$ In fact, if $\varphi_t(s) = k$, then by definition,
\begin{equation}\label{step_13}
	s\in \left[\frac{1}{\lambda^*v^*}\ln \psi(k, -\lambda^*)- \Gamma_0\ln t,\, \frac{1}{\lambda^*v^*}\ln \psi(k+1, -\lambda^*)- \Gamma_0\ln t \right),
\end{equation}
which implies that, for any
$r\in\left(s, \frac{1}{\lambda^*v^*}\ln \psi(k+1, -\lambda^*)- \Gamma_0\ln t\right)$,
 $\varphi_t(r) = k.$ On the other hand, suppose that $s_n\in [0,t]$,  $s_n \uparrow s$, and that $\varphi_t(s) = k$.
 If $s = \frac{1}{\lambda^*v^*}\ln \psi(k, -\lambda^*)-\Gamma_0 \ln t$
 , then $\varphi_t(s_n) \uparrow (k-1)$, if $s > \frac{1}{\lambda^*v^*}\ln \psi(k , -\lambda^*)- \Gamma_0\ln t$, then for $n$ large enough, we have $\frac{1}{\lambda^*v^*}\ln \psi(k, -\lambda^*)- \Gamma_0\ln t < s_n < s$, which implies that $\varphi_t(s_n) = k.$
Thus $\varphi_t(s)$ is a c\`{a}dl\`{a}g  function of $s\in [0,t]$.
The monotonicity follows easily from Remark \ref{remark3}.

Define
\begin{align}\label{def-L}
\mathcal{L}_t:= \# \Big\{ \nu \in N(t):& X_t^\nu \geq V_t +2, \min_{r\in[t-1,t]} X_r^\nu \geq V_t -3,\  X_s^\nu < \varphi_t(s),\ \forall s\in \left[H_{[V_t]}^{\nu}, t\right], \nonumber\\
	&  H_{k}^\nu  \geq \frac{1}{\lambda^*v^*}\ln \psi(k+1, -\lambda^*)- \Gamma_0\ln t, \ \forall k= 1,2,...,[V_t] \Big\},
\end{align}
here for any $y\in\R$ and $\nu\in N(t)$, $H^\nu_y$ is the first time that the particle $\nu$ hits $y$,
i.e.,
$
H^{\nu}_y=\inf\{t\geq0: X^\nu_t=y\}
$.
To prove Theorem \ref{Log-distance}, we will use the following inequality:
$$
\mathrm{P}_y^\xi \left(\mathcal{L}_{t} > 0\right)
\geq \frac {\left(\mathrm{E}_y^\xi\left(\mathcal{L}_{t}\right)\right)^2
}{\mathrm{E}_y^\xi\left((\mathcal{L}_{t})^2\right)},\quad t>0.
$$
We will estimate the first moment of $\mathcal{L}_t$
from below and the second moment of $\mathcal{L}_t$ from
above.
The following lemma gives the estimate of the second moment of $\mathcal{L}_t$.

\begin{lemma}\label{Second-moment}
There exists a non-random constant $\ell_2$ such that for
$t$ large enough,
$$\sup_{y\in[-1,1]} \mathrm{E}_y^\xi \left((\mathcal{L}_t)^2\right)\leq t^{\ell_2}.
$$
\end{lemma}
\textbf{Proof:}
Note when
$\nu \in N(t)$ makes a contribution to $\mathcal{L}_t$, we
have $X_s^{\nu} < \varphi_t(s)
$ for all $s\in[0,t]$. This is because that for
$s\in\left[H_1^\nu,  H_{[V_t]}^\nu\right)$,
there exists  $k\in\{2,...,[V_t]\}$ such that $s\in \left[ H_{k-1}^\nu , H_{k}^\nu \right)$. Therefore,
$$s\geq H_{k-1}^\nu \geq  \frac{1}{\lambda^*v^*}\ln \psi(k, -\lambda^*)- \Gamma_0\ln t,$$
which implies that $\varphi_t(s)\geq k$. But $s < H_{k}^\nu$ means $X_s^\nu < k$, so we get that $X_s^\nu < k \leq \varphi_t(s)$ for
$s\in\left[H_1^\nu, H_{[V_t]}\right).$
If $s\in [0, H_1^\nu)$, then when $t$ is large enough so that $0 \geq  \frac{1}{\lambda^*v^*}\ln \psi(1, -\lambda^*)- \Gamma_0\ln t$,
$$s \geq 0 \geq \frac{1}{\lambda^*v^*}\ln \psi(1, -\lambda^*)- \Gamma_0\ln t,$$
thus $\varphi_t(s)\geq k$ holds for $k=1$. Thus for $s< H_1^\nu$, we have $X_s^\nu < 1\leq \varphi_t(s)$.
In conclusion,
$$
\mathcal{L}_t \leq \# \left\{\nu \in N(t): X_t^\nu \geq V_t+2, X_s^\nu < \varphi_t(s),\ \forall s\in [0,t] \right\}.
$$
By Proposition \ref{Many-to-one},
\begin{align*}
\mathrm{E}_y^\xi \left((\mathcal{L}_t)^2\right)
& \leq \mathrm{E}_y^\xi \left(\mathcal{L}_t\right)
+ (m_2 - m)\int_0^t \Pi_y
\bigg( \exp\left\{\int_0^s (m-1)\xi(B_r)\mathrm{d}r \right\}\xi(B_s)1_{\{B_r < \varphi_t(r),\ \forall 0\leq r\leq s \} }\\
& \qquad \times\left(\Pi_z
 \left(\exp\left\{ \int_0^{t-s}(m-1)\xi(B_r)\mathrm{d} r\right\} 1_{\{B_r \leq \varphi_t(r+s),\ \forall 0\leq r\leq t-s, B_{t-s} \geq V_t+2 \}}\right)
\right)^2_{|z=B_s}
\bigg)\mathrm{d}s.
\end{align*}
Thus, for $z <\varphi_t(s)$, using the monotonicity of
$\psi(\cdot,-\lambda^*)$
(see Remark \ref{remark3}) and  \eqref{F-K-u}, we have
\begin{align}\label{step_2}
&\Pi_z
\left(\exp\left\{ \int_0^{t-s}(m-1)\xi(B_r)\mathrm{d} r\right\} 1_{\{B_r \leq \varphi_t(r+s),\ \forall 0\leq r\leq t-s,B_{t-s} \geq V_t+2 \}}\right)
\nonumber\\
&\leq \Pi_z
\left(\exp\left\{ \int_0^{t-s}(m-1)\xi(B_r)\mathrm{d} r\right\} 1_{\{B_{t-s} \geq V_t +2\}}\right)
\nonumber \\& \leq  \Pi_z
\left(\exp\left\{ \int_0^{t-s}(m-1)\xi(B_r)\mathrm{d} r\right\} \frac{\psi(B_{t-s},-\lambda^*)}{\psi(V_t+2,-\lambda^*)}\right)
 =\frac{\psi(z,-\lambda^*)e^{\gamma(\lambda^*)(t-s)}}{\psi(V_t+2,-\lambda^*)}\nonumber
	\\ &\leq \frac{\psi\left(\varphi_t(s),-\lambda^*\right) e^{\gamma(\lambda^*)(t-s)}}{\psi(V_t,-\lambda^*)} =\psi \left(\varphi_t(s),-\lambda^*\right) e^{-\lambda^* v^* s},
\end{align}
where in the last equality we used the fact that $e^{-\gamma(\lambda^*)t}\psi(V_t,-\lambda^*)=1$ given by \eqref{equality-V_t}.
Since \eqref{step_13} implies that for $s\in[0,t]$,
$$
s\geq  \frac{1}{\lambda^*v^*}\ln \psi(\varphi_t(s), -\lambda^*)- \Gamma_0\ln t,
$$
the last term of \eqref{step_2} is bounded from above  by $t^{\Gamma_0\lambda^* v^*}$ for all $z < \varphi_t(s), s \in [0,t]$.
Similarly as in \eqref{step_2}, for $s=0$ we have
\begin{align*}
&\Pi_y \left(\exp\left\{ \int_0^{t}(m-1)\xi(B_r)\mathrm{d} r\right\} 1_{\{B_{t} \geq V_t+2 \}}\right)\\
& \leq \Pi_y
 \left(\exp\left\{ \int_0^{t}(m-1)\xi(B_r)\mathrm{d} r\right\}\frac{\psi(B_{t}, -\lambda^*)}{\psi(V_t +2, -\lambda^*)}\right)
= \frac{\psi(y,-\lambda^*)e^{\gamma(\lambda^*)t}} {\psi(V_t+2,-\lambda^*)}\\	& \leq \frac{\psi(1,-\lambda^*)e^{\gamma(\lambda^*)t}}{\psi(V_t,-\lambda^*)}= \psi(1,-\lambda^*)\leq e^{K_1},
\end{align*}
where in the last inequality we used \eqref{difference-psi'}.
Therefore  we conclude that
\begin{align*}
\mathrm{E}_y^\xi \left((\mathcal{L}_t)^2\right)
& \leq \mathrm{E}_y^\xi \left(\mathcal{L}_t\right)
+ (m_2 - m)\int_0^t \Pi_y
\bigg( \exp\left\{\int_0^s (m-1)\xi(B_r)\mathrm{d}r \right\}\xi(B_s)1_{\{B_r < \varphi_t(r),\ \forall 0\leq r\leq s \} }\\
& \qquad \times\left(\Pi_z
\left(\exp\left\{ \int_0^{t-s}(m-1)\xi(B_r)\mathrm{d} r\right\} 1_{\{B_r \leq \varphi_t(r+s),\ \forall 0\leq r\leq t-s, B_{t-s} \geq V_t+2 \}}\right)
\right)^2_{|z=B_s}\bigg)\mathrm{d}s\\
& \leq \mathrm{E}_y^\xi
\left(\#\left\{\nu \in N(t): X_t^\nu \geq V_t+2 \right\}\right)
\\ & \qquad+ t^{\Gamma_0\lambda^* v^*}(m_2 - m)\textup{es}\int_0^t \Pi_y
\bigg(
\exp\left\{\int_0^s (m-1)\xi(B_r)\mathrm{d}r \right\}1_{\{B_r < \varphi_t(r),\ \forall 0\leq r\leq s \} }\\
& \qquad \times \Pi_z
\left(\exp\left\{ \int_0^{t-s}(m-1)\xi(B_r)\mathrm{d} r\right\} 1_{\{B_r \leq \varphi_t(r+s),\ \forall 0\leq r\leq t-s, B_{t-s} \geq V_t+2 \}}\right)_{|z=B_s}\bigg)\mathrm{d}s \\
& \leq e^{K_1} + t^{\Gamma_0\lambda^* v^*}(m_2-m)\textup{es} \int_0^t\Pi_y
\left(\exp\left\{ \int_0^{t}(m-1)\xi(B_r)\mathrm{d} r\right\} 1_{\{B_{t} \geq V_t+2 \}}\right)
\mathrm{d}s\\
	& \leq e^{K_1} + t^{\Gamma_0\lambda^* v^*+1}e^{K_1}(m_2-m)\textup{es} .
\end{align*}
Taking $\ell_2 > \Gamma_0\lambda^* v^* +1$,  we arrive at the desired conclusion.
\hfill$\Box$

\bigskip
Now we estimate
$\mathrm{E}_y^{\xi}(\mathcal{L}_t)$ from below.
Define
\begin{align}\label{step_16}
G_t &:= \left\{H_{k}  \geq \frac{1}{\lambda^*v^*}\ln \psi(k
+1, -\lambda^*)- \Gamma_0\ln t, \ \forall k= 1,2,...,[V_t]\right\};\nonumber\\
A(r,t)& := \left\{ B_t \geq V_t +2, \min_{s\in[t-1,t]} B_s \geq V_t -3,\  B _s < \varphi_t(s),\ \forall s\in \left[r, t\right] \right\};\\
\Gamma_1
&:= K_1/(\lambda^* v^*)+2.
\end{align}
For $y\in[-1,1],$ by Proposition \ref{Many-to-one},
 \begin{align}\label{step_15}
&	\mathrm{E}_y^{\xi}(\mathcal{L}_t)
= \Pi_y
\left(\exp\left\{\int_0^t(m-1)\xi(B_s)\mathrm{d}s \right\}; A\left(H_{[V_t]},t \right), G_t\right)
\nonumber \\& \geq \Pi_y
\left(\exp\left\{\int_0^t(m-1)\xi(B_s)\mathrm{d}s \right\}; A\left(H_{[V_t]},t \right), H_{[V_t]} \in [t-\Gamma_1, t-1], G_t\right)\nonumber\\& \geq \Pi_y
 \left(\exp\left\{\int_0^{H_{[V_t]}}(m-1)\xi(B_s)\mathrm{d}s \right\}; A\left(H_{[V_t]},t \right), H_{[V_t]} \in [t-\Gamma_1, t-1], G_t\right).
\end{align}
By the definition of $\varphi_t(s)$, we know that
$$
\frac{1}{\lambda^*v^*}\ln \psi(\varphi_t(s), -\lambda^*)- \Gamma_0\ln t \leq s = \frac{1}{\lambda^*v^*}\ln \psi(V_s, -\lambda^*) < \frac{1}{\lambda^*v^*}\ln \psi(\varphi_t(s)
+1, -\lambda^*)- \Gamma_0\ln t.
$$
It follows from Lemma \ref{Engenfunction-Property1}(2) that  for $s\in[0,t]$,
\begin{align}\label{step_18}
\frac{\Gamma_0 \lambda^* v^* }{K_1}\ln t -1 \leq 	\varphi_t(s)
 - V_s \leq \frac{\Gamma_0 \lambda^* v^* }{K_2}\ln t.
\end{align}
By \eqref{difference-V},  when $r\in [t-\Gamma_1, t-1]$ and $s\in[r,t],$ $0\leq V_t -V_s \leq \frac{\lambda^* v^*}{K_2}(t-s)\leq \frac{\lambda^* v^*}{K_2}\Gamma_1.$
Therefore, applying $\eqref{step_18}$ for $s\in[r,t]$, we have
\begin{align*}
&A(r,t)  \supset \left\{ B_t \geq V_t +2, \min_{s\in[t-1,t]} B_s \geq V_t -3,\  B _s < V_s +\frac{\Gamma_0 \lambda^* v^* }{K_1}\ln t -1,\ \forall s\in \left[r, t\right] \right\}\\
& \supset \left\{ B_t \geq V_t +2, \min_{s\in[t-1,t]} B_s \geq V_t -3,\  B _s < V_t +\frac{\Gamma_0 \lambda^* v^* }{K_1}\ln t -1-\frac{\lambda^*v^*}{K_2}\Gamma_1,\ \forall s\in \left[r, t\right] \right\}
\\& \supset \left\{ B_t \geq [V_t]
+3, \min_{s\in[t-1,t]} B_s \geq [V_t]
-2,\  B _s < [V_t]
+\frac{\Gamma_0 \lambda^* v^* }{K_1}\ln t
-1 -\frac{\lambda^*v^*}{K_2}\Gamma_1,\ \forall s\in \left[r, t\right] \right\}.
\end{align*}
Hence, when $r =H_{[V_t]
} \in [t-\Gamma_1,t-1]$ and $t$ is large enough so that $\frac{\Gamma_0 \lambda^* v^* }{K_1}\ln t
-1 -\frac{\lambda^*v^*}{K_2}\Gamma_1 > 10$,
\begin{align*}
&A(H_{[V_t]},t)  \supset \left\{ B_t \geq B_r +
3, \min_{s\in[t-1,t]} B_s \geq B_r -
2,\  B _s < B_r +10,\ \forall s\in \left[r, t\right] \right\}\Big|_{r=H_{[V_t]} }.
\end{align*}
Using \eqref{step_15} and
the strong Markov property of Brownian motion $B$,
we have for $t$ large,
\begin{align}\label{step_15'}	
\mathrm{E}_y^{\xi}[\mathcal{L}_t] \geq &\Pi_y
\left(\exp\left\{\int_0^{H_{[V_t]}}(m-1)\xi(B_s)\mathrm{d}s \right\}
F_t;
H_{[V_t]} \in [t-\Gamma_1, t-1], G_t\right),
\end{align}
where $F_t$ is defined by
\begin{align*}
F_t:=&\Pi_{0}
\left(B_{t-r} \geq
3, \min_{s\in [t-r-1, t-r]} B_s \geq -
2,\  B_s <10,\ \forall s\in \left[0, t-r\right] \right)_{|r = H_{[V_t]}}.
\end{align*}
For $r =H_{[V_t]} \in [t-\Gamma_1,t-1]$,
\begin{align*}
F_t\geq \inf_{r\in [1,\Gamma_1]} \Pi_{0}
\left(B_{r} \geq
3, \min_{s\in [r-1,r]} B_s \geq -
2, \sup_{s\in [0,r] } B_s \leq 10  \right)
=:	\Gamma_2>0.
\end{align*}
Plugging the lower bound of $F_t$ into \eqref{step_15'},
using \eqref{Change-of-measure} with $\lambda = -\lambda^*$,  for large $t$, we have that
\begin{align}\label{step_17}
& \mathrm{E}_y^{\xi}(\mathcal{L}_t)\geq \Gamma_2 \Pi_y
\left(\exp\left\{\int_0^{H_{[V_t]}}(m-1)\xi(B_s)\mathrm{d}s \right\} ; H_{[V_t]} \in [t-\Gamma_1, t-1], G_t\right)\nonumber \\&= \Gamma_2\widetilde{\Pi}_y^{\xi, -\lambda^*}
\left(e^{\lambda^* v^* H_{[V_t]}} \frac{\psi(y,-\lambda^*)}{\psi([V_t],-\lambda^*)}; H_{[V_t]} \in [t-\Gamma_1, t-1], G_t\right)\nonumber\\
&\geq \Gamma_2\widetilde{\Pi}_y^{\xi, -\lambda^*}
\left(e^{\lambda^* v^* (t-\Gamma_1)} \frac{\psi(-1,-\lambda^*)}{\psi(V_t
	,-\lambda^*)}; H_{[V_t]} \in [t-\Gamma_1, t-1], G_t\right)\nonumber\\
&\geq  \Gamma_2 e^{-\lambda^* v^* \Gamma_1} e^{-K_1}
\widetilde{\Pi}_y^{\xi, -\lambda^*}
\left( H_{[V_t]} \in [t-\Gamma_1, t-1], G_t\right)\nonumber\\
&\geq  \Gamma_3 \widetilde{\Pi}_y^{\xi, -\lambda^*}
\left( H_{[V_t]} \in [t-\Gamma_1, t-1], H_{k}  \geq \frac{1}{\lambda^*v^*}\ln \psi(k, -\lambda^*)- (\Gamma_0-1)\ln t, \ 1\leq k \leq [V_t]\right).
\end{align}
Here in the penultimate inequality,
we used $\ln \psi(V_t, -\lambda^*)=\gamma(\lambda^*)t$ and \eqref{difference-psi'} with $z=0, y=-1$, and
in the last, we used
$\Gamma_3 := \Gamma_2 e^{-\lambda^* v^* \Gamma_1} e^{-K_1}$
and $t$ is large enough so that $\lambda^* v^* \ln t > K_1.$

\bigskip
To continue the estimate of $\mathrm{E}_y^{\xi}(\mathcal{L}_t)
$, we will consider two independent copies
of $(\Xi_t, \widetilde{\Pi}_y^{\xi, -\lambda^*})$.
For $y\in [-1, 1]$, we enlarge the probability space and the corresponding measure $\widetilde{\Pi}_y^{\xi, -\lambda^*}$: let $\Xi_t^j, j=1,2$ be two independent copies of $(\Xi_t, \widetilde{\Pi}_y^{\xi, -\lambda^*})$ and
and
 \begin{equation}\label{define-H^j}
H^j_x:= \inf\left\{t\geq 0:\ \Xi_t^j=x \right\},
\quad j=1,2, x\in \R.
  \end{equation}
For $j=1,2$, define
\begin{equation}\label{def-beta^j}
\beta_0^j:=0,
\quad \beta_k^j := H_{k}^j -\frac{1}{\lambda^*v^*}\ln \psi(k, -\lambda^*),\quad k\geq 1.
\end{equation}
Note that $\lambda^* v^* t - K_1 \leq \ln \psi([V_t], -\lambda^*) \leq \lambda^*v^*t,$
we have
$$
\left[\frac{1}{\lambda^*v^*}\ln \psi([V_t],
-\lambda^*)+\frac{
K_1 }{\lambda^* v^*}-\Gamma_1, \frac{1}{\lambda^*v^*}\ln \psi([V_t], -\lambda^*)
-1\right]\subset  [t-\Gamma_1, t-1].
$$
Together with the definition of $\Gamma_1$ in \eqref{step_16},
we continue the estimate of \eqref{step_17} and get
\begin{align}\label{step_23}
\mathrm{E}_y^{\xi}(\mathcal{L}_t)\geq  \Gamma_3
\widetilde{\Pi}_y^{\xi, -\lambda^*}
\left( \beta_{[V_t]}^j \in
[-2,-1]
, \beta_k^j \geq- (\Gamma_0-1)\ln t, \ 1\leq k \leq [V_t]\right),
\, j=1,2.
\end{align}

We first claim that,  under $\widetilde{\Pi}_y^{\xi,-\lambda^*}$, $\{\beta_k^1, k\geq 0\}$ has independent increments.
The proof of this claim is postponed to the Appendix.
Let $\pi:= \pi_t$ be a uniform random variable on $\{
2,...,[V_t]-1\}$ which is independent of $\Xi_t,\Xi_t^1$ and $\Xi_t^2$.
Define
$$
\beta_{k}^{(t)}:=\left\{\begin{array}{ll}
	\beta_{k}^{1}, & 1 \leq k \leq \pi, \\
	\beta_{\pi}^{1}+\left(\beta_{k}^{2}-\beta_{\pi}^{2}\right), & \pi < k \leq [V_t].
\end{array}\right.
$$
We claim that
\begin{equation}\label{claim-eq-d}
	\left(\beta_{k}^{(t)}, k = 1,...,[V_t]; \, \widetilde{\Pi}_y^{\xi, -\lambda^*}\right)\stackrel{\mathrm{d}}{=}\left(\beta_k^1, k=1,..., [V_t]; \, \widetilde{\Pi}_y^{\xi, -\lambda^*}\right).
\end{equation}
In fact, for real numbers $\alpha_k,k=1,...,[V_t],$ we have
\begin{align*}
&\widetilde{\Pi}_y^{\xi, -\lambda^*} \left(\exp\left\{\mathrm{i} \sum_{k=1}^{[V_t]}\alpha_k\beta_{k}^{(t)} \right\} \Bigg| \beta_k^1, \beta_k^2\right)\\
& = \frac{1}{[V_t]-2}\sum_{m=
2}^{[V_t]-1}\exp\left\{\mathrm{i}\sum_{k=1}^m\alpha_k \beta_k^1 + \mathrm{i}\sum_{k=m+1}^{[V_t]} \alpha_k\left(\beta_m^1 + \beta_k^2-\beta_m^2\right)\right\}.
\end{align*}
Thus, it suffices to show for any $2 \leq m \leq [V_t]-1$,
\begin{equation}\label{step_19}
\widetilde{\Pi}_y^{\xi, -\lambda^*}
\left(\exp\left\{\mathrm{i}\sum_{k=1}^m\alpha_k \beta_k^1 + \mathrm{i}\sum_{k=m+1}^{[V_t]} \alpha_k\left(\beta_m^1 + \beta_k^2-\beta_m^2\right)\right\}\right) 	
= \widetilde{\Pi}_y^{\xi, -\lambda^*}
\left(\exp\left\{\mathrm{i} \sum_{k=1}^{[V_t]}\alpha_k\beta_{k}^1 \right\} \right).
\end{equation}
Since under $\widetilde{\Pi}_y^{\xi, -\lambda^*}$,
$\beta_k^j, j=1,2$, are sums of independent random variables, we know that, for $k>m, \beta_k^1-\beta_m^1$ is independent of $\beta_r^1, r \leq m$. Use this observation, the left-hand side of \eqref{step_19} is
\begin{align*}
&\widetilde{\Pi}_y^{\xi, -\lambda^*}
\left(\exp\left\{\mathrm{i}\sum_{k=1}^m\alpha_k \beta_k^1 + \mathrm{i}\sum_{k=m+1}^{[V_t]} \alpha_k\left(\beta_m^1 + \beta_k^2-\beta_m^2\right)\right\}\right) \\
&= \widetilde{\Pi}_y^{\xi, -\lambda^*}
\left(\exp\left\{\mathrm{i}\sum_{k=1}^m\alpha_k \beta_k^1 + \mathrm{i}\sum_{k=m+1}^{[V_t]} \alpha_k \beta_m^1\right\}\right)\cdot
\widetilde{\Pi}_y^{\xi, -\lambda^*}
\left(\exp\left\{\mathrm{i}\sum_{k=m+1}^{[V_t]} \alpha_k\left( \beta_k^2-\beta_m^2\right)\right\} \right)\\
&= \widetilde{\Pi}_y^{\xi, -\lambda^*}
\left(\exp\left\{\mathrm{i}\sum_{k=1}^m\alpha_k \beta_k^1 + \mathrm{i}\sum_{k=m+1}^{[V_t]} \alpha_k \beta_m^1\right\}\right)\cdot
\widetilde{\Pi}_y^{\xi, -\lambda^*}
\left(\exp\left\{\mathrm{i}\sum_{k=m+1}^{[V_t]} \alpha_k\left( \beta_k^1-\beta_m^1\right)\right\} \right)\\
& =  \widetilde{\Pi}_y^{\xi, -\lambda^*}
\left(\exp\left\{\mathrm{i}\sum_{k=1}^m\alpha_k \beta_k^1 + \mathrm{i}\sum_{k=m+1}^{[V_t]} \alpha_k\left(\beta_m^1 + \beta_k^1-\beta_m^1\right)\right\} \right).
\end{align*}
Hence we obtain \eqref{step_19}, and thus the claim \eqref{claim-eq-d} holds.

Therefore, we have
\begin{equation}\label{step_22}
\mathrm{E}_y^{\xi}(\mathcal{L}_t)
\geq \Gamma_3	\widetilde{\Pi}_y^{\xi, -\lambda^*}
\left( \beta_{[V_t]}^{(t)} \in
[-2,-1], \beta_k^{(t)} \geq- (\Gamma_0-1)\ln t, \ 1\leq k \leq [V_t]\right).
\end{equation}
Now we are going to  give a lower bound of  the right hand side of \eqref{step_22}.
For a constant $C(\ell')$ to be specified in Lemma  \ref{lemma 4.8}(i) below, define
\begin{align}\label{def-I}
I_1 & := \left\{\beta_{k}^1 - \beta_2^1 \geq 0,\ \forall
2\leq k \leq [V_t],\, \beta_{[V_t]}^1 - \beta_2^1 \geq [V_t]^{1/4}\right\}; \nonumber\\
I_2 & := \left\{\beta_{k}^2 - \beta_{[V_t]}^2 \geq 0,\ \forall 2
\leq k \leq [V_t], \,\beta_2^2 - \beta_{[V_t]}^2 \geq [V_t]^{1/4}\right\};\nonumber\\
I_3 &:= \left\{\max_{3\leq k\leq [V_t]} \left|\beta_k^1 - \beta_{k-1}^1\right|\leq C(\ell')\ln [V_t] \right\}; \\
I_4 &:= \left\{\max_{3\leq k\leq [V_t]} \left|\beta_k^2 - \beta_{k-1}^2\right|\leq C(\ell')\ln [V_t] \right\};\\
I_5 &:= \left\{  \beta_2^1 \in [
-2+\beta_J^2 -\beta_{[V_t]}^2 -\beta_J^1 + \beta_2^1 ,
-1 +\beta_J^2 -\beta_{[V_t]}^2 -\beta_J^1 + \beta_2^1 ]\right\}\bigcap \left\{\pi = J \right\}.
\end{align}
Here $J:= \sup\{j\leq [V_t]:
-2 +\beta_j^2 -\beta_{[V_t]}^2 -\beta_j^1 + \beta_2^1 \geq 0\}$.
To give a lower bound of \eqref{step_22}, we will prove that for sufficiently large $t$,
\begin{equation}\label{Claim1}
\bigcap_{n=1}^5 I_n \subset \left\{\beta_{[V_t]}^{(t)} \in [-2,-1], \beta_k^{(t)} \geq-\left(\Gamma_0-1\right) \ln t, \ 1\leq k \leq [V_t] \right\}.
\end{equation}
By the definition of $\beta_k^{(t)} $, we know that, when $\bigcap_{n=1}^5 I_n $ occurs, we have $\beta_k^{(t)} = \beta_k^1$ if $1\leq k\leq \pi=J$ and
$\beta_k^{(t)} = \beta_J^1 + \beta_k^2 -\beta_J^2$ if $J < k\leq [V_t]$. Note that on $I_1\cap I_2$ we have
$$
\beta_{[V_t]}^1 - \beta_2^1 \geq [V_t]^{1/4}, \qquad
\beta_2^2 -\beta_{[V_t]}^2 \geq  [V_t]^{1/4}.
$$
This implies that for $t$ large enough so that $[V_t]^{1/4}> 2$,
$$
-2+\beta_2^2 -\beta_{[V_t]}^2 -\beta_2^1 + \beta_2^1  > 0 \ \textup{and}\
-2 +\beta_{[V_t]}^2 -\beta_{[V_t]}^2 -\beta_{[V_t]}^1 + \beta_2^1<0,
$$
which implies that
$2\leq J <[V_t]$.
When $k\leq J$, then on $I_1\cap I_3$ we have
for $k\geq 2$,
$$
\beta_k^{(t)} = \beta_k^1\geq \beta_2^1
= H_2^1 - \frac{1}{\lambda^* v^*} \ln \psi(2,-\lambda^*) \geq -\frac{1}{\lambda^* v^*} \ln \psi(2,-\lambda^*)\geq   -\frac{2K_1}{\lambda^* v^*}
$$
 and for $k=1$, $\beta_k^{(t)} = \beta_1^1
 =H_1^1 -(\lambda^* v^*)^{-1} \ln \psi(1,-\lambda^*)\geq -K_1/(\lambda^* v^*).
 $
 Choose $\Gamma_0 >2 C(\ell')+1$. For large $t$, we have
$
-2K_1/(\lambda^* v^*)
 > -(\Gamma_0-1)\ln t$ and hence $\beta_k^{(t)} > -(\Gamma_0-1) \ln t.$
Similarly, when $J< k \leq [V_t]$, since $
-2 +\beta_{J+1}^2 -\beta_{[V_t]}^2 -\beta_{J+1}^1 + \beta_2^1 < 0$
and $J+1\geq 3$,
, on $I_2\cap I_3\cap I_4$,
\begin{align*}
\beta_k^{(t)} & = \beta_J^1 + \beta_k^2 -\beta_J^2 \geq \beta_J^1 + \beta_{[V_t]}^2 -\beta_J^2 \geq -2 C(\ell')\ln [V_t] + \beta_{J+1}^1 +\beta_{[V_t]}^2 -\beta_{J+1}^2 \\
& > -2 C(\ell')\ln [V_t] +\beta_2^1
-2\geq -2 C(\ell')\ln [V_t]- \frac{2 K_1}{\lambda^* v^*}-2.
\end{align*}
Since $\Gamma_0 > 2C(\ell')+1$, we get that for large $t$, $\beta_k^{(t)}  > -(\Gamma_0-1) \ln t.$ Note that
$\beta_{[V_t]}^{(t)}  = \beta_J^1 + \beta_{[V_t]}^2 -\beta_J^2$ and by the definition of $I_5$, $\beta_{[V_t]}^{(t)} \in
[-2,-1]$ and hence we get \eqref{Claim1}.

By \eqref{step_22} and \eqref{Claim1}, we have
\begin{equation}\label{step_28}
\mathrm{E}_y^{\xi}(\mathcal{L}_t)
\geq \Gamma_3 \widetilde{\Pi}_y^{\xi, -\lambda^*} \left(\cap^5_{j=1}I_j\right)
= \Gamma_3 \widetilde{\Pi}_y^{\xi, -\lambda^*}
\left(\cap^4_{j=1}I_j\right)
\cdot \widetilde{\Pi}_y^{\xi, -\lambda^*}
\left(I_5 \big| \cap^4_{j=1}I_j\right).
\end{equation}
 Since $I_1$ and $I_2$ are independent, we have	
\begin{align}\label{I1234}
\widetilde{\Pi}_y^{\xi, -\lambda^*} \left(\cap^4_{j=1}I_j\right)
&\geq \widetilde{\Pi}_y^{\xi, -\lambda^*}
\left(I_1\cap I_2\right)+
\widetilde{\Pi}_y^{\xi, -\lambda^*}
\left(I_3\cap I_4\right)-1\\
&\geq\widetilde{\Pi}_y^{\xi, -\lambda^*} \left(I_1\right)
\cdot \widetilde{\Pi}_y^{\xi, -\lambda^*}
\left(I_2\right)+ \widetilde{\Pi}_y^{\xi, -\lambda^*}
\left(I_3\cap I_4\right)-1.
 \end{align}
By our choice of $J$, conditioned on $\cap^4_{j=1}I_j$,
we know that
\begin{equation}\label{def-x1}
x_1:=
-2 +\beta_J^2 -\beta_{[V_t]}^2 -\beta_J^1 + \beta_2^1 \in \left[0, 2C(\ell')\ln [V_t]\right],
\end{equation}
and $x_1$  is independent of $\beta_ 2^1$.
Also note that $\beta_2^1$ is independent of $\cap_{j=1}^4 I_j$ by Lemma \ref{indep-incre}.
 Thus,
\begin{align}\label{step_29}
\widetilde{\Pi}_y^{\xi, -\lambda^*}
\left(I_5 \big| \cap^4_{j=1}I_j\right)
& \geq \frac{1}{[V_t]-2}\inf_{x_1 \in [0, 2C(\ell ') \ln [V_t]]} \widetilde{\Pi}_y^{\xi, -\lambda^*}
\left(\beta_2^1 \in [x_1,x_1+1]\right).
\end{align}
Thus, to continue the estimate in \eqref{step_28},
we need to estimate $\widetilde{\Pi}_y^{\xi, -\lambda^*}
\left(I_3\cap I_4\right)
$, $\widetilde{\Pi}_y^{\xi, -\lambda^*}
\left(\beta_2^1 \in [x, x+1]\right)
$ and $\widetilde{\Pi}_y^{\xi, -\lambda^*}
\left(I_j\right)$ for $j=1, 2$.
Now we first estimate $\widetilde{\Pi}_y^{\xi, -\lambda^*}
\left(I_3\cap I_4\right)
$ and $\widetilde{\Pi}_y^{\xi, -\lambda^*}
\left(\beta_2^1 \in [x, x+1]\right)$.

\begin{lemma}\label{lemma 4.8}
(i) For any $\ell'>0$,
there exists a non-random constant $C(\ell')>0$
so that for large $t$,
\begin{equation}\label{step_20}
\inf_{y\in[-1,1]} \widetilde{\Pi}_y^{\xi, -\lambda^*} \left(I_3\cap I_4\right)\geq 1-[V_t]^{-\ell'}.
 \end{equation}

(ii) There exists a constant
$C_3>0$
such that for all $x >0$,
\begin{equation}
\inf_{y\in[-1,1]}\inf_{x_1 \in [0, x]}\widetilde{\Pi}_y^{\xi, -\lambda^*}\left(\beta_2^1 \in [x_1, x_1+1]\right)
\geq C_3e^{-x/C_3}.
\end{equation}
\end{lemma}
\textbf{Proof:}
(i) By the definition of $\beta_k^j$ and Lemma \ref{Engenfunction-Property1}(2),  for
$k\geq 3$,
\begin{align*}
\left|\beta_k^j - \beta_{k-1}^j\right| \leq  H_{k}^j - H_{k-1}^j  + \frac{K_1}{\lambda^* v^*}.
\end{align*}
For any $0 \leq \eta
< \gamma(-\lambda^*)- (m-1)\textup{es}
=\gamma(\lambda^*)- (m-1)\textup{es}
$, by the strong Markov property of $\Xi^j$,
\begin{align}\label{step_33}
& \widetilde{\Pi}_y^{\xi, -\lambda^*}
\left(\exp\left\{\eta (H_{k}^j - H_{k-1}^j )\right\} \right)
= \widetilde{\Pi}_{k-1}^{\xi, -\lambda^*}
\left(\exp\left\{\eta H_{k}^j \right\} \right) \nonumber\\
&= \Pi_{k-1}
\left(\exp\left\{\int_0^{H_{k
}}\left((m-1)\xi(B_s)-\gamma(-\lambda^*) \right)\mathrm{d}s + \eta H_{k}\right\}\frac{\psi(k,-\lambda^*)}{\psi(k-1,-\lambda^*)}\right)\nonumber\\
& \leq  e^{K_1} \Pi_{k-1}
\left(\exp\left\{\left((m-1)\textup{es}-\gamma(-\lambda^*)+\eta \right)H_{k}\right\}\right)
\nonumber\\
& =  \exp\left\{-\sqrt{2\left(\gamma(-\lambda^*) - (m-1)\textup{es} -\eta\right)} + K_1\right\}
=: c_1
< +\infty.
\end{align}
Hence,
for fixed $\eta\in (0,  \gamma(\lambda^*)- (m-1)\textup{es})$,
we have for all $y \in [-1, 1]$,
\begin{align*}
&\widetilde{\Pi}_y^{\xi, -\lambda^*}
\left((I_3\cap I_4)^c\right)
=\widetilde{\Pi}_y^{\xi, -\lambda^*}
\left(\max_{3\leq k\leq [V_t]}\max_{j=1,2} \left|\beta_k^j - \beta_{k-1}^j\right| > C(\ell')\ln [V_t]\right)\\
 &	 \leq  \sum_{k=3}^{[V_t]}\sum_{j=1}^2\widetilde{\Pi}_y^{\xi, -\lambda^*} \left(\left|\beta_k^j - \beta_{k-1}^j\right| > C(\ell')\ln [V_t]\right)\\
&\leq  \sum_{k=3}^{[V_t]}\sum_{j=1}^2 e^{-\eta C(\ell')\ln[V_t]}\widetilde{\Pi}_y^{\xi, -\lambda^*}
\left(\exp\left\{\eta \left|\beta_k^j -\beta_{k-1}^j \right|\right\}\right)
\leq 2 c_1
[V_t]^{1-\eta C(\ell')}\cdot e^{\eta K_1/(\lambda^*v^*)},
\end{align*}
Taking $C(\ell')$ large so that $\eta C(\ell') -1 > \ell'$, we get that \eqref{step_20} holds for sufficiently large $t$.

(ii) For simplicity, we set $\gamma_0:=\gamma(-\lambda^*)-(m-1)\textup{ei}$. By the definition of $\beta_2^1$ and noting that $0\leq x_0:= \frac{1}{\lambda^*v^*}\ln \psi(
2, -\lambda^*) \leq 2K_1/(\lambda^* v^*)$,
for $x_1 \in [0, x]$,
 we have
\begin{align}\label{step_21}
& \widetilde{\Pi}_y^{\xi, -\lambda^*}
\left(\beta_2^1 \in [x_1, x_1+1]\right)
= \widetilde{\Pi}_y^{\xi, -\lambda^*}
\left(H_{2}\in [x_1+x_0,  x_1+x_0 +1] \right)\nonumber\\
&=  \Pi_y \left(\exp\left\{\int_0^{H_{2}} \left((m-1)\xi(B_s)-\gamma(-\lambda^*) \right)\mathrm{d}s\right\}; H_{2}\in [x_1+x_0,  x_1+x_0 +1]\right)\nonumber\\
&\geq  \Pi_y
\left(e^{-\gamma_0 H_{2}}; H_{2}\in [x_1+x_0,  x_1+x_0 +1]\right) \geq e^{-\gamma_0(x_1+x_0+1)} \Pi_0\left(H_{2-y}\in[x_1+x_0, x_1+ x_0+1]\right)\\
&  \geq e^{-\gamma_0(x+2K_1/(\lambda^* v^*)
+1)} \Pi_0\left(H_{2-y}\in[x_1+x_0 + 1/2, x_1+ x_0+1]\right).
\end{align}
Since $H_{2-y}$ has density $p_{H_
{2-y}}(s):= \frac{(2-y)}{\sqrt{2\pi s^3}}e^{-(2-y)^2/(2s)},$ we have
\begin{align*}
&\Pi_0 \left(H_{2-y}\in[x_1+x_0+1/2, x_1+x_0+1]\right)
= \int_{x_1+x_0+1/2}^{x_1+x_0+1}
\frac{(2-y)}{\sqrt{2\pi s^3}}e^{-(2-y)^2/(2s)}
\mathrm{d}s\\
& \geq \frac{1}{2} \frac{1}{\sqrt{2\pi \left(x_1+1+2K_1/(\lambda^* v^*)
\right)^3}}e^{-9}
 \geq \frac{1}{2} \frac{1}{\sqrt{2\pi \left(x+1+2K_1/(\lambda^* v^*)\right)^3}}e^{-9}.
\end{align*}
Since $z < e^z$  for $z>0$, we have $z^{-3/2}\ge e^{-3z/2}$ for $z>0$, and thus there exist
$c_2,c_3>0$
such that
$$
\inf_{y\in[-1,1]}\inf_{x_1 \in [0,x]}\widetilde{\Pi}_y^{\xi, -\lambda^*}
\left(\beta_1^1 \in [x_1, x_1+1]\right)
\geq c_2 e^{-c_3 x}, \quad x>0.
$$
Taking $C_3= \min\{ c_2, 1/c_3\}
$, we arrive at the desired conclusion.
\hfill$\Box$

The estimates of
$\widetilde{\Pi}_y^{\xi, -\lambda^*} \left(I_1\right)$ and $\widetilde{\Pi}_y^{\xi, -\lambda^*} \left(I_2\right)$ are given below.

\begin{lemma}\label{lemma 4.9}
There exists a non-random $\ell''>0$ such that for $t$ large,
\begin{align*}
\inf_{y\in[-1,1]}\widetilde{\Pi}_y^{\xi, -\lambda^*} \left(I_1\right)
\geq [V_t]^{-\ell''}; \quad
\inf_{y\in[-1,1]}\widetilde{\Pi}_y^{\xi, -\lambda^*} \left(I_2\right) \geq [V_t]^{-\ell''}.
\end{align*}
\end{lemma}

We postpone the proof of Lemma \ref{lemma 4.9}  to
Subsection \ref{proof of lemma 4.9}.

Now we are ready to prove the following estimate on  the first moment of $\mathcal{L}_t$:
\begin{lemma}\label{lemma 4.10}
There exists a non-random constant $\ell_1>0$ such that for $t$ large,
	\begin{equation}
		\inf_{y\in[-1,1]} \mathrm{E}_y^{\xi}(\mathcal{L}_t) \geq t^{-\ell_1}.
	\end{equation}
\end{lemma}
\textbf{Proof:}
Recall the definition of $x_1$ in \eqref{def-x1}.
By \eqref{step_29} and Lemma \ref{lemma 4.8}(ii),
\begin{align}
\widetilde{\Pi}_y^{\xi, -\lambda^*} \left(I_5 \big| \cap^4_{j=1}I_j	\right)
&\geq\frac{1}{[V_2]-2}
C_3e^{-2C(\ell')\ln[V_t]/C_3}
\geq
 C_3[V_t]^{-\frac{2C(\ell')}{C_3}-1}.
\end{align}
By \eqref{I1234}, Lemma \ref{lemma 4.8}(i) and Lemma \ref{lemma 4.9}, we have
\begin{equation}\label{step_30}
	\inf_{y\in[-1,1]} \widetilde{\Pi}_y^{\xi, -\lambda^*}
\left(\cap^4_{j=1}I_j\right)
\geq [V_t]^{-2 \ell ''} - [V_t]^{-\ell'}.
\end{equation}
Note that $\ell'>0$ in Lemma \ref{lemma 4.8} is arbitrary.
Letting $\ell' > 2 \ell''$ and using
\eqref{step_28}, \eqref{step_29} and \eqref{step_30},
we finally get that
$$
\inf_{y\in[-1,1]}\mathrm{E}_y^{\xi}(\mathcal{L}_t) \geq
\Gamma_3 C_3[V_t]^{-2C(\ell')/C_3-1-2 \ell ''}  \cdot
\left(1- [V_t]^{-(\ell'-2\ell'')}\right).
$$
Combining this with the fact that $V_t/t \to v^*$, and letting $\ell_1 > 2\ell'' + 2C(\ell')/C_3+1$, we get the desired result.
\hfill$\Box$

\begin{lemma}\label{domi-prob-number}
There exists a non-random constant
$C_4>1$
such that for sufficiently large $t$,
$$
\sup_{x\in\mathbb{R}} \mathrm{P}_x^\xi\left(\#\left\{\nu\in N(t): X_t^\nu \in [x-1, x+1] \right\} \leq C_4^t\right) \leq C_4^{-t}.
$$
\end{lemma}
\textbf{Proof:} See \cite[Lemma 6.8]{JD} for the case of continuous-time branching random walk in random environment or \cite[Lemma 4.7]{DS} for the case of branching Brownian motion in random environment.
\hfill$\Box$

\textbf{Proof of Theorem \ref{Log-distance}:}
Let $\lambda=-\lambda^*$ and let
$\Gamma_4$ be a constant such that $\Gamma_4 \ln C_4 > \ell_2 + 2\ell_1$, where $C_4 $ is the constant in Lemma \ref{domi-prob-number}.  For any $n\in \mathbb{N}$, define $r_n:=[\Gamma_4 \ln n]$. For simplicity, we will
omit the subscript $n$ from $r_n$ in this proof.
By \eqref{difference-V},  $V_t -V_s \leq \frac{\lambda^* v^*}{K_2}(t-s)$
for $t>s.$
Let $
\Gamma = \lambda^* v^* \Gamma_4/K_2+1$, then $V_n -V_{n-r} \leq (\Gamma -1)\ln n$ and
\begin{align}\label{step_38}
& \mathrm{P}_0^\xi \left(\inf_{n-1\leq t \leq n} M_t - V_n < -\Gamma \ln n \right)\nonumber\\
&\leq  \mathrm{P}_0^\xi\left(\#\left\{\nu\in N(r): |X_r^\nu|\leq 1 \right\} \leq C_4^r\right) + \sup_{|y|\leq 1}\left\{\mathrm{P}_y^\xi\left( \inf_{n-1\leq t \leq n} M_{t-r}- V_{n} < -\Gamma \ln n \right)\right\}^{C_4^{r}}\nonumber \\
&\leq  C_4 n^{-\Gamma_4\ln C_4} + \sup_{|y|\leq 1}\left\{\mathrm{P}_y^\xi\left( \inf_{n-1\leq t \leq n} M_{t-r}- V_{n-r} <-\ln n \right)\right\}^{n^{\Gamma_4  \ln C_4 }/C_4}.
\end{align}
Recall the definition of $\mathcal{L}_t$ given by \eqref{def-L}. For large $n$ such that $\ln n > 3,$ we have
\begin{align}\label{step_31}
\mathrm{P}_y^\xi\left( \inf_{n-1\leq t \leq n} M_{t-r}- V_{n-r} \geq -\ln n \right)& \geq \mathrm{P}_y^\xi\left( \inf_{n-1\leq t \leq n} M_{t-r}- V_{n-r} \geq -3 \right)\nonumber\\
\geq  &\mathrm{P}_y^{\xi}(\mathcal{L}_{n-r} > 0).
\end{align}
On the other hand, by Lemmas \ref{Second-moment} and \ref{lemma 4.10}, there exists $n_0$ such that for all $n\ge n_0$ and  $y\in[-1,1]$,
\begin{equation}\label{step_32}
 \mathrm{P}_y^\xi\left(\mathcal{L}_{n-r} > 0\right)\geq \frac{\left(\mathrm{E}_y^\xi\left(\mathcal{L}_{n-r}\right)\right)^2}{
 	\mathrm{E}_y^\xi\left((\mathcal{L}_{n-r})^2\right)} \geq (n-r)^{-\ell_2 -2\ell_1 }\geq n^{-\ell_2 -2\ell_1 }.
\end{equation}
Combining \eqref{step_38}, \eqref{step_31} and \eqref{step_32}, we get
$$
\mathrm{P}_0^\xi \left(\inf_{n-1\leq t \leq n} M_t - V_n < -\Gamma \ln n \right) \leq C_4 n^{-\Gamma_4 \ln C_4} + \left(1- n^{-\ell_2 -2\ell_1}\right)^{n^{\Gamma_4 \ln C_4 }/C_4}.
$$
Since $1-x\leq e^{-x}$ for $x>0$, we have
$$\mathrm{P}_0^\xi \left(\inf_{n-1\leq t \leq n} M_t - V_n < -\Gamma \ln n \right) \leq C_4 n^{-\Gamma_4 \ln C_4} + \exp \left\{-n^{\Gamma_4\ln C_4 - \ell_2 -2\ell_1}/C_4 \right\}.$$
Applying the Borel-Cantelli lemma, we get
$$\liminf_{n\to\infty}  \frac{\inf_{n-1\leq t \leq n} M_t - V_n}{\ln n} \geq -\Gamma ,\quad \mathrm{P}_0^\xi \mbox{-a.s.}$$
Since $|V_n - V_t|$ is uniformly bounded for $n$ and $t$ satisfying  $|n-t|\leq 1$, we
get that $\mathrm{P}_0^\xi$-almost surely,
$$\liminf_{t\to\infty} \frac{M_t - V_t}{\ln t}\geq -\Gamma , $$
which is equivalent to
$$ \limsup_{t\to\infty} \frac{V_t-M_t }{\ln t}\leq \Gamma .$$
By Theorem \ref{Additive martingale}, $\lim_{t\to\infty}W_t(-\lambda^*)=0$. Since
$W_t(-\lambda^*)\geq\psi(M_t,-\lambda^*)e^{-\gamma(\lambda^*)t},$
we have
$$
\lim_{t\to\infty}\psi(M_t,-\lambda^*)e^{-\gamma(\lambda^*)t}=0.$$
Therefore $ \lim_{t\to\infty}(\ln\psi(M_t,-\lambda^*))-\ln\psi(V_t,-\lambda^*))=\lim_{t\to\infty}(-\gamma(\lambda^*)t+\ln\psi(M_t,-\lambda^*))=-\infty$.
By Remark \ref{remark3},
$M_t - V_t \to -\infty,\ \mathrm{P}_0^\xi$-a.s. as $t\to\infty$,
this completes the proof.
\hfill$\Box$

\subsection{Proof of Theorem \ref{Invariance_of_M_t} }

In this subsection, we also assume that {\bf(H1)}, {\bf(H2)} and {\bf(H3)} hold.

\begin{lemma}\label{Nolen-CLT}
	Let $\kappa>0$. Suppose that $W(x,\omega): [0,\infty) \times \Omega \to \mathbb{R}$ is a continuous process on $(\Omega,\mathcal{F},\mathbb{P})$
such that,  $\mathbb{P}$-almost surely, $W(0,\omega) = 0$ and
	$$\lim_{x\to\infty} \frac{W(x,\omega)}{x} = 0.$$
Suppose further that the family of processes $\left\{n^{-1}W(nx,\omega)\right\}_{n=1}^\infty$ converges weakly as $n\to\infty$ to $\kappa B(x)$, where $B(x)$ is a standard Brownian motion on the interval $[0,M]$, for any $M$,  in the  topology of uniform convergence. Define
	$$h_t(\omega) = \sup\left\{h\geq -\lambda^* v^* t\big| \ W(h+\lambda^* v^* t) = h - 1 \right\}.$$
	Then as $t\to \infty$, the family of processes $\left\{\frac{1}{\kappa  \sqrt{\lambda^* v^* n}} h_{nt} \right\}_{n=1}^\infty$ converges weakly as $n\to\infty$ to a standard Brownian motion on $[0,M]$, in the Skorohod space.
\end{lemma}
\textbf{Proof:} See \cite[Lemma 3.1]{Nolen}.
\hfill$\Box$

\begin{lemma}\label{CLT_of_V_t}
Let  $\tilde{\sigma}_{-\lambda^*}^2$ be the constant defined in Theorem
\ref{Invariance_of_M_t}. Under $\mathbb{P}$, we have
	$$\frac{V_t - v^* t}{\sqrt{t}}\stackrel{\mathrm{d}}{\Rightarrow} \mathcal{N}\left(0, \tilde{\sigma}_{-\lambda^*}^2\right),\, \mbox{ as } t\to\infty.$$
	If $\tilde{\sigma}_{-\lambda^*}^2>0$, then the sequence of processes
	$$
	[0, \infty) \ni t \mapsto \frac{V_{nt}-v^* n t}{\tilde{\sigma}_{-\lambda^*} \sqrt{n}}, \quad n \in \mathbb{N},
	$$
	converges weakly as $n\to\infty$ to a standard Brownian motion on $[0,\infty)$, in the Skorohod topology.
\end{lemma}
\textbf{Proof:} By Lemma \ref{Engenfunction-Property1}(2) and
\eqref{Psi_and_Phi},
$V_t$ is the unique solution to
$\ln \psi(x,-\lambda^*) = \lambda^* v^* t.$ By Remark \ref{CLT2}, for any fixed $y\in\mathbb{R}$, we have
\begin{align*}
\mathbb{P}\left(\frac{V_t - v^* t}{\sqrt{t}} \leq y \right) & = \mathbb{P}\left(V_t \leq v^* t + y\sqrt{t}\right)= \mathbb{P}\left(\lambda^* v^* t \leq \ln \psi(v^*t +y\sqrt{t},-\lambda^*)\right)\\ & = \mathbb{P}\left( -\ln \phi(v^*t +y\sqrt{t},-\lambda^*) \leq \lambda^* y\sqrt{t}  \right) \\ & = \mathbb{P}\left(-\frac{\ln \phi(v^*t +y\sqrt{t},-\lambda^*)  }{\sqrt{v^*t +y\sqrt{t} }} \leq \frac{\lambda^*}{\sqrt{v^*}} y \cdot \sqrt{\frac{t}{t + y\sqrt{t}/v^*}}\right) \\ &\stackrel{t\to\infty}{\longrightarrow} \mathbb{P}\left((\sigma_{-\lambda^*}') \chi \leq \frac{\lambda^*}{\sqrt{v^*}}y\right) =\mathbb{P}((\tilde{\sigma}_{-\lambda^*} )\chi \leq y).
\end{align*}
Here $\chi$ is a standard normal random variable.

If $\tilde{\sigma}_{-\lambda^*}^2>0$,  let $W(x,\omega) : = -\ln \phi\left(\frac{x}{\lambda^*}, -\lambda^*,\omega\right)$ and $\kappa:= \sigma_{-\lambda^*}' /\sqrt{\lambda^*}> 0$.
It is easy to see that
$$h_t(\omega) = \lambda^*
\left(V_{t+ (\lambda^* v^*)^{-1}}(\omega) - v^* t\right).
 $$
Using Lemma \ref{Nolen-CLT}, we have
  that the sequence of processes
$$
[0, \infty) \ni t \mapsto \frac{h_{nt}}{\kappa  \sqrt{\lambda^* v^* n}} = \frac{V_{nt + (\lambda^* v^*)^{-1}}-v^* n t}{\tilde{\sigma}_{-\lambda^*} \sqrt{n}}, \quad n \in \mathbb{N},
$$
converges as $n \to \infty$ weakly to a standard Brownian motion on $[0, M]$, for any $M>0$, in the Skorohod topology. According to \cite[Lemma 3, p.173]{Billingsley}, this is equivalent to the weak convergence to a standard Brownian motion on $[0,\infty)$ in the Skorohod topology.
By Remark \ref{remark3},  $\big|V_{t+(\lambda^* v^*)^{-1}} - V_t\big|\leq 1/K_2$ all $t\geq 0$,  hence the second conclusion of the lemma is valid.
\hfill$\Box$

\textbf{Proof of Theorem \ref{Invariance_of_M_t}:} This follows from Theorem \ref{Log-distance} and Lemma \ref{CLT_of_V_t}.
\hfill$\Box$

\subsection{Proof of Lemma \ref{lemma 4.9}}\label{proof of lemma 4.9}

In this subsection we give the proof of Lemma \ref{lemma 4.9}. We first prove several lemmas.

\begin{lemma}\label{lemma 4.12}
There exists a non-random constant
$C_5>0$
 such that for $j=1,2$,
$$
\inf_{y\in[-1,1], k\in\N, k\geq 2}
\widetilde{\Pi}_y^{\xi, -\lambda^*} \left(\beta_k^j - \beta_{k-1}^j > C_5 \right) > C_5,
\qquad
\inf_{y\in[-1,1], k\in\N, k\geq 2}
\widetilde{\Pi}_y^{\xi, -\lambda^*} \left(\beta_k^j - \beta_{k-1}^j < - C_5 \right) > C_5.
$$
\end{lemma}
\textbf{Proof:}
Since $\beta^1$ and $\beta^2$ are identically distributed,  we only prove the case $j=2$.
Recall that $\beta_k^2 = H_{k}^2 -\frac{1}{\lambda^*v^*}
\ln \psi(k, -\lambda^*)$ and $0 < K_2 \leq \ln \psi(k, -\lambda^*) - \ln \psi(k
-1, -\lambda^*) \leq K_1$ for every $k\in\mathbb{Z}$
and $k\geq 2$.
Note that for any Borel set $B\subset [0, +\infty)$ with $m(B)>0$,
\begin{align*}
\widetilde{\Pi}_y^{\xi, -\lambda^*}\left(H_{k}^2 - H_{k-1
}^2 \in B\right)&  = \widetilde{\Pi}_{k-1}^{\xi, -\lambda^*}\left(H_{k
}^2 \in B\right) \geq \Pi_0 \left( e^{-H_1\left(\gamma(\lambda^*)- (m-1)\textup{ei}\right)}; H_1 \in B\right).
\end{align*}
For $0 < \delta< K_2/(\lambda^* v^*)$, we have
\begin{align*}
\widetilde{\Pi}_y^{\xi, -\lambda^*} \left(\beta_k^2 - \beta_{k-1}^2 < -\delta  \right) & \geq \widetilde{\Pi}_y^{\xi, -\lambda^*}\left(H_{k}^2 - H_{k-1}^2 < -\delta + \frac{K_2}{\lambda^* v^*}\right)
	\\ & \geq \Pi_0 \left( e^{-H_1\left(\gamma(\lambda^*)- (m-1)\textup{ei}\right)}; H_1 < -\delta + K_2/(\lambda^* v^*)\right).
\end{align*}
Similarly,
$$\widetilde{\Pi}_y^{\xi, -\lambda^*} \left(\beta_k^2 - \beta_{k-1}^2 > \delta  \right)\geq \Pi_0 \left( e^{-H_1\left(\gamma(\lambda^*)- (m-1)\textup{ei}\right)}; H_1 > \delta + K_1/(\lambda^* v^*)\right). $$
When $\delta= 0$, it holds that
\begin{align*}
&  \Pi_0 \left( e^{-H_1\left(\gamma(\lambda^*)- (m-1)\textup{ei}\right)}; H_1 <  K_2/(\lambda^* v^*)\right) > 0=\delta,\nonumber\\   &\Pi_0 \left( e^{-H_1\left(\gamma(\lambda^*)- (m-1)\textup{ei}\right)}; H_1 >  K_1/(\lambda^* v^*)\right) > 0=\delta.
\end{align*}
Since $H_1$ has a density,
the maps $\delta\mapsto \Pi_0 \left( e^{-H_1\left(\gamma(\lambda^*)- (m-1)\textup{ei}\right)}; H_1 < -\delta + K_2/(\lambda^* v^*)\right)$ and
$\delta\mapsto \Pi_0 \left( e^{-H_1\left(\gamma(\lambda^*)- (m-1)\textup{ei}\right)}; H_1 >
\delta+ K_1/(\lambda^* v^*)
\right)$
are continuous in $\delta > 0$. Thus there exists $C_5 > 0$ such that
\begin{align*}
&  \Pi_0 \left( e^{-H_1\left(\gamma(\lambda^*)- (m-1)\textup{ei}\right)}; H_1 < -C_5 + K_2/(\lambda^* v^*)\right) > C_5,\\   &\Pi_0 \left( e^{-H_1\left(\gamma(\lambda^*)- (m-1)\textup{ei}\right)}; H_1 > C_5 + K_1/(\lambda^* v^*)\right) > C_5.
\end{align*}
This completes the proof.
\hfill$\Box$

\begin{lemma}\label{Gaussian-scalling}
For any  $\kappa_0>0$,
there exist non-random constants
$N(\kappa_0)>0$ and
$C_6=C_6(\kappa_0)\in(0,1)$
such that for all
$r\geq 2, k\geq N(\kappa_0)$ and $j=1,2,$
\begin{equation*}
\inf_{y\in[-1,1]}\widetilde{\Pi}_y^{\xi, -\lambda^*}
\left(\inf_{r<m  \leq r+k} \left(H_{m}^j -H_{r }^j -\widetilde{\Pi}_y^{\xi, -\lambda^*} \left( H_{m
}^j -H_{r}^j\right)\right) \geq  -\kappa_0 \sqrt{k}\right)\geq
C_6.
\end{equation*}
\end{lemma}
\textbf{Proof:}
We only deal with $j=1$, the case $j=2$ is similar. Define
\begin{equation}\label{step_34}
\Delta H_k:= H_{k}^1 -H_{k-1}^1 -\widetilde{\Pi}_y^{\xi, -\lambda^*} \left( H_{k}^1 -H_{k-1 }^1\right)
\end{equation}
and
$$
f(\eta):= \ln \widetilde{\Pi}_y^{\xi, -\lambda^*} \left(\exp\left\{\eta (H_{k}^1 - H_{k-1}^1 )\right\} \right).
$$

Note that by \eqref{Change-of-measure},
\eqref{step_3-1'}
and \eqref{Psi_2-2}, we have
\begin{align*}
&\widetilde{\Pi}_y^{\xi, -\lambda^*} \left( H_{k
}^1 \right)=\Pi_y \left(H_{k}\exp\left\{-\gamma(-\lambda^*) H_{k} + \int_0^{H_{k}} (m-1)\xi(B_s)\mathrm{d}s \right\} \frac{\psi(k
,-\lambda^*)}{\psi(y,-\lambda^*)} \right) \\
& = \frac{\psi(k,-\lambda^*,\omega)}{\psi(y,-\lambda^*,\omega)} \cdot \frac{1}{-\gamma'(-\lambda^*)} \frac{\partial}{\partial \lambda}\Pi_y \left(\exp\left\{-\gamma(\lambda) H_{k
} + \int_0^{H_{k}} (m-1)\xi(B_s)\mathrm{d}s \right\} \right)\bigg|_{\lambda=-\lambda^*}\\
& = \frac{\psi(k,-\lambda^*,\omega)}{\psi(y,-\lambda^*,\omega)} \cdot \frac{1}{-\gamma'(-\lambda^*)} \frac{\partial}{\partial \lambda}\left(\frac{\psi(y,\lambda,\omega)}{\psi(k,\lambda,\omega)}\right)\bigg|_{\lambda=-\lambda^*}\\& = \frac{\psi(k
-y,-\lambda^*,\theta_y\omega)}{-\gamma'(-\lambda^*)}\frac{\partial}{\partial \lambda}\left(\frac{1}{\psi(
k-y,\lambda,\theta_y\omega)}\right)\bigg|_{\lambda=-\lambda^*} = \frac{1}{\gamma'(-\lambda^*)}\frac{\psi_\lambda(k-y,-\lambda^*,\theta_y\omega)}{\psi(
k-y,-\lambda^*,\theta_y\omega)}
\\ & = \frac{1}{\gamma'(-\lambda^*)}\left[\frac{\psi_\lambda(	k,-\lambda^*,\omega)}{\psi(	k,-\lambda^*,\omega)} - \frac{\psi_\lambda(	y,-\lambda^*,\omega)}{\psi(	y,-\lambda^*,\omega)}\right].
\end{align*}
Thus
\begin{align}
&\widetilde{\Pi}_y^{\xi, -\lambda^*} \left( H_{k}^1 \right)- \widetilde{\Pi}_y^{\xi, -\lambda^*} \left( H_{k-1}^1 \right) = \frac{1}{\gamma'(-\lambda^*)}\left(\frac{\psi_\lambda(k,-\lambda^*,
\omega)}{\psi(k,-\lambda^*,
\omega)} -\frac{\psi_\lambda(k-1,-\lambda^*,
\omega)}{\psi(k-1,-\lambda^*,
\omega)} \right).
\end{align}
Using  Lemma \ref{Engenfunction-Property1}(1), we obtain
\begin{equation}\label{differ-H}
\left|\widetilde{\Pi}_y^{\xi, -\lambda^*} \left( H_{k}^1 \right)- \widetilde{\Pi}_y^{\xi, -\lambda^*} \left( H_{k-1}^1 \right)\right|
\leq
\frac{C_1(-\lambda^*)}{\left|\gamma'(-\lambda^*)\right|},
\end{equation}
here $C_1(-\lambda^*)$ is the constant in Lemma \ref{Engenfunction-Property1}.

By \eqref{step_33}, we have,  uniformly for any $0 \leq \eta \leq \left(\gamma(\lambda^*) -(m-1)\textup{es}\right)/2$
and $k\geq 2$,
\begin{align*}
f'(\eta) &= \frac{\widetilde{\Pi}_y^{\xi, -\lambda^*} \left(\exp\left\{\eta (H_{k}^1 - H_{k-1}^1 )\right\} (H_{k}^1 - H_{k-1}^1 ) \right)}{\widetilde{\Pi}_y^{\xi, -\lambda^*} \left(\exp\left\{\eta (H_{k}^1 - H_{k-1}^1 )\right\} \right)}
\\ & = \frac{\widetilde{\Pi}_{k-1}^{\xi, -\lambda^*} \left(\exp\left\{\eta H_{k}^1 \right\} H_{k}^1  \right)}{\widetilde{\Pi}_{k-1}^{\xi, -\lambda^*} \left(\exp\left\{\eta H_{k}^1 \right\} \right)}
\leq \widetilde{\Pi}_{k-1}^{\xi, -\lambda^*} \left(\exp\left\{\eta H_{k}^1 \right\}H_{k}^1  \right)
\\ & \leq \Pi_0 \left(\exp\left\{-H_1 \left(\gamma(\lambda^*) - (m-1)\textup{es} -\eta \right) \right\}H_1\right)\cdot e^{K_1}\\
	&  \leq  \Pi_0 \left(\exp\left\{-H_1 \left(\gamma(\lambda^*) - (m-1)\textup{es} \right)/2 \right\}H_1\right)\cdot e^{K_1}
= :  c_1< +\infty.
\end{align*}
This implies that $f(\eta) \leq \eta c_1$ for any $0 \leq \eta \leq \left(\gamma(\lambda^*) -(m-1)\textup{es}\right)/2$.
Therefore, by \eqref{differ-H},
\begin{align*}
\widetilde{\Pi}_y^{\xi, -\lambda^*}\left(e^{\eta|\Delta H_k| }\right) &\leq \widetilde{\Pi}_y^{\xi, -\lambda^*} \left(\exp\left\{\eta (H_{k}^1 - H_{k-1}^1 )\right\} \right)\cdot \exp\left\{\eta \left|\widetilde{\Pi}_y^{\xi, -\lambda^*} \left( H_{k}^1 -H_{k-1 }^1\right) \right| \right\} \\ & \leq \exp\left\{\eta c_1 + \eta \left|C_1(-\lambda^*)/\gamma'(-\lambda^*)\right|\right\}.
\end{align*}
Thus, $\Delta H_k$ are sub-exponential random variables for $k\geq 1$. By \cite[Proposition 2.7.1]{RV},
there exists
$c_2$
depending only on
$c_1 +  C_1(-\lambda^*)/\left|\gamma'(-\lambda^*)\right|$
such that
\begin{align}\label{domi-MGF-above}
\widetilde{\Pi}_y^{\xi, -\lambda^*}\left(e^{\eta \Delta H_k }\right)\leq \exp\left\{( c_2\eta)^2\right\},\quad |\eta | \leq 1/c_2, k\geq 1.
\end{align}

By martingale theory and Doob's inequality, for $0 < \eta < 1/(2c_2)$,
\begin{align*}
&\widetilde{\Pi}_y^{\xi, -\lambda^*}\left(
\inf_{r < m  \leq r+k}
 \sum_{\ell=r+1}^m \Delta H_\ell  \geq  -\kappa_0 \sqrt{k}\right)  = 1- \widetilde{\Pi}_y^{\xi, -\lambda^*}\left(\sup_{r<m  \leq r+k} \sum_{\ell=r+1}^m (-\Delta H_\ell)>\kappa_0\sqrt{k}\right)\\
&=  1- \widetilde{\Pi}_y^{\xi, -\lambda^*}\left(
 \sup_{r< m  \leq r+k}
\exp\left\{ \eta \sum_{\ell=r+1}^m (-\Delta H_\ell) \right\}
>
  e^{\eta \kappa_0 \sqrt{k}}\right) \\
&\geq  1-e^{-2\eta \kappa_0 \sqrt{k}}\cdot \widetilde{\Pi}_y^{\xi, -\lambda^*}\left(\exp\left\{ 2\eta \sum_{\ell=r+1}^{r+k} (-\Delta H_\ell) \right\} \right)\\
&=   1-e^{-2\eta \kappa_0 \sqrt{k}}\cdot \prod_{\ell=r+1}^{r+k}\widetilde{\Pi}_y^{\xi, -\lambda^*}\left(\exp\left\{ 2\eta (-\Delta H_\ell) \right\} \right)
	\geq  1- \exp\left\{ 4(c_2 \eta)^2k-2\eta\kappa_0 \sqrt{k} \right\}.
\end{align*}
Taking $\eta : = \kappa_0 /(4c_2^2 \sqrt{k}) \leq 1/(2c_2)$
and $N(\kappa_0):= (\kappa_0/(2c_2))^2$,  then for $k \geq N(\kappa_0)$,
$$
\widetilde{\Pi}_y^{\xi, -\lambda^*}\left(
\inf_{r < m  \leq r+k}
 \sum_{\ell=r+1}^m \Delta H_\ell  \geq  -\kappa_0 \sqrt{k}\right)
\geq C_6(\kappa_0):= 1- \exp\left\{-\kappa_0^2/(4c_2^2) \right\}\in (0,1).
$$
The proof is now complete.
\hfill$\Box$

\begin{lemma}\label{Bernstein-type_ineq}
	Let $\Delta H_k$  be defined by \eqref{step_34}. There exist
positive constants $\alpha^*$ and
$C_7$,  independent of $k$ and $r$,
such that for all $0\leq x \leq \alpha^* k$ and $k,r \in \mathbb{Z}^+,$
$$
\inf_{y\in[-1,1]}\widetilde{\Pi}_y^{\xi, -\lambda^*} \left(\sum_{s= r+1}^{r+k} \Delta H_s \geq x \right)\geq C_7 \exp\left\{- \frac{x^2}{C_7 k}\right\}.
$$
\end{lemma}
\textbf{Proof:}
We continue the constant label in the proof of Lemma \ref{Gaussian-scalling}.
We claim that
\begin{equation}\label{domi-secondmoment}
c_3:=\textup{essinf}_{\xi, k} \widetilde{\Pi}_y^{\xi, -\lambda^*} \left((\Delta H_k)^2\right) > 0.
\end{equation}
In fact, for $x\in(0,1)$, let
$$
\Delta_x H_k:= H_{k-1+x}^1 -H_{k-1 }^1 -\widetilde{\Pi}_y^{\xi, -\lambda^*} \left( H_{k-1+x}^1 -H_{k-1}^1\right).
$$
By the strong Markov property, under $\widetilde{\Pi}_y^{\xi, -\lambda^*}$, $\Delta_x H_k$ and $\Delta H_k -\Delta_x H_k$ are independent. Since $\Delta H_k$  and $\Delta_x H_k$ are both centered, we have
\begin{align*}
&\widetilde{\Pi}_y^{\xi, -\lambda^*} \left((\Delta H_k)^2\right)  > \widetilde{\Pi}_y^{\xi, -\lambda^*} \left((\Delta_x H_k)^2\right)\\
& = {\Pi}_{k-1}\left( \exp\left\{\int_0^{H_{k-1+x}}
\left((m-1)\xi(B_s)-\gamma(-\lambda^*) \right)\mathrm{d}s \right\}\frac{\psi(k-1+x,-\lambda^*)}{\psi(k-1,-\lambda^*)}{(H_{k-1+x})^2}\right)  \\
&\qquad -\left(\widetilde{\Pi}_y^{\xi, -\lambda^*} \left( H_{k-1+x}^1 -H_{k-1}^1\right)\right)^2
	\\ & \geq \Pi_0 \left(  {(H_{x})}^2\exp\left\{-\left(\gamma(-\lambda^*)-(m-1)\textup{ei} \right){H_x }\right\}\right) -
e^{2K_1 x}
	 \left(\Pi_0 \left(  {H_{x}}\exp\left\{-\left(\gamma(-\lambda^*)-(m-1)\textup{es} \right)
{H_x} \right\}\right)\right)^2
	\\ & = x^2 \left(\frac{1}{2u_2}e^{-x\sqrt{2u_2}} -\frac{1}{2u_1}
	e^{-2x(\sqrt{2u_1}-K_1)}
	\right)+\frac{x}{2\sqrt{2u_2^3}}e^{-x\sqrt{2u_2}}=:\widehat{F}(x),
\end{align*}
where $u_1:=\gamma(-\lambda^*)-(m-1)\textup{es}\leq \gamma(-\lambda^*)-(m-1)\textup{ei}=:u_2,$ and in the last equality we used \eqref{Exp-of-H}. Since $\lim_{x\to0} \frac{\widehat F(x)}{x}>0$, there exists $x$ such that for
all $k\in\mathbb{N}$ and all realizations of $\xi$,
$\widetilde{\Pi}_y^{\xi, -\lambda^*} \left((\Delta H_k)^2\right)> \widehat F(x)>0$,
which implies the claim holds.

By the sub-exponential property of $\Delta H_k$ {given by \eqref{domi-MGF-above}},
there exists a constant $K>0$,
independent of $k$, such that
for all $\ell \geq 1$, $\widetilde{\Pi}_y^{\xi, -\lambda^*} \left(\left|\Delta H_k\right|^{\ell}\right)\leq
 (K \ell)^\ell
$
(see \cite[Proposition 2.7.1, (ii)]{RV}). Therefore, by the trivial inequality that $\ell ! \geq (\ell /e)^{\ell},$ when $|\eta|< 1/(2e K),$
\begin{align*}
&\widetilde{\Pi}_y^{\xi, -\lambda^*} \left(e^{\eta \Delta H_k}\right) = 1+ \frac{1}{2!}\widetilde{\Pi}_y^{\xi, -\lambda^*} \left((\Delta H_k)^2\right)\eta^2 +
\sum_{\ell=3}^\infty \frac{1}{\ell!}\widetilde{\Pi}_y^{\xi, -\lambda^*}\left(\left(\Delta H_k\right)^{\ell}\right) |\eta|^\ell  \\
& \geq 1 + \frac{c_{3}}{2}\eta^2
- \sum_{\ell=3}^\infty \frac{1}{(\ell /e)^{\ell}}
(K \ell|\eta|)^\ell
= 1 +\frac{c_{3}}{2}\eta^2
- \sum_{\ell=3}^\infty
(eK |\eta|)^\ell\\
& = 1 + \frac{c_{3}}{2}\eta^2
- \frac{(eK)^3 |\eta|^3}{1-eK| \eta|}\geq 1+ \left(\frac{c_{3}}{2}- \frac{(eK)^3 |\eta|}{2}\right)\eta^2,
\end{align*}
where $c_3$ is given by \eqref{domi-secondmoment}.
Now we choose
$c_4> 0$ and $c_5$
such that
$ (eK)^3 c_4= c_3/2$
and that $1+c_3c_4^2/4 = e^{c_5c_4^2}$, then for all $|\eta|{ \leq} c_4$, we have
$$
1+ \left(\frac{c_3}{2}-
\frac{(eK)^3| \eta|}{2}\right)\eta^2
 \geq 1+ \frac{c_3}{4}\eta^2\geq e^{c_5 \eta^2},
 $$
here in the last inequality we used the fact that the function $h(x):= e^{c_5x}-1-c_3x/4$ is non-positive  for $x\in (0, c_4^2)$, which follows easily from our choice of $c_5$ and $c_4$.
Thus we have
\begin{equation}\label{domi-MGF-below}
\widetilde{\Pi}_y^{\xi, -\lambda^*} \left(e^{\eta \Delta H_k}\right) \geq e^{c_5 \eta^2},\quad {|\eta| \leq c_4, k\geq 1}.
\end{equation}
Now, using  \eqref{domi-MGF-above} and \eqref{domi-MGF-below},
the desired conclusion follows from \cite[Theorem 4]{ZZ} by  taking
$M={c_4\wedge(1/2c_2)}, \alpha =1, C_1 =(c_2)^2, c_1=C_5, c''= {M}^2, u_1=u_2=\cdots=u_k=1$ and $\alpha^* = c' M\alpha  $ with $c'>0$ being the  constant in \cite[Theorem 4]{ZZ},
 and $c_2$ is the constant in \eqref{domi-MGF-above}.
\hfill$\Box$

\bigskip

For $\ell \in \mathbb{Z}$, let
$$
\rho_\ell := \frac{\ln \psi(\ell,-\lambda^*)-\ln \psi(\ell-1,-\lambda^*)}{\lambda^* v^*}
+
 \frac{1}{v^*} \left(\frac{\psi_\lambda(\ell,-\lambda^*)}{\psi(\ell,-\lambda^*)} - \frac{\psi_\lambda(\ell-1,-\lambda^*)}{\psi(\ell-1,-\lambda^*)} \right).
$$
By \eqref{relation-lnpsi-lnphi},
\begin{equation}\label{def-rho}
\rho_\ell= \frac{\ln \phi(\ell,-\lambda^*)-\ln \phi(\ell-1,-\lambda^*)}{\lambda^* v^*}
+
 \frac{1}{v^*} \left(\frac{\phi_\lambda(\ell,-\lambda^*)}{\phi(\ell,-\lambda^*)} - \frac{\phi_\lambda(\ell-1,-\lambda^*)}{\phi(\ell-1,-\lambda^*)} \right).
\end{equation}

Taking  logarithm and differentiating with respect to $\lambda$ in  \eqref{Psi_2-2},  and
 letting $\lambda =-\lambda^*$, we get
\begin{align*}
\lambda^* v^*\rho_\ell =& - \ln \Pi_{\ell-1}\left(\exp\left\{ \int_0^{H_{\ell}}\left((m-1)\xi(B_s) -\gamma(-\lambda^*) \right)\mathrm{d}s \right\}\right)\\
	& -\lambda^* v^* \frac{\Pi_{\ell-1}\left(H_{\ell}\exp\left\{ \int_0^{H_{\ell}}\left((m-1)\xi(B_s) -\gamma(-\lambda^*) \right)\mathrm{d}s \right\}\right)}{\Pi_{\ell-1}\left(\exp\left\{ \int_0^{H_{\ell}}\left((m-1)\xi(B_s) -\gamma(-\lambda^*) \right)\mathrm{d}s \right\}\right)}.
\end{align*}
Define
\begin{align*}
A_1(-\lambda^*)
 &=\Pi_{\ell-1}\left(\exp\left\{ \int_0^{H_{\ell}}\left((m-1)\xi(B_s) -\gamma(-\lambda^*) \right)\mathrm{d}s \right\}; \inf_{0\leq s\leq H_{\ell}} B_s  > j \right),\\
A_2(-\lambda^*) &=\Pi_{\ell-1}\left(\exp\left\{ \int_0^{H_{\ell}}\left((m-1)\xi(B_s) -\gamma(-\lambda^*) \right)\mathrm{d}s \right\}; \inf_{0\leq s\leq H_{\ell}} B_s  \leq j \right),\\
A_1'(-\lambda^*) &=\Pi_{\ell-1}\left(H_{\ell}\exp\left\{ \int_0^{H_{\ell}}\left((m-1)\xi(B_s) -\gamma(-\lambda^*) \right)\mathrm{d}s \right\}; \inf_{0\leq s\leq H_{\ell}} B_s  > j \right),\\
A_2'(-\lambda^*) &=\Pi_{\ell-1}\left(H_{\ell}\exp\left\{ \int_0^{H_{\ell}}\left((m-1)\xi(B_s) -\gamma(-\lambda^*) \right)\mathrm{d}s \right\}; \inf_{0\leq s\leq H_{\ell}} B_s  \leq j \right).
\end{align*}
Then
$$\rho_\ell=-\frac{1}{\lambda^* v^*}\ln\left(A_1(-\lambda^*) + A_2(-\lambda^*)\right)
-
 \frac{A_1'(-\lambda^*)+ A_2'(-\lambda^*)}{A_1(-\lambda^*) + A_2(-\lambda^*)}.
$$
Note that
$A_1(-\lambda^*), A_1'(-\lambda^*) \in  \mathcal{F}_{\ell
} \cap \mathcal{F}^j$.
Using the same notation in \eqref{domi-A1},
\begin{align*}
	&A_1(-\lambda^*) \leq \Pi_{\ell-1}\left(\exp\left\{ \left((m-1)\textup{es} -\gamma(-\lambda^*) \right)H_{\ell} \right\}\right)\\&= \Pi_{0}\left(\exp\left\{ \left((m-1)\textup{es} -\gamma(-\lambda^*) \right)H_{1} \right\}\right) =c_1(-\lambda^*)<1,
\end{align*}
and when $\ell\geq j+
2,$ we have $j-(\ell-1)\leq -1$ and thus
\begin{align*}
&A_1(-\lambda^*) \geq \Pi_{\ell-1}\left(\exp\left\{ \left((m-1)\textup{ei} -\gamma(-\lambda^*) \right)H_{\ell} \right\}; \inf_{0\leq s\leq H_{\ell}} B_s  > j \right)\\
& =  \Pi_{0}\left(\exp\left\{ \left((m-1)\textup{ei} -\gamma(-\lambda^*) \right)H_{1} \right\}; \inf_{0\leq s\leq H_{1}} B_s  > j - (\ell -1) \right) \\
& \geq \Pi_{0}\left(\exp\left\{ \left((m-1)\textup{ei} -\gamma(-\lambda^*) \right)H_{1} \right\}; \inf_{0\leq s\leq H_{1}} B_s  > -1 \right)= c_2(-\lambda^*) >0.
\end{align*}
If $H_{\ell }\geq \ell -j$, then by {\bf (H2)},
$$
\exp\left\{ \int_0^{H_{\ell}}\left((m-1)\xi(B_s) -\gamma(-\lambda^*) \right)\mathrm{d}s \right\} \leq \exp\left\{\left((m-1)\textup{es} -\gamma(-\lambda^*)\right) (\ell -j) \right\};
$$
and if $H_{\ell }\leq \ell -j$, then
\begin{align*}
	\left\{ \inf_{0\leq s\leq H_{\ell}} B_s  \leq j \right\} \subset \left\{\inf_{0\leq s\leq \ell
		-j} B_s  \leq j\right\}.
\end{align*}
Letting
$\gamma_0:=\gamma(-\lambda^*) -(m-1)\textup{es}>0$, when $\ell \geq j+2$, by the reflection principle and Markov's inequality,
\begin{align*}
&A_2(-\lambda^*) \leq e^{-\gamma_0(\ell -j
)}+ \Pi_{\ell-1} \left(\inf_{0\leq s\leq \ell
-j} B_s  \leq j \right)
= e^{-\gamma_0(\ell -j)}+\Pi_{0} \left(\inf_{0\leq s\leq \ell -j} B_s  \leq j- (\ell -1) \right) \\
& = e^{-\gamma_0(\ell -
)}+ 2\Pi_{0} \left(B_{\ell-j}  \leq j- \ell +1 \right) \\
&  \leq e^{-\gamma_0(\ell -j)}+ 2\frac{\Pi_0(e^{-B_{\ell
-j}/2})}{e^{-(j-\ell+1)/2}}= e^{-\gamma_0(\ell -j
)}+ 2\frac{e^{(\ell-j)/8}}{e^{-(j-\ell+1)/2}} \leq (1 + 2e^{1/2})
e^{-(\ell -j)\delta}
\end{align*}
with $\delta:= \min\{3/8, \gamma_0\}.$

Now we take a constant $\lambda_0$ such that  $-\rho>\lambda_0 > -\lambda^*$ and $\gamma(-\lambda^*)>\gamma(\lambda_0)> (m-1)\textup{es}$. Then
there is a constant
$c>0$
such that $\sup_{x>0}\left\{xe^{-(\gamma(-\lambda^*)- \gamma(\lambda_0))x}\right\}\leq c$,
and then
\begin{align*}
&0\leq A_1'(-\lambda^*) \leq c\Pi_{\ell-1}\left(\exp\left
\{\int_0^{H_{\ell}}\left((m-1)\xi(B_s) -\gamma(\lambda_0) \right)\mathrm{d}s \right\};
\inf_{0\leq s\leq H_{\ell}} B_s  > j \right)=  c{A}_1(\lambda_0),\\
&0\leq A_2'(-\lambda^*) \leq c\Pi_{\ell-1}\left(\exp\left
 \{ \int_0^{H_{\ell}}\left((m-1)\xi(B_s) -\gamma(\lambda_0) \right)\mathrm{d}s \right\};
 \inf_{0\leq s\leq H_{\ell}} B_s \leq  j \right)=c{A}_2(\lambda_0).
\end{align*}
Noticing that $\gamma(\lambda_0)> (m-1)\textup{es}$,
using argument similar to those used for $A_1$ and $A_2$, with $-\lambda^*$ replaced by $\lambda_0$,  we have that there
exist positive constants $\tilde{\delta}, c_3(\lambda_0), c_4$ such that
$${A}_1(\lambda_0) < c_3(\lambda_0),
\quad {A}_2(\lambda_0) \leq
c_4e^{-(\ell-j)\tilde{\delta}},\quad \forall\ell \geq j+2.
$$
Therefore, setting $A_i= A_i(-\lambda^*), A_i'= A_i'(-\lambda^*)$ for simplicity, we conclude that
\begin{align*}
&\left|\rho_\ell + \frac{1}{\lambda^* v^*}\ln(A_1) +
\frac{A_1'}{ A_1}\right| = \left|\frac{-\ln(A_1 + A_2)-
	\lambda^*v^* \frac{A_1'+ A_2'}{A_1 + A_2}}{\lambda^* v^*} + \frac{1}{\lambda^* v^*}\ln(A_1)
+\frac{A_1'}{ A_1}\right|\\
&\leq \frac{1}{\lambda^* v^*}\left| \ln (1+ A_2/A_1)\right| + \frac{A_1 A_2' +A_2 A_1'}{(A_1 +A_2)A_1} \\
& \leq \frac{1}{\lambda^* v^*}\frac{A_2}{A_1}+ \frac{c_1(-\lambda^*) cc_4e^{-(\ell -j)\tilde{\delta}}
+(1 + 2e^{1/2})e^{-(\ell -j)\delta}c{c}_3(\lambda_0) }{c_2(-\lambda^*)^2}\\
	& \leq  \frac{1}{\lambda^* v^*}\frac{(1 + 2e^{1/2})
		e^{-(\ell-j)\delta}}{c_2(-\lambda^*)}+ \frac{c_1(-\lambda^*)cc_4e^{-(\ell-j)\tilde{\delta}}+(1 + 2e^{1/2})e^{-(\ell -j)\delta}c
c_3(\lambda_0) }{c_2(-\lambda^*)^2}.
\end{align*}
Thus there exists some constant $C_8>0$ such that for
$\forall\ell \geq j+2$,
\begin{equation}\label{rho-barrho}
\left|\rho_\ell + \frac{1}{\lambda^* v^*}\ln(A_1)
+\frac{A_1'}{ A_1}\right|\leq C_8e^{-(\ell-j)\delta'}
\end{equation}
with $\delta':= \min \{\delta, \tilde{\delta}\}>0$.
Define
$$
\overline{\rho}_\ell^{(j)}: = -\frac{1}{\lambda^* v^*}\ln(A_1(-\lambda^*))
-\frac{A_1'(-\lambda^*)}{ A_1(-\lambda^*)}+ C_8 e^{-(\ell-j)\delta'}.
$$
Then by \eqref{rho-barrho}, for $\ell \geq j+2$, we have that
\begin{align}\label{step_36}
0\leq \overline{\rho}_\ell^{(j)} - \rho_\ell  =& \left(\overline{\rho}_\ell^{(j)} + \frac{1}{\lambda^* v^*}\ln(A_1(-\lambda^*))
+
\frac{A_1'(-\lambda^*)}{ A_1(-\lambda^*)} \right)\nonumber \\
&-\left({\rho}_\ell + \frac{1}{\lambda^* v^*}\ln(A_1(-\lambda^*))
+\frac{A_1'(-\lambda^*)}{  A_1(-\lambda^*)} \right) \nonumber\\
\leq &2C_8e^{-(\ell-j)\delta'}.
\end{align}
Also, $|\overline{\rho}_\ell^{(j)}|$ is uniformly bounded for all $\ell \geq j+2$:
$$
 \Vert \overline{\rho}_\ell^{(j)}\Vert_{\infty}
 \leq \frac{1}{\lambda^* v^*} \ln \frac{1}{c_2(-\lambda^*)} + \frac{cc_3(\lambda_0)}{c_2(-\lambda^*)}
+ C_8 =: M_\rho.
 $$
By the definitions of $A_1(-\lambda^*)$ and  $ A_1'(-\lambda^*)$,  we know that
\begin{align}\label{Measurable-of-rho}
\overline{\rho}_\ell^{(j)} \in \mathcal{F}_{\ell
}\cap \mathcal{F}^j \quad \mbox{when} \quad \ell \geq j+2.
\end{align}

\bigskip

For $ q \geq 0$ and $i\geq 1$, let $t_i:= 2^i$ and define
\begin{align}\label{r-s-L-K}
	r_0  &: = t_{i-1},\quad L_i  :=\left\lceil \frac{C_5 t_i^{1/2}}{16 M_\rho } \right\rceil,
\nonumber \\
	s_0& := t_i \land \inf\left\{
k\geq r_0+1:
\sum_{\ell= r_0+1}^k \overline{\rho}_\ell^{(r_0-[L_i/2])} \geq C_5 t_i^{1/2}/16\right\}, \quad r_{q+1}  : = s_q + L_i,\nonumber \\
	s_{q+1}& : = \inf\left\{
k\geq r_{q+1}+1:
 \sum_{\ell= r_{q+1}+1}^k \overline{\rho}_\ell^{(r_{q+1}-[L_i/2])} \geq C_5 t_i^{1/2}/16\right\} \land \left(r_{q+1} + t_i-t_{i-1}\right),
	\nonumber\\
\mathfrak{K} & := \inf \{ k : s_k \geq t_i\},
\end{align}
here $C_5$ is the constant defined in Lemma \ref{lemma 4.12}.

For $\ell < t_i - t_{i-1}$, define
$$
D(\ell):= \left\{\omega:\, \sum_{r= 1}^k \overline{\rho}_r^{(-[L_i/2])} < C_5 t_i^{1/2}/16, \forall k< \ell, \quad\sum_{r= 1}^\ell \overline{\rho}_r^{(-[L_i/2])} \geq C_5 t_i^{1/2}/16 \right\}.
$$
Also define
$$
D(t_i - t_{i-1}) := \left\{\omega: \sum_{r= 1}^\ell \overline{\rho}_r^{(-[L_i/2])} < C_5 t_i^{1/2}/16, \, \forall \ell < t_i -t_{i-1} \right\}.
$$

\begin{prop}
Fix $i \geq 1$.

(i)	For any $k\geq 1$ and positive integers $\ell_0,...,\ell_{k-1}< t_i-t_{i-1}, \ell_k \leq t_i - t_{i-1}$,
\begin{align}\label{Mixing1-2}
&\mathbb{P}\left(s_j -r_j =\ell_j, 0\leq j\leq k \right)
\leq \left(1+ \zeta\left(C_5 t_i^{1/2}/(32M_\rho)
\right)\right)^{k}\prod_{j=0}^{k}	\mathbb{P}\left(D(\ell_j)\right).
\end{align}

(ii) For any $k\geq 1$ and $\ell_0,...,\ell_k < t_i - t_{i-1}$, when $i$ is sufficiently large so that $t_{i-1}> L_i$,
\begin{align}\label{Measurable-5}
\left\{\mathfrak{K} =k, s_j -r_j= \ell_j , 0\leq j \leq k \right\} \in \mathcal{F}^{t_{i-1}-[L_i/2]}\cap \mathcal{F}_{t_{i+1}}.
\end{align}
\end{prop}
\textbf{Proof: }
(i)
By induction, we only need to prove that
\begin{align}\label{Mixing1-1}
&\mathbb{P}\left(s_j -r_j = \ell_j, 0\leq j\leq k-1,
s_k -r_k = \ell_k \right)\nonumber \\
& \leq \mathbb{P}\left(s_j -r_j = \ell_j, 0\leq j\leq k-1 \right)\mathbb{P}\left(D(\ell_k)\right)\left(1+ \zeta\left(C_5 t_i^{1/2}/(32M_\rho)\right)\right).
\end{align}
For integers $r,s$ such that $r<s$ and $s-r < t_i-t_{i-1}$, define
\begin{align}
D(r,s):= \left\{\omega:\, \sum_{\ell= r+1}^k \overline{\rho}_\ell^{(r-[L_i/2])} < C_5 t_i^{1/2}/16,
\,\forall k\in[r+1,s),
\sum_{\ell= r+1}^s \overline{\rho}_\ell^{(r-[L_i/2])} \geq C_5 t_i^{1/2}/16 \right\}.
\end{align}
Also, when $s-r = t_i - t_{i-1}$, we define
\begin{align}
	D(r,s):= \left\{\omega:\, \sum_{\ell= r+1}^k \overline{\rho}_\ell^{(r-[L_i/2])} < C_5 t_i^{1/2}/16,
\,\forall k\in[r+1, r + t_i - t_{i-1})
\right\}.
\end{align}
Then if $\ell_j < t_i - t_{i-1}$ for all $0\leq j\leq k$,  it is easy to see that
\begin{align}\label{step_56}
&	\left\{s_j -r_j = \ell_j, 0\leq j\leq k\right\}\nonumber \\
 =&D(t_{i-1}, t_{i-1} +\ell_0) \cap D(t_{i-1}+\ell_0 + L_i, t_{i-1}+\ell_0 + L_i+ \ell_1) \nonumber\\
	& \cap\cdots \cap D(t_{i-1}+k L_i +\ell_0+...+\ell_{k-1},t_{i-1}+k L_i +\ell_0+...+\ell_{k} ).
\end{align}
It follows from \eqref{Measurable-of-rho} that $D(r,s) \in \mathcal{F}_{s}\cap \mathcal{F}^{r-[L_i/2]}$ when $s-r \leq t_i - t_{i-1}$. Using \eqref{step_56}, we get
\begin{align}\label{Mesurable-4}
	&\left\{s_j -r_j =\ell_j ,0\leq j \leq k-1 \right\}\in \mathcal{F}_{ t_{i-1}+(k-1) L_i +\ell_0+...+\ell_{k-1}},\nonumber \\
	& \{s_k -r_k = \ell_k\}\in \mathcal{F}^{t_{i-1}+k L_i +\ell_0+...+\ell_{k-1} - [L_i/2]}.
\end{align}
By the definition of $L_i$,
\begin{align*}
&\left(t_{i-1}+k L_i +\ell_0+...+\ell_{k-1} - [L_i/2]\right) - \left( t_{i-1}+(k-1) L_i +\ell_0+...+\ell_{k-1} \right)\\
&= L_i - [L_i/2]
\geq C_5 t_i^{1/2}/(32M_\rho),
\end{align*}
together with the $\zeta$-mixing condition we conclude that
\begin{align}
	&\mathbb{P}\left(s_j -r_j = \ell_j, 0\leq j\leq k\right) \\
&\leq \mathbb{P}\big(s_j -r_j = \ell_j, 0\leq j \leq k-1\big) \\
	& \quad \times \mathbb{P}\big(D(t_{i-1}+k L_i +\ell_0+...+\ell_{k-1},t_{i-1}+k L_i +\ell_0+...+\ell_{k} )\big)\left(1+ \zeta\left(C_5 t_i^{1/2}/(32M_\rho)\right)\right)\\
	& = \mathbb{P}\left(s_j -r_j = \ell_j, 0\leq j\leq k-1\right) \mathbb{P}(D(\ell_k))\left(1+ \zeta\left(C_5 t_i^{1/2}/(32M_\rho)
	\right)\right),
\end{align}
where in the last equality, we used the property that $\mathbb{P}(D(r,s))=\mathbb{P}(D(0,s-r))= \mathbb{P}(D(s-r))$.
Thus \eqref{Mixing1-1} is valid.

(ii)
Note that
\begin{align*}
		&\left\{\mathfrak{K} =k, s_j -r_j= \ell_j , 0\leq j \leq k \right\}  = \left\{s_k \geq t_i, s_{k-1}< t_i, s_j -r_j= \ell_j , 0\leq j \leq k \right\}\\
		& = \left\{s_j - r_j = \ell_j, 0\leq j \leq k, s_{k-1}= t_{i-1}+ (k-1)L_i + \sum_{j=0}^{k-1}\ell_j < t_i, s_k =t_{i-1}+ kL_i + \sum_{j=0}^{k}\ell_j \geq t_i \right\}.
\end{align*}
The event in \eqref{Measurable-5} is not  empty only if $t_{i-1}+ (k-1)L_i + \sum_{j=0}^{k-1}\ell_j < t_i$ and $t_{i-1}+ kL_i + \sum_{j=0}^{k}\ell_j \geq t_i$.  Combining this observation with \eqref{Mesurable-4},  we get that, when $t_{i-1}> L_i$,
\begin{align*}
	&\left\{\mathfrak{K} =k, s_j -r_j= \ell_j , 0\leq j \leq k \right\} \in \left(\bigcup_{j=0}^k \mathcal{F}^{t_{i-1}+j L_i +\ell_0+...+\ell_{j-1} - [L_i/2]}\right) \bigcap \mathcal{F}_{ t_{i-1}+k L_i +\ell_0+...+\ell_k}\\
	& \subset \mathcal{F}^{t_{i-1}-[L_i/2]}\cap \mathcal{F}_{L_i+\ell_{k}+t_i}\subset \mathcal{F}^{t_{i-1}-[L_i/2]}\cap \mathcal{F}_{L_i+(t_i-t_{i-1})+t_i} \subset \mathcal{F}^{t_{i-1}-[L_i/2]}\cap \mathcal{F}_{t_{i+1} }.
\end{align*}
The proof is now complete.
\hfill$\Box$

\bigskip

Put
\begin{equation}\label{def-alpha}
\alpha_1:= 2\alpha^*/(3C_5),\quad \alpha_2:= \alpha^*/(6C_5),
\end{equation}
where $\alpha^*$  is the constant in Lemma \ref{Bernstein-type_ineq},
and $C_5$  the constant in Lemma \ref{lemma 4.12}.
Define
$$
G_i^1:= \left\{\omega:\, \max_{r\in [t_{i-1}, t_{i+1}],
 0<k\leq  t_i^{1/2} /\alpha_1}\
\sum_{\ell=r+1}^{r+k} \overline{\rho}_\ell^{(t_{i-2})} \leq C_5t_{i-1}^{1/2}/32 \right\}.
$$

\begin{lemma}\label{lemma 4.15}	
(i) There exists a constant
$C_9>0$
such that for sufficiently large $i$,
on $G_i^1$,
\begin{align} \label{step_42}
&\inf_{y\in[-1,1]}
\inf_{x\geq C_5 t_{i-1}^{1/2}/2} 	\widetilde{\Pi}_y^{\xi, -\lambda^*}
\left(\beta_{t_i}^1-\beta_2^1
\geq C_5 {t_i}^{1/2}/2,\ \beta_k^1 -\beta_2^1 \geq t_i^{1/4}, k=t_{i-1}+1,...,t_i \Big| \beta_{t_{i-1}}^1-\beta_2^1 =x\right)\nonumber\\
& \geq  C_9^{\mathfrak{K}+1}
\exp\left\{- \sum_{j=0}^{\mathfrak{K}}\frac{t_i}{C_9(s_j - r_j)} \right\} =: Z_i^1,
	\end{align}

(ii) There exist constants $\varepsilon > 0$  and
$K=K(\epsilon)$
such that for $i$ large enough,
\begin{equation}
\mathbb{E}\left(\exp\left\{- \varepsilon 1_{G_i^1}\ln Z_i^1 \right\}\right) \leq K.
\end{equation}
\end{lemma}
\textbf{Proof:}
(i)
{\bf Step 1}
Recall  the definition of $\beta^1_k$ in \eqref{def-beta^j}.
For $j \geq 0$, define
\begin{align*}
	\mathcal{A}_j & :=
	\left\{\beta_{s_j}^1 -\beta_2^1 \geq C_5 t_i^{1/2}  \right\} \bigcap \bigcap_{k=r_j+1}^{s_j} \left\{ \beta_k^1 -\beta_{r_j}^1 +\sum_{\ell = r_j+1}^k \rho_\ell
\geq -\frac{1}{16}C_5 t_i^{1/2} \right\}
	\\
	\widetilde{\mathcal{A}}_j  &
	:= \bigcap_{k=s_j+1}^{r_{j+1}}\left\{\beta_k^1 -\beta_{s_j}^1 + \sum_{\ell = s_j+1}^k \rho_\ell
\geq -\frac{1}{16}C_5t_i^{1/2} \right\},\quad
	\mathcal{B}  : =\left( \bigcap_{j=0}^\mathfrak{K} \mathcal{A}_j \right)\bigcap \left(\bigcap_{j=0}^{\mathfrak{K}-1} \widetilde{\mathcal{A}}_j \right).
\end{align*}
Put
\begin{align*}
	\left(Z_i^1\right)^*:=&\inf_{x\geq C_5 t_{i-1}^{1/2}} 	\widetilde{\Pi}_y^{\xi, -\lambda^*}  \left(\beta_{t_i}^1-\beta_2^1
 \geq C_5 {t_i}^{1/2}/2,\
\beta_k^1 -\beta_2^1 \geq t_i^{1/4},
 k=t_{i-1}+1,...,t_i \Big| \beta_{t_{i-1}}^1-\beta_2^1 =x\right)\\
	 = &\widetilde{\Pi}_y^{\xi, -\lambda^*}  \left(\beta_{t_i}^1-\beta_2^1
\geq C_5 {t_i}^{1/2}/2,\
\beta_k^1 -\beta_2^1 \geq t_i^{1/4},
 k=t_{i-1}+1,...,t_i \Big| \beta_{t_{i-1}}^1-\beta_2^1
  =C_5 t_{i-1}^{1/2}/2\right),
\end{align*}
here the last equality follows from the following argument: since $\{\beta_k^1, k\geq 1\}$ has independent increments,
\begin{align*}
&	\widetilde{\Pi}_y^{\xi, -\lambda^*}  \left(\beta_{t_i}^1-\beta_2^1 \geq C_5 {t_i}^{1/2}/2,\ \beta_k^1 -\beta_2^1 \geq t_i^{1/4},
k=t_{i-1}+1,...,t_i \Big| \beta_{t_{i-1}}^1-\beta_2^1 =x\right)\\
&= \widetilde{\Pi}_y^{\xi, -\lambda^*}  \left(\beta_{t_i}^1-\beta_{t_{i-1}}^1 \geq C_5 {t_i}^{1/2}/2-x,\ \beta_k^1 -\beta_{t_{i-1}}^1 \geq t_i^{1/4}-x, k=t_{i-1}+1,...,t_i \Big| \beta_{t_{i-1}}^1-\beta_2^1 =x\right)\\
	& = 	\widetilde{\Pi}_y^{\xi, -\lambda^*}  \left(\beta_{t_i}^1-\beta_{t_{i-1}}^1 \geq C_5 {t_i}^{1/2}/2 -x,\ \beta_k^1 -\beta_{t_{i-1}}^1 \geq t_i^{1/4} -x,
k=t_{i-1}+1,...,t_i\right)
\end{align*}
and the last quantity is decreasing in $x$.

 In this step we will show that
\begin{equation}\label{Claim2}
 \left(Z_i^1\right)^*\geq	\widetilde{\Pi}_y^{\xi, -\lambda^*}  \left(\mathcal{B}\, \Big|\beta_{t_{i-1}}^1-\beta_2^1 =C_5 t_{i-1}^{1/2} /2\right).
\end{equation}
Given that $\mathcal{B}$ occurs, under $\widetilde{\Pi}_y^{\xi, -\lambda^*}  (\ \cdot  \ |\beta_{t_{i-1}}^1-\beta_2^1 =C_5 t_{i-1}^{1/2} /2)$, we have,  for $0\leq j \leq \mathfrak{K}$ and $r_j \leq k \leq s_j$,
$$
\beta_k^1 -\beta_{r_j}^1 +\sum_{\ell = r_j+1}^k \rho_\ell  \geq -C_5 t_i^{1/2}/16
\quad \mbox{ on } \mathcal{A}_j.
$$
By the definition of $s_j$ we have that
for all $r_j < k < s_j$,
$$
\sum_{\ell = r_j+1}^{k} \rho_\ell \leq \sum_{\ell = r_j+1}^{k} \overline{\rho}_\ell^{(r_j-[L_i/2])} < C_5 t_i^{1/2} /16,
$$
and for $k = s_j$,
$$
\sum_{\ell = r_j+1}^{s_j} \rho_\ell  \leq \sum_{\ell = r_j+1}^{s_j-1} \rho_\ell + M_\rho \leq  M_\rho +  C_5 t_i^{1/2} /16.
$$
Thus,  for $0\leq j \leq \mathfrak{K}$ and $r_j \leq k \leq s_j$, we have
\begin{equation}\label{domi-below-beta}
	\beta_k^1 -\beta_{r_j}^1 \geq -C_5 t_i^{1/2}/8 -M_\rho.
\end{equation}
Note that
\begin{equation}\label{domi-below-beta0}\beta_{r_0}^1 - \beta_2^1 =\beta_{t_{i-1}}^1 - \beta_2^1= C_5 t_{i-1}^{1/2} /2= C_5t_i^{1/2}/(2\sqrt{2}).
\end{equation}
Also for $\mathfrak{K}-1\geq j \geq 0$,
on $\widetilde{\mathcal{A}}_j $ we have
\begin{equation}\label{step_60}
	\beta_{r_{j+1}}^1 -\beta_{s_j}^1 \geq - \sum_{\ell = s_j+1}^{r_{j+1}} \rho_\ell -C_5 t_i^{1/2}/16 \geq -\left\lceil \frac{C_5 t_i^{1/2}}{16M_\rho } \right\rceil M_\rho -C_5 t_i^{1/2}/16 \geq -C_5 t_i^{1/2}/8- M_\rho
\end{equation}
and on $\mathcal{A}_j$, $\beta_{s_j}^1 -\beta_2^1 \geq  C_5 t_i^{1/2} $. Hence, together with \eqref{step_60}, we have
\begin{equation}\label{domi-below-beta1}
	\beta_{r_j}^1 - \beta_2^1 = (\beta_{r_j}^1 - \beta_{s_{j-1}}^1) + (\beta_{s_{j-1}}^1- \beta_2^1) \geq 7C_5t_i^{1/2}/8 -M_\rho,\quad  \mathfrak{K} \geq j\geq 1.
\end{equation}
So, combining \eqref{domi-below-beta}, \eqref{domi-below-beta0} and \eqref{domi-below-beta1}, we have, for $ 0\leq j \leq \mathfrak{K}$ and $r_j \leq k \leq s_j$,
$$
\beta_k^1 -\beta_2^1 \geq \min\left\{C_5t_i^{1/2}/(2\sqrt{2}), 7C_5t_i^{1/2}/8 -M_\rho \right\} - C_5 t_i^{1/2}/8-M_\rho.
$$
Let $i$ be sufficiently large so that
$$\min\left\{C_5t_i^{1/2}/(2\sqrt{2}), 7C_5t_i^{1/2}/8 -M_\rho \right\} - C_5 t_i^{1/2}/8-M_\rho > t_i^{1/4}.
$$
Then $\beta_k^1 - \beta_2^1 \geq t_i^{1/4}$ for all $0\leq j \leq \mathfrak{K}$ and $r_j \leq k \leq s_j$ when $i$ is sufficiently large.
Similarly, when $0\leq j \leq \mathfrak{K}-1$ and $s_j < k < r_{j+1}$, on $\widetilde{\mathcal{A}}_j $,
$$\beta_k^1 -\beta_{s_j}^1 \geq - \sum_{\ell = s_j+1}^k \rho_\ell -C_5 t_i^{1/2}/16\geq -\sum_{\ell = s_j+1}^{r_{j+1}} M_\rho -C_5 t_i^{1/2}/16\geq -C_5 t_i^{1/2}/8-M_\rho,$$
which  implies $\beta_k^1 -\beta_2^1 \geq 7C_5 t_i^{1/2}/8-M_\rho$ since on $\mathcal{A}_j$, $\beta_{s_j}^1 -\beta_2^1 \geq  C_5 t_i^{1/2} $. For large $i$ such that $7C_5 t_i^{1/2}/8-M_\rho > t_i^{1/4}$, we also have $\beta_k^1 -\beta_2^1 \geq t_i^{1/4}$ holds for all $0\leq j\leq \mathfrak{K}-1$ and $s_j< k< r_{j+1}$.
We have also proved that, for all $s_0\leq k \leq s_{\mathfrak{K}}$,
$$\beta_k^1 - \beta_2^1 \geq 7C_5t_i^{1/2}/8 -M_\rho - C_5 t_i^{1/2}/8-M_\rho.$$
 Finally, to show $\beta_{t_i}^1 - \beta_2^1 \geq C_5 t_i^{1/2}/2$ holds for large $i$, it suffices to prove that $s_0 \leq t_i\leq s_{\mathfrak{K}}$, which is trivial by the definitions of $s_j$ and $\mathfrak{K}$,.
Thus \eqref{Claim2} is valid.

{\bf Step 2} In this step we will show that
there exists a constant
$C_9>0$ such that
\begin{equation}\label{e:s1}
\inf_{y\in[-1,1]}
\widetilde{\Pi}_y^{\xi, -\lambda^*}  \left(\mathcal{B}\, \Big|\beta_{t_{i-1}}^1-\beta_2^1 =C_5 t_{i-1}^{1/2} /2\right) \geq C_9^{\mathfrak{K}+1}\exp\left\{- \sum_{j=0}^{\mathfrak{K}}\frac{t_i}{C_9(s_j - r_j)} \right\}.
\end{equation}
By the strong Markov property,
\begin{align}\label{step_50}
 & \widetilde{\Pi}_y^{\xi, -\lambda^*}\left(\mathcal{B} \Big|\beta_{t_{i-1}}^1-\beta_2^1 =C_5 t_{i-1}^{1/2}/2 \right)\nonumber \\
&\geq \widetilde{\Pi}_y^{\xi, -\lambda^*} \left(\mathcal{A}_0 \Big|\beta_{r_0}^1-\beta_2^1 =C_5 t_{i-1}^{1/2} /2\right) \cdot \widetilde{\Pi}_y^{\xi, -\lambda^*}\left(\widetilde{\mathcal{A}}_0\, \Big| \beta_{s_0}^1 -\beta_2^1 = C_5 t_i^{1/2}\right) \nonumber\\ &\quad \times \widetilde{\Pi}_y^{\xi, -\lambda^*}\left(\mathcal{A}_1 \,\Big| \beta_{r_1}^1 -\beta_2^1 = 7C_5 t_i^{1/2}/8 - M_\rho\right)\times ... \times \widetilde{\Pi}_y^{\xi, -\lambda^*} \left(\mathcal{A}_\mathfrak{K}\, \Big|\beta_{r_\mathfrak{K}}^1-\beta_2^1 =7C_5 t_i^{1/2}/8 - M_\rho \right)\nonumber \\
&=  \widetilde{\Pi}_y^{\xi, -\lambda^*} \left(\mathcal{A}_0 \Big|\beta_{r_0}^1-\beta_2^1 =C_5 t_{i-1}^{1/2}/2 \right) \prod_{j=0}^{\mathfrak{K}-1}\widetilde{\Pi}_y^{\xi, -\lambda^*}(\widetilde{\mathcal{A}}_j ) \prod_{j=1}^\mathfrak{K} \left(\mathcal{A}_j\Big| \beta_{r_j}^1 -\beta_2^1 =7C_5 t_i^{1/2}/8 - M_\rho  \right).
\end{align}
Note that $\widetilde{\mathcal{A}}_j $ can be rewritten as
$$
\widetilde{\mathcal{A}}_j   := \left\{
\inf_{s_j < k\leq r_{j+1} }
\left(H_{k}^1 -H_{s_j }^1 -\widetilde{\Pi}_y^{\xi, -\lambda^*} \left( H_{k}^1 -H_{s_j }^1\right) \right) \geq -C_5 t_i^{1/2}/16\right\}.
$$
This is because
\begin{align*}
&\beta_k^1 -\beta_{s_j}^1 + \sum_{\ell = s_j+1}^k \rho_\ell
= \sum_{\ell =s_j +1}^k \left(\beta_{\ell}^1 - \beta_{\ell-1}^1 + \rho_\ell\right) \\ &= \sum_{\ell =s_j +1}^k \left(H_{\ell}^1 - H_{\ell-1}^1
+ \frac{1}{v^*} \left(\frac{\psi_\lambda(\ell,-\lambda^*)}{\psi(\ell,-\lambda^*)} - \frac{\psi_\lambda(\ell-1,-\lambda^*)}{\psi(\ell-1,-\lambda^*)} \right) \right)\\
	& = \sum_{\ell =s_j +1}^k \left(H_{\ell}^1 - H_{\ell-1}^1
	-\frac{\Pi_{\ell-1}\left(H_{\ell}\exp\left\{ \int_0^{H_{\ell}}\left((m-1)\xi(B_s) -\gamma(-\lambda^*) \right)\mathrm{d}s \right\}\right)}{\Pi_{\ell-1}\left(\exp\left\{ \int_0^{H_{\ell}}\left((m-1)\xi(B_s) -\gamma(-\lambda^*) \right)\mathrm{d}s \right\}\right)} \right)\\
	& = \sum_{\ell =s_j +1}^k \left(H_{\ell}^1 - H_{\ell-1}^1
	-\widetilde{\Pi}_{\ell-1}^{\xi, -\lambda^*} \left( H_{\ell}^1 \right)  \right) = \sum_{\ell =s_j +1}^k \left(H_{\ell}^1 - H_{\ell-1}^1
	-\widetilde{\Pi}_{y}^{\xi, -\lambda^*} \left( H_{\ell}^1- H_{\ell-1}^1 \right)  \right) .
\end{align*}
Also note that $r_{j+1}- s_j < C_5 t_i^{1/2}/(16 M_\rho)+1 < t_i$ for sufficiently large $i$ and $r_{j+1}- s_j \geq C_5 t_i^{1/2}/(16 M_\rho)$.
Taking $\kappa_0 =C_5/16$
in Lemma \ref{Gaussian-scalling},
$C_6:=C_6(\kappa_0)\in (0,1)$ with $\kappa_0=C_5/16$
such that for large $i$,
\begin{equation}\label{domi-tilde-A}
\widetilde{\Pi}_y^{\xi, -\lambda^*} \left(\widetilde{\mathcal{A}}_j \right)\geq
C_6.
\end{equation}
Define $d_0=C_5 t_{i-1}^{1/2}/2$ and $d_j = 7C_5 t_i^{1/2}/8 - M_\rho$ for $j\geq 1$.  Let $f(X)= x_1 +...+ x_{s_j -r_j}, g(X) = \inf_{1\leq k\leq s_j-r_j} (x_1+...+x_k)$ and $X= (x_1,...,x_{s_j-r_j}).$ By Harris' inequality, see \cite[Theorem 2.15]{BLM},  we have
for all $j\geq0$,
\begin{align}\label{domi-expect-Aj}
& \widetilde{\Pi}_y^{\xi, -\lambda^*}\left(\mathcal{A}_j\Big| \beta_{r_j}^1 -\beta_2^1 =d_j  \right) \geq \widetilde{\Pi}_y^{\xi, -\lambda^*} \left( \beta_{s_j}^1 -\beta_2^1 \geq C_5 t_i^{1/2} \Big| \beta_{r_j}^1 -\beta_2^1 =d_j \right)\nonumber\\
	& \quad \times \widetilde{\Pi}_y^{\xi, -\lambda^*}\left(\inf_{r_j<k \leq s_j} \left(H_{k}^1 -H_{r_j }^1 -\widetilde{\Pi}_y^{\xi, -\lambda^*} \left( H_{k}^1 -H_{r_j }^1\right)\right) \geq  -C_5 t_i^{1/2}/16\right) \nonumber\\
	 &=  \widetilde{\Pi}_y^{\xi, -\lambda^*} \left( \beta_{s_j}^1 -\beta_{r_j}^1 \geq C_5 t_i^{1/2}-d_j\right)\nonumber\\
	& \quad \times \widetilde{\Pi}_y^{\xi, -\lambda^*}\left(\inf_{r_j<k \leq s_j} \left(H_{k}^1 -H_{r_j }^1 -\widetilde{\Pi}_y^{\xi, -\lambda^*} \left( H_{k}^1 -H_{r_j }^1\right)\right) \geq  -C_5 t_i^{1/2}/16\right)\nonumber\\
	 &\geq \widetilde{\Pi}_y^{\xi, -\lambda^*} \left( \beta_{s_j}^1 -\beta_{r_j}^1 \geq C_5 t_i^{1/2}\right)\nonumber\\
	& \quad \times \widetilde{\Pi}_y^{\xi, -\lambda^*}\left(\inf_{r_j<
		 k \leq s_j} \left(H_{k}^1 -H_{r_j }^1 -\widetilde{\Pi}_y^{\xi, -\lambda^*} \left( H_{k}^1 -H_{r_j }^1\right)\right) \geq  -C_5 t_i^{1/2}/16\right)\nonumber\\
&=: (I)\times (II).
\end{align}

We  treat $(I)$ first.
Note that $C_5 t_i^{1/2}/16 +M_\rho \leq C_5 t_i^{1/2}/2$ for large $i$. By the  definition of $s_j$, we have
$\sum_{\ell = r_j+1}^{s_j} \rho_\ell \leq C_5 t_i^{1/2}/2$ and
\begin{align}\label{step_49}
(I)=	&\widetilde{\Pi}_y^{\xi, -\lambda^*} \left( \beta_{s_j}^1 -\beta_{r_j}^1 + \sum_{\ell=r_j+1}^{s_j} \rho_\ell \geq C_5 t_i^{1/2}+\sum_{\ell=r_j+1}^{s_j} \rho_\ell \right)\nonumber\\
	\geq&\widetilde{\Pi}_y^{\xi, -\lambda^*} \left(H_{s_j}^1 -H_{r_j }^1 -\widetilde{\Pi}_y^{\xi, -\lambda^*} \left( H_{s_j}^1 -H_{r_j }^1\right)\geq  3C_5t_i^{1/2}/2\right).
\end{align}
By the definitions of $s_j$ and $r_j$, we have $t_i \geq s_j -r_j \geq C_5 t_i^{1/2}/(16M_\rho)$.
By \eqref{step_36} and the definition of $s_j$, for large $i$ such that
$2C_8 t_i e^{- C_5 t_i^{1/2} \delta' /(32 M_\rho)} < C_5 t_i^{1/2}/32$,
when $s_j-r_j<t_i-t_{i-1},$
\begin{align*}
	C_5 t_i^{1/2}/16 & \leq \sum_{\ell= r_{j}+1}^{s_j} \overline{\rho}_\ell^{(r_j-[L_i/2])} \leq \sum_{\ell= r_{j}+1}^{s_j} {\rho}_\ell + 2 C_8(s_j -r_j)\exp\left\{-\left(r_j+1-\left(r_j- [L_i/2]\right)\right) \delta '\right\}\\
	& \leq \sum_{\ell= r_{j}+1}^{s_j} {\rho}_\ell + 2C_8 t_i e^{-\left([L_i/2]+1\right) \delta'}\leq \sum_{\ell= r_{j}+1}^{s_j} {\rho}_\ell + 2C_8 t_i e^{- C_5 t_i^{1/2} \delta' /(32 M_\rho)}\\ & \leq \sum_{\ell= r_{j}+1}^{s_j} {\rho}_\ell + C_5 t_i^{1/2}/32\leq \sum_{\ell= r_{j}+1}^{s_j} \overline{\rho}_\ell^{(t_{i-2})} + C_5 t_i^{1/2}/32.
\end{align*}
This implies that,
on $G_i^1$,
$s_j -r_j \geq t_i^{1/2}/\alpha_1$.
When $s_j-r_j=t_i -t_{i-1}$,
we can also take $i$ large so that
$s_j -r_j=t_i-t_{i-1}> t_i^{1/2}/\alpha_1$.
By the definitions of $\alpha_1$ (see \eqref{def-alpha}) and $\alpha^*$, we have $\alpha^* (s_j - r_j) \geq 3 C_5 t_i^{1/2}/2$ and hence by
Lemma \ref{Bernstein-type_ineq} (with $k=s_j-r_j$ and $x=3C_5 t_i^{1/2}/2$) and \eqref{step_49},
the constant $C_7$ defined in Lemma \ref{Bernstein-type_ineq} satisfies that on $G_i^1$ and for any $j$,
\begin{equation}\label{domi-(I)}
     (I)\geq C_7\exp\left\{ - \frac{9C_5^2t_i}{4C_7(s_j - r_j)} \right\}.
\end{equation}

Now we treat $(II)$.
Taking $\kappa_0 = C_5/16$ in Lemma \ref{Gaussian-scalling}, we  get
\begin{equation}\label{step_39}
(II)=	\widetilde{\Pi}_y^{\xi, -\lambda^*}\left(
\inf_{r_j < k \leq s_j}
\left(H_{k}^1 -H_{r_j }^1 -\widetilde{\Pi}_y^{\xi, -\lambda^*} \left( H_{k}^1 -H_{r_j }^1\right)\right) \geq  -C_5 t_i^{1/2}/16\right) \geq C_6,
\end{equation}
where $C_6=C_6(\kappa_0)$ with $\kappa_0=C_5/16$.
If $s_j = r_j + t_i-t_{i-1}$, then $s_j -r_j \geq C_5t_i^{1/2}/(16M_\rho)$ for large $i$ and \eqref{step_39} also holds.
Plugging \eqref{domi-(I)} and \eqref{step_39}  into \eqref{domi-expect-Aj}, we obtain
\begin{equation}\label{domi-expect-Aj2}
	\widetilde{\Pi}_y^{\xi, -\lambda^*}\left(\mathcal{A}_j\Big| \beta_{r_j}^1 -\beta_2^1 =d_j  \right) \geq
C_7\cdot C_6
\exp\left\{ - \frac{9C_5^2t_i}{4C_7(s_j - r_j)} \right\}.
\end{equation}
Combining \eqref{step_50}, \eqref{domi-tilde-A} and \eqref{domi-expect-Aj2}, we conclude that
on $G_i^1$,
\begin{align*}
	Z_i^1&\geq C_7 \cdot
	C_6(C_7\cdot (C_6)^2)^{\mathfrak{K}}
 \exp\left\{- \sum_{j=0}^{\mathfrak{K}}\frac{9C_5^2t_i}{4C_7(s_j - r_j)}\right\}
 \\  &\geq (C_7\cdot (C_6)^2)^{\mathfrak{K}+1}
 \exp\left\{- \sum_{j=0}^{\mathfrak{K}}\frac{9C_5^2t_i}{4C_7(s_j - r_j)} \right\}.
\end{align*}
Since $\mathfrak{K}$ and $s_j-r_j$ does not depend on $y\in[-1,1]$, \eqref{e:s1} holds with $C_9:=\left(C_7\cdot C_6^2\right)\wedge\left(4C_7/(9C_5^2)\right)$.

Combining  Step 1 and Step 2, we get \eqref{step_42}.

(ii)
Define
$$
\mathfrak{J}: = \inf\{j: s_j-r_j \geq t_i - t_{i-1} \}.
$$
Then $\mathfrak{K} \leq \mathfrak{J}$. In fact,  on the event $\mathfrak{K}=k$, we have $s_k \geq t_i$ but $s_{k-1}< t_i$, and then for all $0\leq j< k,$ $s_j - r_j \leq s_{k-1}-t_{i-1}<t_i - t_{i-1}$, which implies that $\mathfrak{J} \geq k$.
For $\varepsilon\in(0, 1]$, by \eqref{step_42},
\begin{align}
&\mathbb{E}\left(\exp\left\{- \varepsilon
1_{G_i^1}
\ln Z_i^1 \right\}\right) \leq \mathbb{E}\left( \exp\left\{\varepsilon
1_{G_i^1}
\left(-\left(\mathfrak{J} +1\right)
\ln C_9 + \sum_{j=0}^{\mathfrak{J}}\frac{t_i}{C_9(s_j - r_j)} \right) \right\}\right) \nonumber\\
&\leq  \mathbb{E}\left( \exp\left\{\varepsilon  \left(-\left(\mathfrak{J}+1\right)\ln C_9 + \sum_{j=0}^{\mathfrak{J}}\frac{t_i}{C_9(s_j - r_j)} \right) \right\}\right)\nonumber \\
		& \leq  \sum_{k=0}^\infty \left(\frac{1}{C_9 }\right)^{\varepsilon (k+1)} \mathbb{E}\left(\exp\left\{\varepsilon\frac{t_i}{C_9(s_k - r_k)} + \varepsilon  \sum_{j=0}^{{k-1}}\frac{t_i}{C_9(s_j - r_j)} \right\} 1_{\{\mathfrak{J} = k\}} \right).
\end{align}
By \eqref{Mixing1-2},
\begin{align*}
& \mathbb{E}\left(\exp\left\{\varepsilon\frac{t_i}{C_9(s_k - r_k)} + \varepsilon  \sum_{j=0}^{{k-1}}\frac{t_i}{C_9(s_j - r_j)} \right\} 1_{\{\mathfrak{J} = k\}} \right)\\
&=  \mathbb{E}\Bigg(\exp \left\{\varepsilon\frac{t_i}{C_9(s_k - r_k)} \right\} 1_{\{s_k-r_k = t_i-t_{i-1}\}} \prod_{j=0}^{k-1} \exp\left\{ \varepsilon  \frac{t_i}{C_9(s_j - r_j)} \right\} 1_{\{s_j -r_j < t_i -t_{i-1} \}}\Bigg) \\
	& \leq e^{2\varepsilon /C_9} \sum_{1\leq \ell_0,...,\ell_{k-1} < t_i-t_{i-1}} \prod_{j=0}^{k-1} \exp\left\{ \varepsilon  \frac{t_i}{C_9\ell_j} \right\} \mathbb{P}\left(s_j-r_j = \ell_j, 0\leq j \leq k-1\right)\\
	&  \leq \left(1+\zeta\left(C_5 t_i^{1/2}/(32M_\rho)
	\right)
\right)^{k}
e^{2\varepsilon /C_9} \sum_{1\leq \ell_0,...,\ell_{k-1} < t_i-t_{i-1}} \prod_{j=0}^{k-1} \left(\exp\left\{ \varepsilon  \frac{t_i}{C_9\ell_j} \right\} \mathbb{P}\left(D(\ell_j)\right)\right)\\
& = e^{2\varepsilon /C_9}\left(1+\zeta\left(C_5 t_i^{1/2}/(32M_\rho)
\right)
\right)^{k} \left(\sum_{\ell = 1}^{t_i - t_{i-1}-1}\exp\left\{ \varepsilon  \frac{t_i}{C_9\ell} \right\} \mathbb{P}\left(D(\ell)\right)\right)^k.
\end{align*}
Therefore,
\begin{align}\label{domi-exp-G-Z}
&\mathbb{E}\left(\exp\left\{- \varepsilon
1_{G_i^1}
\ln Z_i^1 \right\}\right) \nonumber\\
&\leq  \sum_{k=0}^\infty
\frac{e^{2\varepsilon /C_9}}{(C_9)^{\varepsilon (k+1)} }
\left(1+\zeta\left(C_5 t_i^{1/2}/(32M_\rho)
\right)
\right)^{k}
\left(\sum_{\ell = 1}^{t_i - t_{i-1}-1}\exp\left\{ \varepsilon  \frac{t_i}{C_9\ell} \right\} \mathbb{P}\left(D(\ell)\right)\right)^k.
\end{align}

Note that $t_i - t_{i-1} = t_i/2$, and $D(\ell)$ are disjoint for all $\ell< t_i/2$.
Therefore, for all $q<t_i/2$,
\begin{align}\label{step_57}
	\sum_{\ell = 1}^{q}\mathbb{P}(D(\ell))= \mathbb{P}\left(\bigcup_{\ell=1}^q D(\ell)\right)= \mathbb{P}\left(\max_{1\leq k\leq q} \sum_{ \ell = 1}^k \overline{\rho}_\ell^{(-[L_i/2])}\geq C_5 t_i^{1/2}/16\right).
\end{align}
Using \eqref{step_36}, for any $q  < t_i/2$ and large $i$ such that $C_5 t_i^{1/2}/16 -2C_8e^{-([L_i/2]+1) \delta'}\cdot t_i/2 \geq C_5 t_i^{1/2}/32$, we have
\begin{align}\label{step_54}
&	\mathbb{P}\left(\max_{1\leq k\leq q} \sum_{ \ell = 1}^k \overline{\rho}_\ell^{(-[L_i/2])}\geq C_5 t_i^{1/2}/16\right)\nonumber  \leq \mathbb{P}\left(\max_{1\leq k\leq q} \sum_{ \ell = 1}^k {\rho}_\ell \geq C_5 t_i^{1/2}/16 -2C_8e^{-(1+[L_i/2]) \delta'}\cdot q  \right)\\ &  \leq \mathbb{P}\left(\max_{1\leq k\leq q } \sum_{ \ell = 1}^k {\rho}_\ell \geq C_5 t_i^{1/2}/32\right).
\end{align}
Recall that $\rho_\ell \in \mathcal{F}_{\ell
}$ and that $\mathbb{E}(\rho_\ell)=0$, so for all $k\leq \ell
$, by $\zeta$-mixing, we have that
$\left| \mathbb{E}\left(\rho_k| \mathcal{F}^\ell\right)\right|\leq \mathbb{E}(|\rho_k|)\zeta(\ell-k)$ and therefore, for all $\ell \geq 1$,
\begin{align*}
	&\sup_{1\leq j\leq \ell} \left(\left\Vert \rho_j^2\right\Vert_{\infty} + 2 \left\Vert \rho_\ell \sum_{k=j}^{\ell-1}
	\mathbb{E} \left(\rho_k| \mathcal{F}^\ell\right)
	 \right\Vert_{\infty}\right)\\
	& \leq M_\rho^2 + \sup_{1\leq j\leq \ell} 2\Vert\rho_\ell\Vert_{\infty}\sum_{k=j}^{\ell-1} \mathbb{E}(|\rho_k|)\zeta(\ell-k)\leq M_\rho^2 + 2M_\rho^2 \sum_{k=0}^\infty \zeta(k)<\infty.
\end{align*}
By \cite[Theorem 2.4]{Rio}, or more generally,  \cite[Theorem 1]{KM},
there exists a constant
$c_1$
 such that for any $x>0$ and $q \in \mathbb{Z}^+$,
\begin{equation}\label{step_35}
	\mathbb{P}\left(\max_{1\leq k\leq q}\sum_{\ell=1}^k \rho_\ell >x\right)\leq c_1e^{-x^2/(c_1q)}.
\end{equation}
Using \eqref{step_57}, \eqref{step_54} and \eqref{step_35}, we conclude that, when $i$ is large enough, for all $q< t_i/2$,
\begin{align}\label{step_58}
	\sum_{\ell = 1}^{q}\mathbb{P}(D(\ell)) \leq \mathbb{P}\left(\max_{1\leq k\leq q } \sum_{ \ell = 1}^k {\rho}_\ell \geq C_5 t_i^{1/2}/32\right) \leq c_1e^{-C_5^2 t_i/(32^2 c_1 q)}
= : c_1e^{-c_2t_i /q}.
\end{align}
Also, for $q = t_i/2-1$, we have that
$$
\sum_{\ell = 1}^{t_i/2 -1}\mathbb{P}(D(\ell)) \leq \mathbb{P}\left(\max_{1\leq k\leq t_i/2 } \sum_{ \ell = 1}^k {\rho}_\ell \geq C_5 t_i^{1/2}/32\right).
$$
Recall the definition of $\rho_l$ given by \eqref{def-rho}.
If $\sigma_{-\lambda^*}^2 > 0$, then by Lemma \ref{Engenfunction-Property2} and the continuous mapping theorem,
$$
 \mathbb{P}\left(\max_{1\leq k\leq t_i/2} \sum_{ \ell= 1}^k {\rho}_\ell \geq C_5 t_i^{1/2}/32\right)
\stackrel{i \to\infty}{\longrightarrow}\Pi_0\left(\sup_{0\leq t\leq 1/2} \sigma_{-\lambda^*} B_t/(\lambda^* v^*) \geq C_5/32\right)<1.
$$
Otherwise, if $\sigma_{-\lambda^*}^2 =0$, by Lemma \ref{Engenfunction-Property2}, when $C_5 t_i^{1/2}/32 > \sup_{k\in\mathbb{Z}_+}\Vert 
\sum_{ \ell = 1}^k
 {\rho}_\ell\Vert_\infty$, we have
\begin{align*}
\mathbb{P}\left(\max_{1\leq k\leq t_i/2} \sum_{ \ell = 1}^k
 {\rho}_\ell \geq C_5 t_i^{1/2}/32\right)=0.
\end{align*}
In conclusion,
there exists a constant
$c_3 \in (0,1)$
 such that for large $i$,
\begin{equation}\label{step_44}
\sum_{\ell = 1}^{t_i/2 -1}\mathbb{P}(D(\ell))  \leq  \mathbb{P}\left(
\max_{1\leq k\leq t_i/2} \sum_{ \ell = 1}^k
 {\rho}_\ell \geq C_5 t_i^{1/2}/32\right)\leq c_3.
\end{equation}

For $i$ so that  $t_i < 4 (t_i/2 -1)$, by \eqref{step_58} and \eqref{step_44}, applying Abel's equation $\sum_{k=1}^q a_k b_k = \sum_{k=1}^{q-1}(a_k - a_{k+1})\sum_{\ell=1}^k b_\ell + a_q \sum_{\ell=1}^q b_\ell$,  we get
\begin{align}\label{step_37}
		& \sum_{\ell = 1}^{t_i - t_{i-1}-1}\exp\left\{ \varepsilon  \frac{t_i}{C_9\ell} \right\} \mathbb{P}\left(D(\ell)\right) \nonumber \\
		&= \sum_{\ell=1}^{t_i/2 -2}\left(\exp\left\{ \varepsilon  \frac{t_i}{C_9\ell} \right\} - \exp\left\{ \varepsilon  \frac{t_i}{C_9(\ell+1)} \right\}\right)\sum_{q=1}^\ell \mathbb{P}(D(q)) + \exp\left\{ \varepsilon  \frac{t_i}{C_9(t_i/2 -1)} \right\}\sum_{\ell = 1}^{t_i/2 -1}\mathbb{P}(D(\ell))\nonumber\\
		& \leq \sum_{\ell=1}^{t_i/2 -2}\left(\exp\left\{ \varepsilon  \frac{t_i}{C_9\ell} \right\} - \exp\left\{ \varepsilon  \frac{t_i}{C_9(\ell+1)} \right\}\right)c_1e^{-c_2t_i /\ell } + c_3\exp\left\{ \varepsilon  \frac{4}{C_9} \right\}.
\end{align}
Noting that $(e^{a/x})' = -a x^{-2}e^{a/x},$ we have
\begin{align*}
	&\left(\exp\left\{ \varepsilon  \frac{t_i}{C_9\ell} \right\} - \exp\left\{ \varepsilon  \frac{t_i}{C_9(\ell+1)} \right\}\right)c_1e^{-c_2t_i /\ell }\\
	& = c_1e^{-c_2t_i /\ell } \int_\ell^{\ell +1}\frac{\varepsilon t_i}{C_9 x^2}e^{\varepsilon \frac{t_i}{C_9 x}}\mathrm{d}x \leq c_1\int_\ell^{\ell +1}\frac{\varepsilon t_i}{C_9 x^2}e^{\varepsilon \frac{t_i}{C_9 x}}e^{-c_2t_i /x }\mathrm{d}x.
\end{align*}
So by \eqref{step_37}, we conclude that
\begin{align}\label{step_59}
	&\sum_{\ell = 1}^{t_i - t_{i-1}-1}\exp\left\{ \varepsilon  \frac{t_i}{C_9\ell} \right\} \mathbb{P}\left(D(\ell)\right) \leq c_1 \int_1^{t_i/2-1} \frac{\varepsilon t_i}{C_9 x^2}e^{\varepsilon \frac{t_i}{C_9 x}}e^{-c_2t_i /x }\mathrm{d} x + c_3\exp\left\{ \varepsilon  \frac{4}{C_9} \right\}\nonumber \\
	& \stackrel{z:= t_i/x}{=} c_{1}\frac{\varepsilon}{C_9}\int_{t_i/(t_i/2 -1)}^{t_i} e^{z(\varepsilon/C_9 - c_2)}\mathrm{d}z+
c_3 \exp\left\{ \varepsilon  \frac{4}{C_9} \right\}\nonumber \\
	&\leq c_1\frac{\varepsilon}{C_9}\int_{2}^{\infty} e^{z(\varepsilon/C_9 - c_2)}\mathrm{d}z+ c_3\exp\left\{ \varepsilon  \frac{4}{C_9} \right\} =: F(\varepsilon).
\end{align}
Since $F(0) < 1$,
we may take $i_1$ sufficiently large and $\varepsilon$ sufficiently small such that
$$
F(\varepsilon) \cdot \frac{1+\zeta\left(C_5t_{i_1}^{1/2}/(32M_\rho)
	\right)}{(C_9)^\varepsilon} < 1.
$$
Since $\zeta$ is decreasing, for sufficiently small $\varepsilon$ and $i\geq i_1$, using \eqref{domi-exp-G-Z} and \eqref{step_59}, we get
 \begin{align*}
&\mathbb{E}\left(\exp\left\{- \varepsilon
1_{G_i^1}
\ln Z_i^1 \right\}\right) \\& \leq \frac{e^{2\varepsilon /C_9}}{(C_9)^\varepsilon}\sum_{k=0}^\infty  \left(\frac{1+\zeta\left(C_5 t_i^{1/2}/(32M_\rho)
	\right)}{(C_9)^\varepsilon}\right)^{k} \left(\sum_{\ell = 1}^{t_i - t_{i-1}-1}\exp\left\{ \varepsilon  \frac{t_i}{C_9\ell} \right\} \mathbb{P}\left(D(\ell)\right)\right)^k\\
& \leq \frac{e^{2\varepsilon /C_9}}{(C_9)^\varepsilon}\sum_{k=0}^\infty  \left(\frac{1+\zeta\left(C_5 t_i^{1/2}/(32M_\rho)
	\right)}{(C_9)^\varepsilon} \cdot F(\varepsilon)\right)^{k}\\
 	& \leq \frac{e^{2\varepsilon /C_9}}{(C_9)^\varepsilon}\sum_{k=0}^\infty  \left(\frac{1+\zeta\left(C_5 t_{i_1}^{1/2}/(32M_\rho)
 		\right)}{(C_9)^\varepsilon} \cdot F(\varepsilon)\right)^{k}
 =: K< \infty.
 \end{align*}
\hfill$\Box$

Using similar arguments,
we have the following lemmas \ref{remark2}-\ref{lemma 4.19}.
Since the arguments are similar, we will only sketch the proofs.

In the following,  $\bar{C}$ is a positive constant which will be specified later in the proof of Lemma \ref{lemma 4.9}.
Define
$k(n):=\left[ \log_2 \left([\bar{C}\ln n]^2\right)\right]$ and
\begin{align*}
	r_0^1 & : = [\bar{C}\ln n],\qquad
L_n^1:=\left\lceil \frac{ C_5 \bar{C} \ln n}{4M_\rho } \right\rceil,\\
	s_0^1&:=t_{k(n)} \land \inf\left\{
k\geq r_0^1+1:
\sum_{\ell= r_0^1+1}^k \overline{\rho}_\ell^{(r_0^1-[L_n^1/2])}
\geq \frac{1}{4}C_5 \bar{C} \ln n \right\},\quad
	r_{j+1}^1 : = s_j^1 + L_n^1, \\
	s_{j+1}^1& : = \inf\left\{
k\geq r_{j+1}^1+1:
\sum_{\ell= r_{j+1}^1 +1}^k \overline{\rho}_\ell^{(r_{j+1}^1 -[L_n^1/2])}
 \geq\frac{1}{4} C_5 \bar{C} \ln n\right\}
\land \left(r_{j+1}^1 + t_{k(n)} - [\bar{C}\ln n]\right),
	\\ \mathfrak{K}^1 & : = \inf \{ k : s_k^1  \geq t_{k(n)}\}.
\end{align*}
Define
\begin{align*}
\widetilde{G}_n^1
&:= \left\{\max_{[\bar{C}\ln n] \leq r \leq 2t_{k(n)}, 0< k\leq  [\bar{C}\ln n] /\alpha_2}\ \sum_{m=r+1}^{r+k} \overline{\rho}_m^{\left(\left[2^{-1}\bar{C} \ln n
	\right]\right)}
\leq \frac{1}{8}C_5 \bar{C} \ln n \right\},
	\end{align*}
where $\alpha_2$ is defined in \eqref{def-alpha}.
Recall that $C_9$  is the constant  in Lemma \ref{lemma 4.15}.

\begin{lemma}\label{remark2}
(i)
When $n$ is large enough, on $\widetilde{G}_n^1$,
we have that
\begin{align}\label{step_61}
&\inf_{y\in[-1,1]}
 \inf_{x\geq C_5[\bar{C}\ln n]}\widetilde{\Pi}_y^{\xi, -\lambda^*} \Big\{\beta_k^1 -\beta_2^1 \geq 0,\ \forall k \mbox{ such that } [\bar{C}\ln n]+1\leq k\leq t_{k(n)},\nonumber\\
		&\hspace{3cm} \beta_{t_{k(n)}}^1 - \beta_2^1 \geq C_5 t_{k(n)}^{1/2}/2 \, \Big|\, \beta_{[\bar{C}\ln n]}^1-\beta_2^1 = x \Big\}\nonumber\\
		& \geq C_9^{\mathfrak{K}^1+1}\exp\left\{-\sum_{j=0}^{\mathfrak{K}^1 } \frac{t_{k(n)}}{C_9(s_j^1 - r_j^1)} \right\}=: Y_n^1.
	\end{align}

(ii)
There exist positive constants $\varepsilon$  and $K= K(\varepsilon)$
such that for $n$ large enough,
$$	
\mathbb{E}\left(\exp\left\{- \varepsilon 1_{\tilde{G}_n^1}\ln Y_n^1 \right\}\right) \leq K.
$$
\end{lemma}
\textbf{Proof: }
For $\ell < t_{k(n)} - [\bar{C}\ln n]$, define
$$
D^1(\ell):= \left\{\omega:\, \sum_{r= 1}^k \overline{\rho}_r^{(-[L_n^1/2])} <
\frac{1}{4}C_5 \bar{C} \ln n
,
\,\forall k\in[1, \ell),\,
\sum_{r= 1}^\ell \overline{\rho}_r^{(-[L_n^1/2])}
\geq \frac{1}{4} C_5 \bar{C} \ln n \right\}.
$$
Using the $\zeta$-mixing condition,  similar to \eqref{Mixing1-1}, we have
\begin{align}\label{Mixing2-1}
&\mathbb{P}\left(s_j^1 -r_j^1 = \ell_j, 0\leq j\leq k-1, s_k -r_k = \ell_k \right) \nonumber \\
&\leq \mathbb{P}\left(s_j^1 -r_j^1 = \ell_j, 0\leq j\leq k-1 \right)\mathbb{P}\left(D^1(\ell_k)\right)
\left(1+\zeta\left(C_5 \bar{C} \ln n/(8M_\rho)
\right)\right).
\end{align}
As variants of $\mathcal{A}_j$, $(\mathcal{A}_j)'$ and $\mathcal{B}$ in  Lemma \ref{lemma 4.15}, we define
\begin{align*}
	&\mathcal{A}_j^1  :=
	\left\{\beta_{s_j^1}^1 -\beta_2^1 \geq 4C_5 \bar{C} \ln n \right\}\bigcap \bigcap_{k=r_j^1 +1}^{s_j^1}
	 \left\{ \beta_k^1 -\beta_{r_j^1}^1 +\sum_{\ell = r_j^1+1}^k \rho_\ell
\geq -\frac{1}{4} C_5 \bar{C} \ln n  \right\}, \\
	&(\mathcal{A}_j^1 )'  := \bigcap_{k=s_j^1+1}^{r_{j+1}^1}\left\{\beta_k^1 -\beta_{s_j^1}^1 + \sum_{\ell = s_j^1+1}^k \rho_\ell
 \geq -\frac{1}{4} C_5 \bar{C} \ln n\right\},\quad
	\mathcal{B}^1  : = \left(\bigcap_{j=0}^{\mathfrak{K}^1} \mathcal{A}_j^1\right) \bigcap \left(\bigcap_{j=0}^{\mathfrak{K}^1-1} (\mathcal{A}_j^1 )'\right).
\end{align*}
Then repeating the arguments in  Step 1 and Step 2 in the proof of Lemma \ref{lemma 4.15},
also note that $C_5 [\bar{C}\ln n]/2 \geq C_5 t_{k(n)}^{1/2}/2$,
we get  \eqref{step_61} .

Define
$$\mathfrak{J}^1: = \inf\{j: s_j^1-r_j^1 \geq t_{k(n)} - [\bar{C}\ln n] \},$$
and
$$ b_n:= t_{k(n)}-[\bar{C}\ln n].$$
By repeating the argument of the proof of Lemma \ref{lemma 4.15} (ii), we get that
when $n$ is large enough so that
$t_{k(n)}\leq 2\left(t_{k(n)} - [\bar{C}\ln n] \right)$,
\begin{align*}
&\mathbb{E}\left(\exp\left\{- \varepsilon 1_{\tilde{G}_n^1}\ln Y_n^1 \right\}\right)\nonumber\\
	&\leq  \sum_{k=0}^\infty \left(\frac{1}{C_9 }\right)^{\varepsilon (k+1)}\exp\left\{\varepsilon  \frac{t_{k(n)}}{C_9(t_{k(n)} - [\bar{C}\ln n])}\right\}
	\left(1+\zeta\left(C_5 \bar{C} \ln n/(8M_\rho)
	\right)\right)^{k}\\
	& \qquad \times
	\left(\sum_{\ell = 1}^{t_{k(n)} - [\bar{C}\ln n]-1}\exp\left\{ \varepsilon  \frac{t_{k(n)}}{C_9\ell} \right\} \mathbb{P}\left(D^1(\ell)\right)\right)^k\\
	& \leq \sum_{k=0}^\infty \left(\frac{1}{C_9 }\right)^{\varepsilon (k+1)}e^{2\varepsilon / C_9}
	 \left(1+\zeta\left(C_5 \bar{C} \ln n/(8M_\rho)
	 \right)\right)^{k}
	\left(\sum_{\ell = 1}^{b_n-1}\exp\left\{ \varepsilon  \frac{t_{k(n)}}{C_9\ell} \right\} \mathbb{P}\left(D^1(\ell)\right)\right)^k.
\end{align*}
Since $b_n\sim t_{k(n)}$,  we have $\sqrt{b_n} \sim t_{k(n)}^{1/2} \leq \bar{C} \ln n  $.
Similar to  \eqref{step_58} and \eqref{step_44}, we have,  for $q< b_n$,
\begin{align*}
	\sum_{\ell = 1}^{q}\mathbb{P}(D^1(\ell)) &\leq c_1e^{-c_2t_{k(n)} /q},\\
	\sum_{\ell = 1}^{b_n -1}\mathbb{P}(D^1(\ell))  & \leq\mathbb{P}\left(
	\max_{1\leq k \leq b_n} \sum_{m=1}^k
	 \rho_m \geq
	\frac{1}{8} C_5 \bar{C} \ln n
	\right)\leq  c_3 < 1.
\end{align*}
The same arguments as in \eqref{step_37} and \eqref{step_59} show that for large $n$ such that $t_{k(n)}\leq 2(b_n-1)$,
\begin{align*}
	& \sum_{\ell = 1}^{b_n-1}\exp\left\{ \varepsilon  \frac{t_{k(n)}}{C_9\ell} \right\} \mathbb{P}\left(D^1(\ell)\right) \\ & \leq \sum_{\ell=1}^{b_n -2}\left(\exp\left\{ \varepsilon  \frac{t_{k(n)}}{C_9\ell} \right\} - \exp\left\{ \varepsilon  \frac{t_{k(n)}}{C_9(\ell+1)} \right\}\right)c_1e^{-c_2t_{k(n)} /\ell } + c_3 \exp\left\{ \varepsilon  \frac{t_{k(n)}}{C_9(b_n -1)} \right\}\\
	& \leq c_1\frac{\varepsilon}{C_9}\int_{1}^{\infty} e^{z(\varepsilon/C_9 - c_2)}\mathrm{d}z+ c_3\exp\left\{ \varepsilon  \frac{2}{C_9} \right\}=: F^1(\varepsilon).
\end{align*}
Thus, there exists $n_0>0$ such that when $n\geq n_0$, it holds that
\begin{align*}
\mathbb{E}\left(\exp\left\{- \varepsilon
1_{\tilde{G}_n^1}
\ln Y_n^1 \right\}\right)
	\leq\frac{ e^{2\varepsilon / C_9}}{(C_9)^\varepsilon} \sum_{k=0}^\infty \left(\frac{
	\left(1+\zeta\left(C_5 \bar{C} \ln n_0/(8M_\rho)
	\right)\right)
}{(C_9)^\varepsilon } F^1(\varepsilon)\right)^{k}.
\end{align*}
The rest of the proof is the same as the proof of (ii) in Lemma \ref{lemma 4.15}.
\hfill$\Box$

Recall that $t_{i}= 2^i$, and  $L_i$ is defined in \eqref{r-s-L-K}. For $n\in \mathbb{N}$,
define
\begin{align*}
		r_0^{(2,n)}& := t_{i-1},\\
		s_0^{(2,n)}& := t_i \land \inf\left\{ k\geq r_0^{(2,n)}+1:\ \sum_{\ell = r_0^{(2,n)}+1}^{k}\overline{\rho}_{n-\ell}^{(n-k-[L_i/2]))} \geq C_5 t_i^{1/2}/16 \right\}, \quad r_{j+1}^{(2,n)} := s_j^{(2,n)}  + L_i, \\
		s_{j+1}^{(2,n)} & : = \inf \left\{k\geq r_{j+1}^{(2,n)}+1:\ \sum_{\ell = r_{j+1}^{(2,n)}+1}^{k}\overline{\rho}_{n-\ell}^{(n-k-[L_i/2])} \geq C_5 t_i^{1/2}/16 \right\}\land (r_{j+1}^{(2,n)} + t_i - t_{i-1}),\\
\mathfrak{K}^{(2,n)} &:= \inf\{ k: s_k^{(2,n)} \geq t_i \}.
\end{align*}
Also define
\begin{align*}
G_i^{(2,n)}
:= \left\{\max_{r\in [t_{i-1}, t_{i+1}],
0< k\leq  t_i^{1/2}/\alpha_1 }\
\sum_{m=r+1}^{r+k} \overline{\rho}_{n-m}^{(n-t_{i+2})} \leq C_5t_{i-1}^{1/2}/32 \right\},
   \end{align*}
where $\alpha_1$ is defined in \eqref{def-alpha}.

When considering $I_2$, we need to be more careful
since the definition of $\beta_k^2$ only makes sense for $k\geq 1$. Note that our definition of $s_j^{(2,n)}$ and $r_j^{(2,n)}$ does not need any restriction for $n$. Thus, the proof of (ii)
for Lemma \ref{lemma 4.18} below
is the same as Lemma \ref{lemma 4.15}. Also note that our proof for (i) in Lemma \ref{lemma 4.15} only relies on
the property that the sequence of independent random variables $H_{\ell}^1 - H_{\ell -1}^1 - \widetilde{\Pi}_y^{\xi, -\lambda^*} \left[ H_{\ell}^1 - H_{\ell -1}^1\right], \ell \geq 2$ is a sequence of centered sub-exponential random variables satisfying \eqref{domi-MGF-above} and \eqref{domi-MGF-below}.
The strong Markov property of $\Xi^1$ only used for the independence of $H_\ell^1$.
Note that for $k\geq 2$, we have
$$\beta_k^2 -\beta_{k-1}^2 =-\rho_k +H_{k}^2 - H_{k -1}^2 - \widetilde{\Pi}_y^{\xi, -\lambda^*} \left[ H_{k}^2 - H_{k -1}^2\right] =:-\rho_k +\Delta H_k^2,$$
and for $|\eta|< \gamma(-\lambda^*)-(m-1)\textup{es}$,
\begin{align}\label{Extend-Delta-H}
	& \widetilde{\Pi}_y^{\xi, -\lambda^*} \left(e^{\eta \Delta H_k^2}\right)=\widetilde{\Pi}_{k-1}^{\xi, -\lambda^*} \left(\exp\left\{\eta \left(H_k^2 +  \frac{1}{v^*} \left(\frac{\psi_\lambda(k,-\lambda^*,\omega )}{\psi(k,-\lambda^*,\omega)} -\frac{\psi_\lambda(k-1,-\lambda^*,\omega)} {\psi(k-1,-\lambda^*,\omega)} \right) \right)\right\}\right)\nonumber \\
	& = \widetilde{\Pi}_{k-1}^{\xi, -\lambda^*} \left(\exp\left\{\eta \left(H_k^2 +  \frac{1}{v^*} \frac{\psi_\lambda(1,-\lambda^*,\theta_{k-1}\omega )}{\psi(1,-\lambda^*,\theta_{k-1}\omega )}\right)\right\}\right)\nonumber\\
	& = \exp\left\{ \frac{\eta}{v^*} \frac{\psi_\lambda(1,-\lambda^*,\theta_{k-1}\omega )}{\psi(1,-\lambda^*,\theta_{k-1}\omega )}\right\} \Pi_0\left(\exp\left\{\int_{0}^{H_1} \left((m-1)\xi(B_s,\theta_{k-1}\omega) -\gamma(-\lambda^*)+\eta\right)\mathrm{d}s \right\} \right)
\end{align}
and the right-hand side of \eqref{Extend-Delta-H} also makes sense for $k\leq 1$.
Therefore, to extend the case of $\beta_k^2$ for $k\leq 1$,
we ignore the definition of $\beta_0^2$ and $\beta_1^2$ defined in \eqref{def-beta^j},  and
we take a sequence of independent random variables $\left\{\Delta H_{k}^2: k\leq 1\right\}$ independent to everything else with Laplace function given in \eqref{Extend-Delta-H} under $\widetilde{\Pi}_y^{\xi, -\lambda^*}$,
then, or $k\leq 1$, define
$$\beta_k^2:= \beta_2^2 +\sum_{\ell=k+1}^2
 \left(\rho_\ell - \Delta H_\ell^2\right).$$
Replacing the definition of $\beta_0^2$ and $\beta_1^2$ defined in \eqref{def-beta^j}, by the above,  Lemma \ref{lemma 4.18}(i) can be proven by an argument similar to that of Lemma \ref{lemma 4.15}(i).

\begin{lemma}\label{lemma 4.18}
(i) For sufficiently large $i$ and
$n\geq 2^{i-1}$,
on $G_i^{(2,n)}$,
\begin{align} \label{step_62}
&\inf_{y\in[-1,1]}
	\inf_{x\geq C_5 t_{i-1}^{1/2}/2} 	\widetilde{\Pi}_y^{\xi, -\lambda^*}  \left(
\{\beta_{n -t_i}^2-\beta_n^2 \geq C_5 {t_i}^{1/2}/2\}\cap \cap^{t_i}_{k=t_{i-1}+1}\{\beta_{n-k}^2 -\beta_n^2 \geq t_i^{1/4}\}
\Big| \beta_{n-t_{i-1}}^2-\beta_n^2 =x\right)\nonumber\\
	& \geq  C_9^{\mathfrak{K}^{(2,n)}+1}
	\exp\left\{- \sum_{j=0}^{\mathfrak{K}^{(2,n)}}\frac{t_i}{C_9\left(s_j^{(2,n)} - r_j^{(2,n)}\right)} \right\} =: Z_i^{(2,n)},
\end{align}
here $C_9$ is the constant in Lemma \ref{lemma 4.15}.

(ii) There exist constants $\varepsilon > 0$  and
$K= K(\varepsilon)$ such that for $n$ large enough,
\begin{equation}
\mathbb{E}\left(\exp\left\{- \varepsilon
1_{G_i^{(2,n)}}
\ln Z_i^{(2,n)} \right\}\right)
\leq K.
	\end{equation}
\end{lemma}
\textbf{Proof: } (i) For $r<s$ define
\begin{align}
	D^{(2,n)}(r,s):= \left\{\omega:\, \sum_{\ell = r+1}^{k}\overline{\rho}_{n-\ell}^{(n-k-[L_i/2])} < C_5 t_i^{1/2}/16,
\forall k\in[r+1,s),
 \sum_{\ell = r+1}^{s}\overline{\rho}_{n-\ell}^{(n-s-[L_i/2])} \geq C_5 t_i^{1/2}/16 \right\}.
\end{align}
Then by \eqref{Measurable-of-rho} we see that $D^{(2,n)}(r,s)\in \mathcal{F}_{n-r}\cap \mathcal{F}^{n-s-[L_i/2]}$. Define
$$
D^{(2,n)}(\ell):= \left\{\omega:\, \sum_{\ell = 1}^{k}\overline{\rho}_{n-\ell}^{(n-k-[L_i/2])} < C_5 t_i^{1/2}/16,
\forall k\in[1, s),
\sum_{\ell = 1}^{s}\overline{\rho}_{n-\ell}^{(n-s-[L_i/2])} \geq C_5 t_i^{1/2}/16 \right\}.
$$
Similar to \eqref{Mixing1-1} and \eqref{Mixing2-1}, we have
\begin{align}\label{Mixing3-1}
	&\mathbb{P}\left(s_j^{(2,n)} -r_j^{(2,n)} = \ell_j, 0\leq j\leq k-1, s_k^{(2,n)} -r_k^{(2,n)} = \ell_k \right) \nonumber \\
&\leq 	\mathbb{P}\left(s_j^{(2,n)} -r_j^{(2,n)} = \ell_j, 0\leq j\leq k-1 \right)\mathbb{P}\left(D^{(2,n)}(\ell_k)\right)\left(1+ \zeta\left(C_5 t_i^{1/2}/(32M_\rho)
\right)\right).
\end{align}
Define
\begin{align*}	
	&\mathcal{A}_j^{(2,n)}  := \left\{
	\beta_{n-s_j^{(2,n)}}^2 -\beta_n^2 \geq  C_5 t_i^{1/2}\right\}
	\bigcap \bigcap^{s_j^{(2,n)} }_{k=r_j^{(2,n)}+1}
	\left\{\beta_{n-k}^2 -\beta_{n-r_j^{(2,n)}}^2 +\sum_{\ell = r_j^{(2,n)}+1}^k \rho_{n-\ell}  \geq - C_5 t_i^{1/2}/16\right\}, \\
	&\left(\mathcal{A}_j^{(2,n)}\right)' :=\bigcap_{k=s_j^{(2,n)}+1}^{r_{j+1}^{(2,n)}} \left\{ \beta_{n-k}^2 -\beta_{n-s_j^{(2,n)}}^2 + \sum_{\ell = s_j^{(2,n)}+1}^k \rho_{n-\ell}\geq - C_5 t_i^{1/2}/16\right\},\\
	&\mathcal{B}^{(2,n)}  : = \left(\bigcap_{j=0}^{\mathfrak{K}^{(2,n)}} \mathcal{A}_j^{(2,n)} \right)\bigcap \left(\bigcap_{j=0}^{\mathfrak{K}^{(2,n)}-1} \left(\mathcal{A}_j^{(2,n)}\right)'\right).
\end{align*}
Then repeating the proof of Lemma \ref{lemma 4.15} (i), we get \eqref{step_62}.

(ii) Define
$$
\mathfrak{J}^{(2,n)}: = \inf\{j: s_j^{(2,n)}-r_j^{(2,n)} \geq t_i - t_{i-1} \}.
$$
Although $s_j^{(2,n)}-r_j^{(2,n)}$ may depend on $n$, we see that the upper bound of \eqref{step_58} and \eqref{step_44} still hold if we replace $D(\ell_k)$ by $D^{(2,n)}(\ell_k)$ because of the stationarity of $\rho$.
Repeating the argument in the proof of Lemma \ref{lemma 4.15} (ii), we get that (ii) holds.
The proof is now complete.
\hfill$\Box$

Let
$L_n^1$ be defined in Lemma \ref{remark2}.
Define
\begin{align*}
		r_0^{(2)} & : = [\bar{C}\ln n],\\
	s_0^{(2)} &:= t_{k(n)} \land \inf\left\{
k\geq r_0^{(2)}+1:
\sum_{\ell= r_0^{(2)}+1}^k \overline{\rho}_{n-\ell}^{(n-k- [L_n^1/2])}
\geq\frac{1}{4}C_5 \bar{C} \ln n \right\},\quad
	r_{j+1}^{(2)} : = s_j^{(2)} + L_n^1 \\
	s_{j+1}^{(2)}& : = \inf\left\{
k\geq r_{j+1}^{(2)}+1:\sum_{\ell= r_{j+1}^{(2)}+1}^k
\overline{\rho}_{n-\ell}^{(n-k-[L_n^1/2])}
\geq\frac{1}{4} C_5 \bar{C} \ln n\right\}
\land \left(r_{j+1}^{(2)} + t_{k(n)} - [\bar{C}\ln n]\right),
	\\ \mathfrak{K}^{(2)} & : = \inf \{ k : s_k^{(2)}  \geq t_{k(n)} \}.
\end{align*}
Define
\begin{align*}
\tilde{G}_n^{2} :=
\left\{\max_{[\bar{C}\ln n] \leq r \leq 2t_{k(n)},
0< k\leq  [\bar{C}\ln n]
 /\alpha_2}\ \sum_{m=r+1}^{r+k} \overline{\rho}_{n-m}^{(n-4t_{k(n)})}
 \leq \frac{1}{8}C_5 \bar{C} \ln n \right\},
\end{align*}
here $\alpha_2$ is defined in \eqref{def-alpha}.

\begin{lemma}\label{lemma 4.19}
(i) For sufficiently large $n$,  on $\tilde{G}_n^{2}$,
\begin{align*}
&\inf_{y\in[-1,1]}
	\inf_{x\geq C_5[\bar{C}\ln n]} \widetilde{\Pi}_y^{\xi, -\lambda^*}
\left\{ \{\beta_{n -t_{k(n)}}^2 - \beta_n^2 \geq C_5 t_{k(n)}^{1/2}/2\}\cap \cap^{t_{k(n)}}_{k=[\bar{C}\ln n]+1}\{\beta_{n-k}^2 -\beta_n^2 \geq 0\}
\mid \beta_{n-[\bar{C}\ln n]}^2-\beta_n^2 =x
\right\}\\&\geq   C_9^{\mathfrak{K}^{(2)}+1}\exp\left\{-\sum_{j=0}^{\mathfrak{K}^{(2)} } \frac{t_{k(n)}}{C_9\left(s_j^{(2)} - r_j^{(2)}\right)} \right\}=: Y_n^2,
\end{align*}
here $C_9$ is the constant in Lemma \ref{lemma 4.15}.

(ii) There exist
positive constants $\varepsilon$  and $K= K(\varepsilon)$
such that for $n$ large enough,
$$	
\mathbb{E}\left(\exp\left\{- \varepsilon 1_{\tilde{G}_n^{2}}\ln Y_n^{2} \right\}\right) \leq K.
$$
\end{lemma}
\textbf{Proof: } (i) Define
\begin{align*}
	&\mathcal{A}_j^{(2)}  := \left\{
	\beta_{n-s_j^{(2)}}^2 -\beta_n^2 \geq 4C_5 \bar{C} \ln n \right \}
	\bigcap \bigcap^{s_j^{(2)}}_{k=r_j^{(2)}+1}
	\left\{\beta_{n-k}^2 -\beta_{n-r_j^{(2)}}^2 +\sum_{\ell = r_j^{(2)}+1}^k \rho_{n-\ell} \geq
-\frac{1}{4} C_5 \bar{C} \ln n\right\}, \\
	&\left(\mathcal{A}_j^{(2)}\right)'  :=\bigcap_{k=s_j^{(2)}+1}^{r_{j+1}^{(2)}} \left\{
	\beta_{n-k}^2 -\beta_{n-s_j^{(2)}}^2 + \sum_{\ell = s_j^{(2)}+1}^k \rho_{n-\ell}
\geq -\frac{1}{4}C_5 \bar{C} \ln n\right\},\\
	&\mathcal{B}^{(2)}  : = \left(\bigcap_{j=0}^{\mathfrak{K}^{(2)}} \mathcal{A}_j^{(2)} \right)\bigcap \left(\bigcap_{j=0}^{\mathfrak{K}^{(2)}-1} \left(\mathcal{A}_j^{(2)}\right)'\right).
\end{align*}
Repeating the proof of in Lemma \ref{lemma 4.15} (i), we get that (i) holds.

(ii)
Define
$$\mathfrak{J}^{(2)}: = \inf\{j: s_j^{(2)}-r_j^{(2)} \geq t_{k(n)} - [\bar{C}\ln n] \}.$$
Using an argument similar to that in the proof  of Lemma \ref{lemma 4.15} (ii), we get that (ii)  holds.
\hfill$\Box$

\textbf{Proof of Lemma \ref{lemma 4.9}:}
Recall the definitions of $I_1$ and $I_2$  in \eqref{def-I}.
Recall that
$k(n)=\left[ \log_2 \left([\bar{C}\ln n]^2\right)\right]$.
Define
 $j(n) : = \lceil \log_2 n \rceil$.
Let $\bar C$ be a constant to be specified later.
By \eqref{step_36}, there exists a large $N_1$ such that for $i,n \geq N_1$ and $\ell \geq t_{i+2}= 2^{i+2}$, we have that
\begin{align*}
 \left\{\max_{r\in [t_{i-1}, t_{i+1}], 0< k\leq  t_i^{1/2} /\alpha_1}\ \sum_{m=r+1}^{r+k} {\rho}_m \leq C_5t_{i-1}^{1/2}/64 \right\} & \subset
G_i^1 ,\\
 \left\{\max_{[\bar{C}\ln n] \leq r \leq 2[\bar{C}\ln n]^2, 0< k\leq  [\bar{C}\ln n] /\alpha_2}\ \sum_{m=r+1}^{r+k} {\rho}_m \leq
 \frac{1}{16} C_5 \bar{C} \ln n
  \right\} & \subset
\widetilde{G}_n^1,\\
 \left\{\max_{r\in [t_{i-1}, t_{i+1}], 0<k\leq  t_i^{1/2}/\alpha_1 }\ \sum_{m=r+1}^{r+k} {\rho}_{\ell-m}\leq C_5t_{i-1}^{1/2}/64 \right\} & \subset
G_i^{(2,\ell)},\\
 \left\{\max_{[\bar{C}\ln n] \leq r \leq 2[\bar{C}\ln n]^2, 0< k\leq  [\bar{C}\ln n] /\alpha_2}\ \sum_{m=r+1}^{r+k}
 {\rho}_{n-m} \leq
 \frac{1}{16} C_5 \bar{C} \ln n
 \right\} & \subset
\widetilde{G}_n^2.
\end{align*}
Take  an integer $N_2$ such that $N_2 \geq 2^{N_1 +2}$.
By Lemma \ref{lemma 4.12}, \eqref{step_35} and the stationarity of $\rho_m$,
we have
\begin{align}\label{sum-1}
& \sum_{n=N_2}^\infty
\sum_{i=k(n)}^{j(n)}
\mathbb{P}
	\left((G_i^1)^c\right)
+ \sum_{\ell=N_2}^\infty
\sum_{i=k(\ell)}^{j(\ell)}
\mathbb{P}\left((G_i^{(2,\ell)})^c\right)\nonumber\\
&\leq 2\sum_{n=N_2}^\infty\sum_{i=k(n)}^{j(n)} t_{i+1} \mathbb{P}\left(\max_{ 0\leq k\leq  t_i^{1/2}/\alpha_1 }\ \sum_{m=1}^{k} \rho_m > C_5t_{i-1}^{1/2}/64 \right) \nonumber\\
	& \leq 2\sum_{n=N_2}^\infty 2n \sum_{i=k(n)}^{j(n)} c_1 \exp\left\{-\alpha_1 C_5^2 t_{i-1}/(64^2 c_1t_i^{1/2})\right\} \nonumber\\
&\leq  4c_1\sum_{n=N_2}^\infty n\log_2(n) \exp\left\{-\alpha_1 C_5^2 \bar{C} \ln n/(64^2 c_1)\right\} < \infty,
\end{align}
and
\begin{align}\label{sum-2}
&\sum_{n=N_1}^\infty \mathbb{P}
\left(\left(\tilde{G}_n^1\right)^c\right)
+ 	\sum_{n=N_1}^\infty \mathbb{P}
\left(\left(\tilde{G}_n^{2}\right)^c\right)
\leq 2\sum_{n=N_1}^\infty 2[\bar{C}\ln n]^2 \mathbb{P}\left(\max_{0\leq k\leq  [\bar{C}\ln n]/\alpha_2 }\ \sum_{m=1}^{k} \rho_m>	
\frac{1}{16} C_5 \bar{C} \ln n
\right) \nonumber\\
	& \leq 4\sum_{n=N_1}^\infty [\bar{C}\ln n]^2 \cdot c_1 \exp\left\{-\alpha_2 C_5^2[\bar{C} \ln n]/(16^2 c_1)\right\} < \infty,
\end{align}
by taking $\bar{C}$ large enough.

Recall the definition of $Z_i^1$ and
$G_i^1$
in Lemma \ref{lemma 4.15}.
By \eqref{Measurable-5}, \eqref{Measurable-of-rho}, we have that for all $i\in\mathbb{Z}_+$ large enough such that $t_{i-1}> L_i$, $t_i > t_i^{1/2}/\alpha_1
$ and $t_{i-1}- [L_i/2]> t_{i-2}$,
\begin{align}\label{Measurable-2}
&	G_i^1
\in \mathcal{F}^{t_{i-2}} \cap \mathcal{F}_{t_{i+1}+ t_i^{1/2}/\alpha_1
} \subset \mathcal{F}^{t_{i-2}} \cap \mathcal{F}_{t_{i+1}+ t_i},\nonumber \\ &Z_i^1 \in\mathcal{F}^{t_{i-1}-[L_i/2]}\cap \mathcal{F}_{t_{i+1}}\subset \mathcal{F}^{t_{i-2}} \cap \mathcal{F}_{t_{i+1}+ t_i}.
	\end{align}
Then when $i$ is large enough, $ Z_{4i}^1 ,
G_{4i}^1
\in \mathcal{F}_{t_{4i+1}+ t_{4i}}\cap \mathcal{F}^{t_{4i -2}}$.
Together with Lemma \ref{lemma 4.15}, for large $i$, we have
\begin{align*}
&\mathbb{E}\left(\exp\left\{-\varepsilon
1_{G_{4i}^1}
\ln Z_{4i}^1 \right\}\Big| \mathcal{F}_{t_{4i-3	}+t_{4i-4}}\right)\leq \left(1+\zeta(2^{4i-4} )\right)\mathbb{E}\left(\exp\left\{-\varepsilon
1_{G_{4i}^1}
\ln Z_{4i}^1 \right\}\right)
 \leq K \left(1+\zeta(0)\right).
\end{align*}
Thus, we have, for $k=0,1,2,3$,
\begin{align*}
\mathbb{E}\left(\exp\left\{-\varepsilon \sum_{i= \lceil k(n)/4\rceil}^{[j(n) /4]-1}
1_{G_{4i}^1}
\ln Z_{4i+k} \right\}\right) & \leq \prod_{i=\lceil k(n)/4\rceil}^{[j(n)/4]-1}
K\left(1+\zeta(0)\right)\leq
c_1^{j(n)}
\end{align*}
for some constant $c_1>0.$ Now for a constant $c_2 > 0$, we have
\begin{align*}
&\mathbb{P} \left(\sum_{i=k(n)+1}^{j(n)}1_{G_i^1}
\ln Z_i^1 < -c_2 j(n) \right)  \leq\sum_{k=0}^3 \mathbb{P}\left( \sum_{i= \lceil k(n)/4\rceil}^{[j(n) /4]-1}
1_{G_{4i+k}^1}\ln Z_{4i+k} < -c_2 j(n) / 4    \right) \\
& \leq 4 e^{-\varepsilon  c_2 j(n)/4}\max_{k=1,2,3,4} \mathbb{E}\left(\exp\left\{-\varepsilon \sum_{i= \lceil k(n)/4\rceil}^{[j(n) /4]-1}
1_{G_{4i+k}^1}\ln Z_{4i+k} \right\}\right)
\leq  4e^{-\varepsilon  c_2 j(n)/4}\cdot c_1^{j(n)}.
\end{align*}
Letting $ c_2 $ be large enough so that $-\varepsilon c_2 /4 + \ln  c_1 <-\ln 2$, there exists  $N_3$ such that
\begin{equation}\label{sum-3}
\sum_{n= N_3}^\infty \mathbb{P} \left(\sum_{i=k(n)+1}^{j(n)}1_{G_i^1}
\ln Z_i^1 < - c_2 j(n) \right) \leq \sum_{n= N_3}^\infty 4e^{-\varepsilon  c_2 j(n)/4}\cdot c_1^{j(n)} < \infty.
\end{equation}
Similarly, by Lemma \ref{remark2}, Lemma \ref{lemma 4.18} and Lemma \ref{lemma 4.19}, we have
\begin{align}\label{sum-4}
&\sum_{n= N_3}^\infty \mathbb{P}\left(1_{\tilde{G}_n^1}
\ln Y_n^1 < - c_2  \ln n \right) + \sum_{n= N_3}^\infty\P\left(\sum_{i=k(n)+1}^{j(n)}
1_{G_i^{(2,n)}}
\ln Z_i^{(2,n)}  < - c_2 j(n)\right) \nonumber\\
&\quad+ \sum_{n= N_3}^\infty\P \left(
1_{\tilde{G}_n^2}
\ln Y_n^2< - c_2 \ln n \right)< \infty.
\end{align}

 Now define
\begin{align*}
\widetilde{\Omega}:=&\Omega_2 \bigcap
\left\{(G_n^1 )^c\
\textup{i.o.} \right\}^c \bigcap \left\{
(\widetilde{G}_n^1)^c \
\textup{i.o.} \right\}^c \bigcap \left\{\left(\bigcap_{i=k(n)+1}^{j(n)}
G_i^{(2,n)}
\right)^c\ \textup{i.o.} \right\}^c \bigcap \left\{(
\widetilde{G}_n^2)^c \
\textup{i.o.} \right\}^c\\
&\bigcap
\left\{ \left\{\sum_{i=k(n)+1}^{j(n)}
1_{G_i^1}
\ln Z_i^1 < - c_2 j(n) \right\} \ \textup{i.o.}\right\}^c \bigcap
\left\{ \left\{
1_{\widetilde{G}_n^1}
\ln Y_n^1 < - c_2  \ln n \right\} \ \textup{i.o.}\right\}^c\\
&\bigcap
\left\{ \left\{\sum_{i=k(n)+1}^{j(n)}
1_{G_i^{(2,n)}}
\ln Z_i^{(2,n)}  < - c_2 j(n)\right\} \ \textup{i.o.}\right\}^c \bigcap
\left\{ \left\{
1_{\widetilde{G}_n^2}
 \ln Y_n^2< - c_2 \ln n\right\} \ \textup{i.o.}\right\}^c.
\end{align*}
Then by \eqref{sum-1}, \eqref{sum-2}, \eqref{sum-3} and \eqref{sum-4},
$\mathbb{P}(\widetilde{\Omega}) = 1$.
 By  the construction of $\widetilde{\Omega}$, for every $\omega \in \widetilde{\Omega}$, there exists $N= N(\omega)$ such that when $n\geq N$, we have
 $$
1_{G_n^1}= 1_{ \tilde{G}_n^1} = \prod_{i=k(n)+1}^{j(n)} 1_{ G_i^{(2,n)}} = 1_{\tilde{G}_n^2 }=1,
$$
 and
\begin{align*}
\sum_{i=k(n)+1}^{j(n)}1_{G_i^1}
 \ln Z_i^1 \geq - c_2 j(n), \quad 1_{\widetilde{G}_n^1}\ln Y_n^1 \geq - c_2  \ln n ,\\
\sum_{i=k(n)+1}^{j(n)}
1_{G_i^{(2,n)}}
\ln Z_i^{(2,n)}  \geq - c_2 j(n),\quad
1_{\widetilde{G}_n^2}
\ln Y_n^2\geq  - c_2 \ln n.
\end{align*}
Now for $\omega \in \widetilde{\Omega}$ and $n\geq N$, it holds that uniformly for $y\in[-1,1]$,
\begin{align*}
& \widetilde{\Pi}_y^{\xi, -\lambda^*} \left(\beta_{k}^1 - \beta_1^1 \geq 0,\ \forall 1\leq k \leq n, \beta_{n}^1 - \beta_1^1 \geq n^{1/4}\right)\\
&\geq  \prod_{k=1}^{[\bar{C}\ln n]}\widetilde{\Pi}_y^{\xi, -\lambda^*}\left(\beta_k^1 -\beta_{k-1}^1 > C_5\right)\cdot Y_n^1  \cdot
\prod_{i=k(n) + 1}^{j(n)} Z_i^1\\
&\geq  C_5^{[\bar{C}\ln n]}\cdot \exp\left\{
1_{\widetilde{G}_n^1}
\ln Y_n^1\right\}\cdot \exp\left\{\sum_{i=k(n)+1}^{j(n)}
1_{G_i^1}
 \ln Z_i^1 \right\}\\
& \geq C_5^{[\bar{C}\ln n]} \cdot e^{- c_2 \ln n} \cdot e^{- c_2 j(n)}\geq \frac{1}{C_5} n^{\bar{C}\ln C_5 - c_2 - c_2 /\ln 2},
\end{align*}
where in the second inequality we used Lemma \ref{lemma 4.12}.
Taking  $[V_t]= n$, and $\ell'' > -\bar{C} \ln C_5 +  c_2  + c_2 /\ln 2$,
we get the first result.

The second result can be proved similarly. Note that $n \geq 2^{j(n)-1}= t_{j(n)-1}$ holds for every $n>2$. When $n\geq N(\omega)$, by Lemma \ref{lemma 4.18}, uniformly for $y\in [-1,1]$,
\begin{align*}
	& \widetilde{\Pi}_y^{\xi, -\lambda^*} \left(\beta_{n-k}^2 - \beta_n^2 \geq 0,\ \forall 1\leq k \leq n, \beta_2^2 - \beta_n^2 \geq n^{1/4}\right)\\
&\geq  \prod_{k=1}^{[\bar{C}\ln n]}\widetilde{\Pi}_y^{\xi, -\lambda^*}\left(\beta_{n-k}^2 -
\beta_{n-k+1}^2
> C_5\right)\cdot
Y_n^2  \cdot  \prod_{i=k(n) + 1}^{j(n)} Z_i^{(2,n)}\\
& \geq  C_5^{[\bar{C}\ln n]}\cdot \exp\left\{1_{\widetilde{G}_n^2}
\ln Y_n^2\right\}\cdot \exp\left\{\sum_{i=k(n)+1}^{j(n)}1_{G_i^{(2,n)}}
\ln Z_i^{(2,n)} \right\}\\
& \geq \frac{1}{C_5} n^{\bar{C}\ln C_5 - c_2 - c_2 /\ln 2}.
\end{align*}
Taking $[V_t]= n$ and $\ell'' > -\bar{C} \ln C_5 +  c_2  + c_2 /\ln 2$ completes the proof.
\hfill$\Box$

\section{Appendix}

\begin{lemma}\label{indep-incre}
For $j=1,2$, let $\beta_k^{j}$ be define by \eqref{def-beta^j}. Under $\widetilde{\Pi}_y^{\xi,-\lambda^*}$, $\{\beta_k^{j}, k\geq
1 \}$ has independent increments.
\end{lemma}
\textbf{Proof: }
We only need to prove the result for $j=1$.
For any $n\geq 2$, $
1
\leq k_0 < k_1 < ...< k_n$, we only need to prove that $H_{k_1 }^1 - H_{k_0 }^1 ,...,H_{k_n }^1 - H_{k_{n-1}}^1$ are independent.
For any bounded measurable functions $f_1,...,f_n$, by the strong Markov property of  $\left\{\Xi_t, t\geq 0;\widetilde{\Pi}_y^{\xi,-\lambda^*}\right\}$, we have
\begin{align*}
	&\widetilde{\Pi}_y^{\xi, -\lambda^*} \left(\prod_{j=1}^nf_j\left(H_{k_j }^1 - H_{k_{j-1}}^1\right)\bigg| \mathcal{F}_{H_{k_{0}}^1}\right)= \widetilde{\Pi}_{k_0}^{\xi, -\lambda^*} \left(f_1\left(H_{k_1 }^1 \right)
\prod_{j=2}^nf_j\left(H_{k_j }^1 - H_{k_{j-1}}^1\right)\right)\\
	& = \widetilde{\Pi}_{k_0}^{\xi, -\lambda^*} \left(f_1\left(H_{k_1 }^1 \right)
	\widetilde{\Pi}_{k_0}^{\xi, -\lambda^*}\left(\prod_{j=2}^nf_j\left(H_{k_j }^1 - H_{k_{j-1}}^1\right)\bigg| \mathcal{F}_{H_{k_1}^1}\right)\right)\\
	& = \widetilde{\Pi}_{k_0}^{\xi, -\lambda^*} \left(f_1\left(H_{k_1}^1 \right)
	\widetilde{\Pi}_{k_1}^{\xi, -\lambda^*}\left(f_2\left(H_{k_2}^1\right)\prod_{j=3}^nf_j\left(H_{k_j}^1 - H_{k_{j-1}}^1\right)\right)\right)\\
	& = \widetilde{\Pi}_{k_0}^{\xi, -\lambda^*} \left(f_1\left(H_{k_1 }^1 \right)
	\right)\widetilde{\Pi}_{k_1}^{\xi, -\lambda^*}\left(f_2\left(H_{k_2}^1\right)\prod_{j=3}^nf_j\left(H_{k_j}^1 - H_{k_{j-1}}^1\right)\right)\\
	& = \widetilde{\Pi}_{k_0}^{\xi, -\lambda^*} \left(f_1\left(H_{k_1 }^1 \right)
	\right)\widetilde{\Pi}_{k_1}^{\xi, -\lambda^*}\left(f_2\left(H_{k_2}^1\right)\widetilde{\Pi}_{k_1}^{\xi, -\lambda^*}\left(\prod_{j=3}^nf_j\left(H_{k_j}^1 - H_{k_{j-1}}^1\right)\bigg| \mathcal{F}_{H_{k_2}^1}\right)\right)\\
	& =  \widetilde{\Pi}_{k_0}^{\xi, -\lambda^*} \left(f_1\left(H_{k_1 }^1 \right)
	\right)\widetilde{\Pi}_{k_1}^{\xi, -\lambda^*}\left(f_2\left(H_{k_2}^1\right)\right)\widetilde{\Pi}_{k_2}^{\xi, -\lambda^*}\left(f_3 \left(H_{k_3 }^1\right)\prod_{j=4}^nf_j\left(H_{k_j}^1 - H_{k_{j-1}}^1\right)\right)\\
	& =\cdots = \prod_{j=1}^n \widetilde{\Pi}_{k_{j-1}}^{\xi, -\lambda^*} \left(f_j\left(H_{k_j}^1 \right)
	\right).
\end{align*}
Taking expectation in the display above, we get
\begin{align}\label{strong-Markov-2}
	\widetilde{\Pi}_y^{\xi, -\lambda^*} \left(\prod_{j=1}^nf_j\left(H_{k_j }^1 - H_{k_{j-1}}^1\right)\right)=\prod_{j=1}^n \widetilde{\Pi}_{k_{j-1}}^{\xi, -\lambda^*} \left(f_j\left(H_{k_j }^1 \right)
	\right).
\end{align}
In particular, taking $f_1=...=f_{n-1}=1$, we get
$$
 \widetilde{\Pi}_y^{\xi, -\lambda^*} \left(f_n\left(H_{k_n }^1 - H_{k_{n-1}}^1\right)\right)= \widetilde{\Pi}_{k_{n-1}}^{\xi, -\lambda^*} \left(f_n\left(H_{k_n }^1 \right)
\right).
$$
Thus, \eqref{strong-Markov-2} can be rewritten as
\begin{align}\label{strong-Markov-3}
	\widetilde{\Pi}_y^{\xi, -\lambda^*} \left(\prod_{j=1}^nf_j\left(H_{k_j }^1 - H_{k_{j-1}}^1\right)\right)=\prod_{j=1}^n \widetilde{\Pi}_{y}^{\xi, -\lambda^*} \left(f_j\left(H_{k_j }^1 - H_{k_{j-1}}^1\right)
	\right),
\end{align}
which says  that  $\left(H_{k_1 }^1 - H_{k_0 }^1 \right),...,\left( H_{k_n }^1 - H_{k_{n-1}}^1\right)$ are independent.
\hfill$\Box$

\noindent

\begin{singlespace}

\end{singlespace}

\vskip 0.2truein

\noindent{\bf Haojie Hou:} School of Mathematical Sciences,
Peking
University,  Beijing, 100871, P.R. China. Email: {\texttt
houhaojie@pku.edu.cn}

\smallskip

\noindent{\bf Yan-Xia Ren:} LMAM School of Mathematical Sciences \& Center for
Statistical Science, Peking
University,  Beijing, 100871, P.R. China. Email: {\texttt
yxren@math.pku.edu.cn}

\smallskip
\noindent {\bf Renming Song:} Department of Mathematics,
University of Illinois,
Urbana, IL 61801, U.S.A.
Email: {\texttt rsong@math.uiuc.edu}

\end{document}